\documentclass[11pt,final]{article}
\usepackage{authblk}
\usepackage[letterpaper,margin=1in]{geometry}
\usepackage[utf8]{inputenc}
\usepackage[english]{babel}

\usepackage{epic,eepic}
\usepackage{amssymb,amsthm,amscd,latexsym,mathrsfs,units,enumerate,bm,bbm,cancel,physics,mathtools,braket}
\usepackage{enumitem}
\usepackage{color,xcolor}
\usepackage{graphicx,wrapfig}

\usepackage[colorlinks]{hyperref}
\usepackage[notref, notcite]{showkeys}

\renewcommand{\geq}{\geqslant}
\renewcommand{\leq}{\leqslant}

\theoremstyle{plain}
\newtheorem{theorem}{Theorem}[section]
\newtheorem{remark}[theorem]{Remark}
\newtheorem{lemma}[theorem]{Lemma}
\newtheorem{proposition}[theorem]{Proposition}
\newtheorem{corollary}[theorem]{Corollary}
\newtheorem*{thm*}{Theorem}

\theoremstyle{definition} 

\newtheorem{defn}[theorem]{Definition}

\theoremstyle{remark}
\newtheorem*{remark*}{Remark}

\numberwithin{equation}{section}

\newcommand{\N}{\mathbb{N}}

\newcommand{\R}{\mathbb{R}}

\DeclareMathOperator{\Span}{span}

\newcommand{\nulow}{0.2533}
\newcommand{\nuup}{3.9}
\newcommand{\ourint}{[\nulow,\nuup]}

\newcommand{\eq}[2]{\begin{equation}\begin{split} #1\end{split}#2\end{equation}}

\newcommand{\mat}[1]{\begin{pmatrix}#1\end{pmatrix}}

\title{Asymptotic stability of sharp fronts: \protect\\ 
Analysis and rigorous computation}

\author[1]{Blake~Barker\thanks{E-mail:blake@mathematics.byu.edu}}
\author[2]{Jared~C.~Bronski\thanks{E-mail:bronski@illinois.edu}}
\author[2]{Vera~Mikyoung~Hur\thanks{E-mail:verahur@illinois.edu}}
\author[3]{Zhao~Yang\thanks{E-mail:yangzhao@amss.ac.cn}}
\affil[1]{Department of Mathematics, Brigham Young University, Provo, UT 84604, USA}
\affil[2]{Department of Mathematics, University of Illinois at Urbana-Champaign, \protect\\
Urbana, IL 61801, USA}
\affil[3]{State Key Laboratory of Mathematical Sciences, \protect\\ Academy of Mathematics and Systems Science, Chinese Academy of Sciences,\protect\\ 
 Beijing, 100190, China}

\begin{document}

\maketitle

\begin{abstract}
We investigate the stability of traveling front solutions to nonlinear diffusive-dispersive equations of Burgers type, with a primary focus on the Korteweg-de Vries--Burgers (KdVB) equation, although our analytical findings extend more broadly. Manipulating the temporal modulation of the translation parameter of the front and employing the energy method, we establish asymptotic, nonlinear, and orbital stability, provided that an auxiliary Schr\"odinger equation possesses precisely one bound state. Notably, our result is independent of the monotonicity of the profile and does not necessitate the initial condition to be close to the front. We identify a sufficient condition for stability based on a functional that characterizes the `width' of the traveling wave profile. Analytical verification for the KdVB equation confirms that this sufficient condition holds for the relative dispersion parameter within an open interval $\supset [-0.25, 0.25]$, encompassing all monotone profiles. Utilizing validated numerics or rigorous computation, we present a computer-assisted proof demonstrating that the stability condition itself holds for parameter values within the interval $\ourint$.
\end{abstract}

\tableofcontents

\section{Introduction}\label{sec:intro}

We consider the nonlinear diffusive-dispersive equation of Burgers type
\begin{equation}\label{eqn:main}
u_t+uu_x=u_{xx}+{\mathcal L}u,
\end{equation}
subject to the boundary conditions
\begin{equation}\label{eqn:BC}
u(x,t)\rightarrow u_-~\text{as}~x\rightarrow-\infty\quad\text{and}\quad
u(x,t)\rightarrow u_+~\text{as}~x\rightarrow+\infty,
\end{equation}
where $x\in\mathbb{R}$, $t\geq0$, $u(x,t)$ is a real-valued function, and $u_->u_+$ are constants; ${\mathcal L}$ is defined through a Fourier multiplier as 
\begin{equation}\label{def:L}
\widehat{\mathcal{L} v} =\widehat{\mathcal{L}}(k) \hat v(k),\quad k\in\mathbb{R},
\end{equation}
with Fourier symbol $\widehat{\mathcal{L}}(k)$. We assume that  
\begin{equation}\label{cond:L}
\operatorname{Re}(\widehat{\mathcal{L}}(k))= 0
\quad\text{and}\quad
\widehat{\mathcal{L}}(0)=0.
\end{equation}
The former ensures that $\int_{-\infty}^{\infty} v\mathcal{L}v~dx=0$, and the latter guarantees that the left kernel of ${\mathcal L}$ includes constant functions. It should be noted that some physically relevant but mathematically irrelevant constants have been scaled out from \eqref{eqn:main}. 

Our primary example involves ${\mathcal L} = \nu \partial_x^3$ for some $\nu\in\mathbb{R}$, whereby \eqref{eqn:main} becomes the Korteweg-de Vries--Burgers (KdVB) equation
\begin{equation}\label{eqn:KdVB}
u_t + u u_x = u_{xx} + \nu u_{xxx},
\end{equation}
although our analytical findings extend beyond to encompass a wide range of dispersive operators. Whitham \cite{Whitham} and others (see, for instance, \cite{johnson1972shallow,su1969korteweg}) introduced \eqref{eqn:KdVB} as a fundamental mathematical model for undular bores propagating upstream in rivers, such as the Severn River in England and the Dordogne River in France. Applications in plasmas and various physical systems have been explored in \cite{grad1967unified,NS1970} among others. Notably, \eqref{eqn:KdVB} is among the simplest PDEs that incorporate nonlinearity, dissipation, and dispersion.

For the sake of exposition, we assume that the initial value problem for \eqref{eqn:main} and \eqref{eqn:BC} is well-posed for all time in a suitable function space, where $\mathcal{L}$ is defined in \eqref{def:L} and satisfies \eqref{cond:L}, where $u_->u_+$. See, for instance, \cite{BRS,AmickBonaSchonbek} for more details on the well-posedness for \eqref{eqn:KdVB}. 

We assume that \eqref{eqn:main} has traveling wave solutions of the form $u(x,t)=\phi(x-ct)$ for some $c\in\mathbb{R}$, the wave speed, which satisfy
\begin{equation}\label{eqn:TravWave}
-c\phi_x+\phi\phi_x=\phi_{xx}+{\mathcal L}\phi.
\end{equation}
Our focus is on front solutions, which additionally satisfy 
\[
\phi(x)\rightarrow u_-~\text{as}~x\rightarrow-\infty\quad\text{and}\quad
\phi(x)\rightarrow u_+~\text{as}~x\rightarrow+\infty
\]
for some constants $u_{-} > u_{+}$. We remark that $c=\frac{u_{-}+ u_{+}}{2}$. This can be easily seen by integrating \eqref{eqn:TravWave} over $(-\infty,\infty)$ and noting that $\int_{-\infty}^\infty{\mathcal L}\phi~dx = \widehat{\mathcal{L}\phi}(0)=0$ by assumption, whereby $c(u_+-u_-)=\frac12(u_+^2-u_-^2)$. Moreover, we assume that
\[
\phi_x,\phi_{xx}\in L^2(\mathbb{R})\quad\text{and}\quad
\int^{\infty}_{-\infty}(1+|x|)\phi_x(x)~dx<\infty.
\]
These conditions hold for traveling front solutions of \eqref{eqn:KdVB}. 

Bona and Schonbek \cite{BS1985} established the existence of traveling front solutions of \eqref{eqn:KdVB}, uniquely determined up to translation. Of particular relevance to our purposes is that the profile is monotone if and only if $|\nu|\leq \frac14$. 

After suitable rescaling and a Galilean transformation, we set
\begin{equation}\label{eqn:BC1}
u_{-} = 1\quad\text{and}\quad\quad u_{+}=-1,
\end{equation}
so that $c=0$. To interpret, traveling front solutions of \eqref{eqn:main} become stationary. 
Our primary interest is in the stability of such stationary front solutions. 

In order to motivate our results, we first outline the most naive way in which one might attempt to make an energy argument, explain why it fails, how earlier works have addressed the issue, and how our approach differs. 

The most straightforward approach to stability would be to consider solutions of \eqref{eqn:main} and \eqref{eqn:BC} of the form
\begin{equation}\label{def:v}
u(x,t)=\phi(x-x_0)+v(x,t)
\end{equation}
for some constant $x_0\in\mathbb{R}$, where $v$ belongs to a suitable function space. Substituting $u(x,t)=\phi(x-x_0)+v(x,t)$ into \eqref{eqn:main} gives 
\[
v_t + v \phi_x + \phi v_x + v v_x = v_{xx} + \mathcal{L}v.
\]
Here we have used the fact that $\phi$ satisfies \eqref{eqn:TravWave}. 
Note that this is exact: there is no smallness assumption or linearization, and all terms have been retained. Multiplying this by $v$ and integrating, we arrive at 
\begin{align}
    \frac12 \frac{d}{dt}\int v^2 dx &= -\int v^2 \phi_x~dx - \int \phi v v_x~dx - \int v^2 v_x~dx + \int v v_{xx}~dx + \int v \mathcal{L} v~dx \notag \\
     &= -\int \phi_x v^2~dx + \frac12 \int \phi_x v^2~dx - \int v_x^2~dx \notag  \\
    &= -\frac12\int \phi_x v^2~dx - \int v_x^2~dx = -\langle v,-v_{xx}+\tfrac12\phi_x v \rangle = - \langle v,H_0 v\rangle \label{eqn:EnergyEstimate3}
\end{align}
Here and in the following, $\langle f,g\rangle=\int \!fg~\!dx$ represents the $L^2(\mathbb{R})$ inner product. 
The second equality uses an integration by parts for $\int \phi v v_x dx$, the fact that $\int v^2 v_x~dx=0$ and $\int v \mathcal{L} v~dx = 0.$ 
This is the case for the main case of interest where $\mathcal{L}=\nu\partial_{xxx}.$ In the case that $\mathcal{L}$ has a dissipative component, the above would become an inequality. 

This calculation is appealing as it shows that the rate of change of the $L^2$ norm of the perturbation is given by a quadratic form in $v$ determined by the operator $-\frac{d^2}{dx^2} + \frac12\phi_x(x)$, but unfortunately such a quadratic form cannot be definite. Since the potential term $\phi_x(x)$ has a negative integral 
\begin{equation}
\frac12\int^{\infty}_{-\infty}\phi_x(x)~dx=\frac12(\phi(\infty)-\phi(-\infty))=-1
\label{eqn:decreasing}
\end{equation}
it follows from the Rayleigh-Ritz theorem (see, for instance, \cite[Volume IV, Theorem XIII.3 p. 82]{RS}) that $\mathbf{H}_0$ must have a negative eigenvalue. To see this note that the existence of a test function $\Upsilon(x)\in L^2(\mathbb{R})$ such that $\int\Upsilon_x^2 + \frac12 \phi_x \Upsilon^2 dx <0 $ implies that the Schr\"odinger operator ${\mathbf{H}}_0=-\frac{d^2}{dx^2}+\frac12\phi_x(x)$ has at least one negative eigenvalue. To construct such a test function one defines $\Upsilon(x)$ in $L^2(\mathbb{R})$ that is equal to one over a large but compact set and decays slowly to zero outside this set. For such a function the term $\int v_x^2 dx$ can be made arbitrarily small and $\frac12\int \phi_x v^2 dx $ can be made arbitrarily close to $\frac12\int \phi_x dx=-1$. As a consequence the quadratic form cannot be definite and  $\|v(t)\|_{L^2(\mathbb{R})}$ is potentially an exponentially increasing function of $t$. 

An approach to making an energy argument work is due to Pego \cite{Pego}, expanding on the work of Goodman \cite{Goodman86}. Rather than working with the perturbation itself, Pego considered the anti-derivative of the perturbation defined by  
\begin{equation}\label{def:V}
V(x,t): = \int_{-\infty}^x (u(y,t) -\phi(y-y_0))~dy, 
\end{equation}
where $y_0$ is chosen such that $V(x,0)\to0$ as $x\rightarrow \infty$, and showed that 
\begin{equation}\label{eqn:energy(V)}
\frac12 \frac{d}{dt} \int V^2~dx = - \int V_x^2~dx + \frac12 \int \phi_xV^2~dx - \int V V_x^2~dx.   
\end{equation}
This is similar to \eqref{eqn:EnergyEstimate3} with two important differences. First, the term involving $\phi_x$ comes in with the opposite sign, whereby quadratic terms are strictly negative when $\phi_x(x)<0$ for all $x\in\mathbb{R}$. Second, the nonlinear terms no longer integrate to zero. 

To briefly summarize the remaining analysis, the first and second terms on the right side of \eqref{eqn:energy(V)} define the negative of a quadratic form that is manifestly negative semi-definite, provided $\phi_x(x)<0$ for all $x\in\mathbb{R}$---that is, $\phi(x)$ is a monotonically decreasing function for all $x\in\mathbb{R}$. The last term on the right side, on the other hand, is cubic and requires control. In essence, Pego \cite{Pego} demonstrated that $\|V(t)\|_{L^2(\mathbb{R})}^2$ is a decreasing function of $t$ so long as $\|V\|_{L^\infty(\mathbb{R})}$ remains sufficiently small, leading to a stability result for monotone profiles when subjected to small perturbations of the front.

Several other works (see, for instance, \cite{ACH2014,Wang,Engler}) address stability and related problems and they all seem to employ Pego's strategy \cite{Pego} of working with the anti-derivative of the perturbation of the front (see \eqref{def:V}), rather than directly on the perturbation itself (see \eqref{def:v}). Such approach typically necessitates the monotonicity of the profile. For further insights into stability of traveling front solutions, see \cite{Humpherys, ZumbrunHoward,ZumbrumHoward2K}.

Khodja \cite{Khodja} and, later, Naumkin and Shishmarev \cite{naumkinShishmarev,NaumkinShishmarevII} (see also \cite{NaumkinShishmarevIII}) adopted a perturbative approach to relax the monotonicity. However, their stability conditions, while less restrictive, are not particularly explicit and seem to be applicable for `slightly' non-monotone profiles. Our own numerical experiments indicate that the stability condition in \cite{naumkinShishmarev} holds for $|\nu|\lesssim 0.3$ for the KdVB equation. Despite these advancements, numerical studies (see, for instance, \cite{CanosaGazdag}) suggest the stability of all traveling front solutions to \eqref{eqn:KdVB}, irrespective of the monotonicity of the profile. Thus far, a rigorous mathematical treatment that accommodates strongly oscillatory profiles has remained elusive. 

The primary issue here is the translation invariance. The reason why the straightforward $L^2$ estimate does not work is that there is a one-parameter family of traveling wave solutions, and one really wants to look at the distance of $u(x,t)$ to {\em the set of all translates of} $\phi(x)$, rather than the distance to any particular translate. Pego's calculation addresses this problem by looking at the anti-derivative: the function $\int_{-\infty}^x (u(y,t)-\phi(y))~dy$ only tends to zero as $x\rightarrow\infty$ if an appropriate translate of $\phi$ is chosen, depending on $u$. But this has the unfortunate consequence of introducing an additional cubic term that must be estimated. We instead apply what is a fairly old idea to deal with the translation invariance by allowing the position of the front to vary in time in such a way as to make an energy estimate possible.  

For a broad class of nonlinear diffusive-dispersive equations of Burgers type, our objective here is to establish a stability condition that is independent of the monotonicity of the profile and does not require that the initial condition lies close to the front. It is noteworthy that some traveling front solutions of the KdVB equation, when extended with higher-order power-law nonlinearities, are known to exhibit instability \cite{PegoSmerekaWeinstein,Gardner}.

In Section~\ref{sec:analysis}, we work with the perturbation of the front (see \eqref{def:v}), rather than its anti-derivative (see \eqref{def:V}). In order to deal with the translational invariance, we introduce a temporal modulation of the translation parameter of the front (see \eqref{def:x0}), allowing the perturbed solution to `chase' the front, so as to `minimize' their distance. A judicious choice of the modulation introduces a rank-one perturbation into the Schr\"odinger operator that arises during the course of the energy method (see \eqref{eqn:H0+rank_1}).  
This perturbation has the potential to eliminate one negative eigenvalue in the spectrum of the unperturbed operator, enabling the application of the energy method in scenarios where the unperturbed Schr\"odinger equation (see \eqref{eqn:H0}) has precisely one bound state. 

It is worth noting that the idea of employing modulation to mitigate issues associated with translational invariance to enhance the stability properties has a long history, drawing on foundational ideas in Whitham's modulation theory \cite{Whitham}. Collet, Eckmann, Epstein, and Stubbe \cite{CEES93} and, independently, Goodman \cite{Goodman94} utilized this idea when bounding the dimension of the asymptotic attractor for the Kuramoto--Sivashinsky equation in one dimension (see \cite[Section 4 and (4.1)]{CEES93} in particular). Goodman and Miller \cite{GoodmanMiller} applied a similar idea when studying the stability of two-dimensional viscous shocks. Additionally, Gu\`es, Metivi\'er, Williams, and Zumbrun \cite{Gues2005} applied this when studying the vanishing viscosity limit, offering valuable historical insights and interpretations of various equations governing modulation (see \cite[Section 1]{Gues2005} in particular). 

In Section~\ref{sec:suff condition}, we tackle the challenge of verifying the non-existence of a second bound state of \eqref{eqn:H0}. Utilizing results by Jost and Pais \cite{JostPais} and Bargmann \cite{Bargmann}, we establish a sufficient condition for stability, based on a functional in Definition~\ref{def:tau} that can be interpreted as an $L^1$ distance between the monotonized front and the corresponding piece-wise constant (ideal) shock or, alternatively, as a measure of the front `width'. 
Analytical verification confirms that a traveling front solution of \eqref{eqn:KdVB} satisfies the sufficient condition for $|\nu|\leq \tfrac14$, signifying stability for $\nu$ in an open interval $\supset [-1/4,1/4]$ that encompasses all monotone profiles \cite{BS1985}. 
Building on our analytical findings, in Section~\ref{sec:numerics}, we conduct numerical experiments for \eqref{eqn:KdVB}, which indicate that the sufficient condition for stability holds for $|\nu|\lesssim 1.2$ and the stability condition itself for $|\nu|\lesssim 4.1$. 

Turning our attention to the KdVB equation, we embark on verifying the stability condition through validated numerics or rigorous computation techniques, culminating in a computer-assisted proof (CAP), demonstrating that the stability condition holds for $|\nu|$ within the significant range of $\ourint$, a notable achievement compared to the numerical prediction of $|\nu|\lesssim 4.1$. Consider 
\begin{equation}\label{eqn:Ricatti1}
-\frac{d^2v}{dx^2}+\frac12\phi_x(x)v=0\quad\text{and}\quad 
\text{$v(x)\rightarrow 1$ as $x\rightarrow\infty$}.
\end{equation}
It follows from the Sturm oscillation theorem that the number of zeros of $v$ is equal to the number of bound state of $-\frac{d^2}{dx^2}+\frac12\phi_x(x)$. Our objective here is to show that $v$ has precisely one zero, an endeavor well-suited to CAP techniques. This effort aligns our research with the growing application of CAP techniques to address open problems in differential equations. See, for instance, \cite{T2011,VL2015,HLMJ2016,CT2016,BJ2022,HLMJ2022,VdBMJR2016} for further insights into CAP techniques for differential equations. 

Section~\ref{sec:BigPicture} delineates the methodology underpinning our CAP approach, setting the stage, with further details elaborated in Section~\ref{sec:Supporting}. Notably, all computations are performed using interval arithmetic, ensuring rigorous control over rounding errors. This is facilitated by the use of Rump's INTLAB package for MATLAB \cite{Ru99a}.

We begin by establishing rigorous error bounds for numerical approximations of the traveling front solution of \eqref{eqn:KdVB}. Suppose
\[
\phi(x)=1+\sum_{n=1}^\infty \phi_n\theta^ne^{n\mu(\nu)x}, \quad x\leq 0,
\]
where $\phi_n$ are to be determined and $\mu(\nu)$ is in \eqref{def:mu}, representing the eigenvalue of the linearization about the fixed point $(\phi,\phi_x)=(1,0)$. The convergence of this series is guaranteed through comparison with a suitable geometric series. Numerical approximate solutions are computed over the interval $[x_N,\infty)$ for some $x_N\gg1$, along with solutions of \eqref{eqn:Ricatti1}. We then select $x_0(=0)<x_1<\cdots <x_N$ and compute numerical approximate solutions at each $x_j$, $j=0,1,\ldots,N$, as convergent Taylor series. Applying the Newton--Kantorovich theorem (see, for instance, \cite{NK1,NK2,Castelli2018, newton2, CHURCH2022133072} within the CAP context), we establish the existence of the unique exact solution in a neighborhood of our numerical approximate solution, accompanied by strict error bounds at each each $x_j$, $j=0,1,\ldots,N$, and, in turn, over the interval $[0,x_N]$. 

This rigorous computation allows us to determine how many times the solution of \eqref{eqn:Ricatti1} crosses zero within the interval $[0,\infty)$, ensuring that the derivative does not vanish near a zero, thereby confirming that the zero is simple. We then analytically show that $v$ has no zeros over the interval $(-\infty,0]$ provided that $\int^0_{-\infty}\sqrt{-\phi_x(x)}~dx$ is sufficiently small compared to $-v(0)$. This resembles Calogero's bound \cite{Calogero, Calogero2, Simon} on the number of bound states of Schr\"odinger equations with some potentials. We leverage the rigorous enclosure of the front to estimate $\int^0_{-\infty}\sqrt{-\phi_x(x)}~dx$. 

Our methodology, focused on a fixed value of $\nu$, is extended to an interval of $\nu$, using the Chebyshev interpolation of analytic functions \cite{reddy,Tadmor1986}. See Theorem~\ref{theorem:analytic_interpolation} for details. This enables approximations of various functions over small, overlapping intervals of $\nu$, while accounting for errors introduced by varying $\nu$.

Our CAP establishes that the stability result holds within the range $|\nu|\in \ourint$, covering nearly $95\%$ of the targeted range. Recall that our analytical proof has already established stability for $|\nu|\leq 0.25$. Presently, our effort is devoted to extending this result further to encompass the entire range $0.25\leq |\nu|\lesssim 4.09$. However, as we venture into higher values of $|\nu|$, we anticipate challenges arising from the widening of the profile, leading to increased computational time and memory issues. As $|\nu|$ approaches $0.25$, we encounter non-analytic behavior in the eigenvalues of the linearization about the fixed point $(\phi,\phi_x)=(-1,0)$. 

While our computations indicate the emergence of a second negative eigenvalue for $-\frac{d^2}{dx^2}+\frac12\phi_x(x)$ beyond the range $|\nu|\leq 4.1$, intriguingly, numerical evidence \cite{CanosaGazdag} suggests that spectral stability may extend to all $\nu\in\mathbb{R}$. It will be interesting to employ CAP techniques to rigorously establish spectral stability for $|\nu|\gtrsim 4.1$, particularly, as one approaches the KdV limit as $\nu\to\infty$.

 
\section{One bound state implies stability}\label{sec:analysis}

We begin with an important proposition, which shows that a Schr\"odinger operator with a single bound state can be made positive definite by an appropriate perturbation. We introduce some notation. 
\begin{defn}
 For $q \in L^2(\mathbb{R})$, let $\Pi_q$ denote the rank-one operator with range spanned by  $q$, defined by 
 \[
 \Pi_q v = q \int_{-\infty}^\infty q(x) v(x)~dx 
 \]
 Note that $\Pi_q$ is compact, self-adjoint, and positive definite -- up to normalization it is a projection onto a one-dimensional subspace. 
\end{defn}

We now demonstrate that a Schr\"odinger operator plus a particular rank-one perturbation can be made positive semi-definite, provided that the unperturbed Schr\"odinger equation has exactly one bound state. 

\begin{proposition}\label{prop:gamma}
Suppose $q\in L^2(\mathbb{R})$ is a real-valued function, ${\displaystyle \int^{\infty}_{-\infty}}(1+|x|)q(x)~dx <\infty$ and ${\displaystyle \int^{\infty}_{-\infty} q(x)~dx<0}$. If 
\[
H_0:=-\frac{d^2}{dx^2}+q(x)
\]
has precisely one negative eigenvalue, and zero is a regular value, in the sense that there do not exist bounded solutions to ${{H}_0}\psi=0$ then 
\[
H_\gamma:=-\frac{d^2}{dx^2}+q(x)+\gamma\Pi_q 
\]
is positive semi-definite for $\gamma>0$ sufficiently large. More specifically, $H_\gamma$ is non-negative for ${\displaystyle \gamma>-\frac{1}{\int^{\infty}_{-\infty} q(x)~dx}}$.
\end{proposition}

\begin{proof}
It is worth noting at the start that the projection term is a compact perturbation of $H_0$, whence it follows from the Weyl criteria \cite{RS} that the essential spectrum is unchanged. 

Recall that $(1+|x|)q(x) \in L^1(\mathbb{R})$ is a well-known condition which  ensures that $H_0$ has good spectral properties. In particular it guarantees that $H_0$ has a finite number of point eigenvalues that are negative, along with the continuous spectrum $[0,\infty)$, and eigenfunctions decays to zero like $e^{-\sqrt{-\lambda}|x|}$ as $|x|\to\infty$, where $\lambda$ is the corresponding negative eigenvalue. 

Let $n_-(H_\gamma)$ denote the number of negative eigenvalues of $H_{\gamma}$, counted with multiplicity. Since $\Pi_q$ is positive definite, $n_-(H_\gamma)$ is a non-increasing function of $\gamma$ over the interval $(0,\infty)$. Particularly, 
\[
n_-(H_0)\geq n_-(H_\gamma)\quad \text{for all $\gamma>0$}.
\]
More generally, 
\[
n_-(H_{\gamma_1})\geq n_-(H_{\gamma_2})\quad\text{whenever $\gamma_1<\gamma_2$}.
\]
Since $n_-(H_0)=1$ by hypothesis, it follows that $n_-(H_\gamma)\leq1$ for all $\gamma>0$. 
Note that
\[
H_{\gamma}1 = \left(1+\gamma \int_{-\infty}^{\infty} q(x)~dx\right) q(x).
\]
When $\gamma =-(\int \phi_x)^{-1}(>0)$, clearly, $H_\gamma v=0$ has a bounded solution $v=1$---namely, a zero-energy resonance. 
Since $H_0$ has one bound state and since $\Pi_q$ is positive definite, this suggests (and, indeed, we prove) that at this value of $\gamma$, a bound state is absorbed into the essential spectrum.

We proceed with some preliminaries. 
Let $\lambda_0$ denote the unique negative eigenvalue of $H_0$, and $v_0$ the corresponding eigenfunction, that is, $v_0$ is the ground state. Recall that $v_0(x)\rightarrow 0$ exponentially as $|x|\rightarrow\infty$. We may assume that $v_0(x)\geq 0$ for all $x\in\mathbb{R}$. We claim that
\[
\int^{\infty}_{-\infty} q(x)v_0(x)~dx <0.
\]
Indeed, we integrate
\[
-\frac{d^2v_{0}}{dx^2}+qv_0=\lambda_0 v_0,
\]
to obtain 
\[ 
\int_{-\infty}^{\infty} q(x)v_0(x)~dx = \lambda_0\int_{-\infty}^{\infty}v_0(x)~dx. 
\]
Since $\int v_0>0$ and $\lambda_0<0$, the left side must be negative. This proves the claim. 

Moreover, let $\lambda_\gamma$ be a negative eigenvalue of $H_\gamma$, and $v_\gamma$ the corresponding eigenfunction. We claim that 
\[
\int^{\infty}_{-\infty}q(x) v_\gamma(x)~dx\neq0.
\]
Suppose on the contrary that $\int qv_\gamma=0$, so that
\[
H_\gamma v_\gamma=H_0 v_\gamma=\lambda_\gamma v_\gamma.
\]
Since $H_0$ has exactly one negative eigenvalue, it follows that $\lambda_\gamma=\lambda_0$ and $v_\gamma$ is a constant multiple of $v_0$. Since $\int qv_0<0$, however, a contradiction proves the claim.
Additionally, the Hellmann--Feynman theorem \cite{FH} gives 
\[
\frac{d\lambda_\gamma}{d\gamma}=|\langle q,v_{\gamma}\rangle|^2 > 0.
\]
In other words, $\lambda_\gamma$ is a strictly increasing function of $\gamma$. 

We are now in a position to prove that $H_\gamma$ has no negative eigenvalues for $\gamma>-(\int q)^{-1}$. It is noteworthy that the result is related to the Aronszajn--Krein formula for finite-rank perturbation theory (see, for instance, \cite{Simon}).

Let $\gamma>0$. Suppose $\lambda_\gamma$ is a negative eigenvalue of $H_\gamma$, and $v_\gamma$ the corresponding eigenfunction, that is, 
\[
H_\gamma v_\gamma=H_0 v_\gamma+\gamma \langle q,v_\gamma\rangle q=\lambda_\gamma v_\gamma,
\]
so that 
\[
(H_0-\lambda_\gamma)v_\gamma=-\gamma \langle q,v_\gamma\rangle q.
\]
We can write 
\begin{equation}\label{eqn:u_gamma}
v_\gamma=-\gamma \langle q,v_\gamma\rangle (H_0-\lambda_\gamma)^{-1}q.
\end{equation}
Since $\lambda_\gamma$ is a strictly increasing function of $\gamma$, it follows that $\lambda_0<\lambda_\gamma$. Moreover, since $\lambda_0$ is the unique negative eigenvalue of $H_0$, it follows that $\lambda_\gamma$ is in the resolvent set of $H_0$, thereby $(H_0-\lambda_\gamma)^{-1}$ is well-defined.

Taking the inner product of \eqref{eqn:u_gamma} and $q$, we obtain
\[
\langle q,v_\gamma\rangle=-\gamma \langle q,v_\gamma\rangle\langle q,(H_0-\lambda_\gamma)^{-1}q \rangle.
\]
Since $\langle q,v_\gamma\rangle=\int qv_\gamma\neq0$, 
\[
1+\gamma\langle q,(H_0-\lambda_\gamma)^{-1}q\rangle=0.
\]
Let $R(\lambda)=1+\gamma\langle q,(H_0-\lambda)^{-1}q\rangle$, and observe the following.

Firstly $R(\lambda)$ is a Herglotz function, in that it maps the upper half-place $\Im(\lambda)>0$ to the upper half-plane and the lower half plane $\Im(\lambda)<0$ to the lower half plane. It it standard (see\cite{OPUC}) that such functions have their singularities on the real axis and are monotone when $f'(\lambda)$ is defined. 

Secondly since $H_0$ has a single negative eigenvalue we have that 
\[
R(\lambda) = 1 + \frac{\gamma |\langle q,v_0\rangle|^2}{\lambda_0-\lambda} + \gamma \langle q,\mathbf{P} (H_0-\lambda)^{-1} \mathbf{P} q \rangle 
\]
where $v_0$ is the ground state of $H_0$ and $\mathbf{P}$ is the operator of orthogonal projection onto the essential spectrum of $H_0.$ In particular  $\langle q,\mathbf{P} (H_0-\lambda)^{-1} \mathbf{P} q \rangle$ is analytic for $\lambda<0$ and thus $\lim_{\lambda\rightarrow \lambda_0^+}=-\infty$ due to the term $\frac{\gamma |\langle q,v_0\rangle|^2}{\lambda_0-\lambda}.$

Finally we have $\lim_{\lambda\rightarrow 0^-} R(\lambda) = \lim_{\lambda\rightarrow 0^-} 1+ \gamma\langle q,(H_0-\lambda)^{-1}q\rangle = 1 + \int q~\!dx$. To see this note that we can construct the Jost functions $\psi_\lambda^L(x)$ and $\psi_\lambda^R(x)$ which solve $(H_0-\lambda)\psi_\lambda^{L/R}=0$ together with the boundary conditions $\psi_\lambda^L(x)\rightarrow e^{\sqrt{-\lambda}x}$ as $x\rightarrow-\infty$ and  $\psi_\lambda^R(x)\rightarrow e^{-\sqrt{-\lambda}x}$ as $x\rightarrow+\infty.$ By variation of parameters we have the formula 
\[
(H_0-\lambda)^{-1} q = 1/W(\lambda) \left(\psi_\lambda^L(x) \int _x^\infty q(y) \psi_R(y) dy + \psi_\lambda^R(x) \int_{-\infty}^x q(y) \psi_\lambda^L(y) dy  
\right)\]
where $W(\lambda)$ is the Wronskian of the Jost solutions $W(\lambda)=\psi_\lambda^L\frac{d\psi_\lambda^R}{dx}-\psi_\lambda^R\frac{d\psi_\lambda^L}{dx}$. The Wronskian $W(\lambda)$ is non-zero for $\lambda \in (\lambda_0,0)$ since $H_0$ has only the one bound state at $\lambda_0$ and is non-zero at $\lambda=0$ by assumption, since $H_0\psi=0$ is assumed to have no bounded solutions. The function $(H_0-\lambda)^{-1} q$ exists as a function in $L_2(\mathbb{R})$ for $\lambda \in (\lambda_0,0).$ The limit $\lim_{\lambda\rightarrow 0^-} (H_0-\lambda)^{-1}q $ exists as a bounded function but is not, of course in $L_2(\mathbb{R}).$ It is straightforward to check that $H_0 1 = q$, and therefore we must have $\lim_{\lambda\rightarrow 0^-} (H_0-\lambda)^{-1}q =1$, since if we had a second bounded solution $\tilde\psi$ to $H_0 \tilde\psi = q$ then the difference $\tilde \psi-1$ would give a bounded solution to $H_0\psi =0$. By dominated convergence we have that $\lim_{\lambda \rightarrow 0^-}\langle q,(H_0-\lambda)^{-1} q\rangle = \langle q,1\rangle =\int_{-\infty}^\infty q(x)dx <0.$ Therefore since $R(\lambda)=1+\gamma\langle q,(H_0-\lambda)^{-1} q\rangle$ we have that $\lim_{\lambda\rightarrow 0^-}R(\lambda) = 1 + \gamma \int_{-\infty}^\infty q(x) dx$, which can be chosen to be strictly negative for $\gamma$ sufficiently large, since the potential by assumption has negative integral. Since $R(\lambda)$ is monotone increasing with $\lim_{\lambda\rightarrow\lambda_0^+} R(\lambda)=-\infty$ and $\lim_{\lambda\rightarrow 0^-}R(\lambda)=1+\gamma\int_{-\infty}^\infty q(x) dx<0 $ it follows that $R(\lambda)$ has no roots in $(-\lambda_0,0)$ for $\gamma$ sufficiently large. Therefore $H_\gamma$ has no negative eigenvalues for $\gamma$ sufficiently large.  

\end{proof}

We next state our first main analytical result. 

\begin{theorem}\label{thm:2.1}
Consider \eqref{eqn:main}-\eqref{eqn:BC}, where $\mathcal{L}$ is defined in \eqref{def:L}, satisfying \eqref{cond:L}, and $u_\mp=\pm1$. Let 
\begin{equation}\label{eqn:phi0}
\phi\phi_x=\phi_{xx}+{\mathcal L}\phi
\quad\text{and}\quad
\phi(x)\rightarrow\pm1~\text{as $x\rightarrow\mp\infty$}.
\end{equation}
Suppose $\phi_x,\phi_{xx}\in L^2(\mathbb{R})$ and ${\displaystyle \int^{\infty}_{-\infty}(1+|x|)\phi_x(x)~dx<\infty}$.
If 
\begin{equation}\label{eqn:H0}
-(1-\epsilon)\frac{d^2v}{dx^2}+\frac12\phi_x(x)v=\lambda v
\end{equation}
has precisely one bound state (corresponding to a negative eigenvalue) for some $\epsilon>0$ sufficiently small, then any solution of \eqref{eqn:main}-\eqref{eqn:BC}, subject to $u(\cdot,0)-\phi \in H^2(\mathbb{R})$, satisfies
\begin{align}
&\frac{d}{dt} \Vert u(\cdot,t)-\phi(\cdot-\xi_0(t)) \Vert_{L^2(\mathbb{R})}^2 \leq  0 \quad\text{for all $t>0$}\label{eqn:stabiliy1}
\intertext{and}
& \|u_x(\cdot,t)-\phi_x(\cdot-\xi_0(t))\|_{L^2(\mathbb{R})}^2\rightarrow 0 \quad\text{as $t\rightarrow\infty$}\label{eqn:stability2}
\end{align}
for some function $\xi_0$, thereby 
\begin{equation}\label{eqn:stability3}
\|u(\cdot,t)-\phi(\cdot-\xi_0(t))\|_{L^p(\mathbb{R})} \rightarrow 0\quad\text{as $t\rightarrow \infty$}
\end{equation}
for all $2<p\leq\infty$. 
\end{theorem}

Consequently, for a broad class of nonlinear diffusive-dispersive equations of Burgers type, a traveling front solution is asymptotically, nonlinearly, and orbitally stable, provided that an auxiliary Schr\"odinger equation has exactly one bound state. Importantly, our result has the potential to extend to non-monotone profiles, provided that the spectral condition is met. Particularly, we will show that the result applies to oscillatory profiles for the KdVB equation for a range of the dispersion parameter. 

\begin{remark*}
Our result depends only on the front solution, irrespective of the initial condition, as long as the initial condition differs from the front solution by a function in $H^2(\mathbb{R})$. Specifically, it {\em does not} require  that $u(\cdot,0)-\phi$ be small, representing a substantial improvement over earlier findings even for monotone profiles. In contrast, the results in \cite{Pego} and \cite{Khodja} require that $\|V(\cdot,0)\|_{H^3(\mathbb R)}$ is small, where $V$ is defined in \eqref{def:V}. Similarly the results in \cite{naumkinShishmarev, NaumkinShishmarevII} require $\|V(x,0)\Vert^2_{H^3(\mathbb R)}+\|x V(x,0) \Vert^2_{L^2(\mathbb R)}\ll 1$. 
\end{remark*}

\begin{proof}
Consider
\begin{equation}\label{def:x0}
u(x,t)=\phi(x-\xi_0(t)) + v(x,t),
\end{equation}
where $\phi$ represents a traveling front solution of \eqref{eqn:main}-\eqref{eqn:BC}, satisfying \eqref{eqn:phi0}, and $\xi_0(t)$ is the temporal modulation of the translation parameter, to be specified later (see \eqref{eqn:x0}). Comparing this with \eqref{def:v} we see that, rather than considering  the $L^2$ distance to a fixed translate of the front profile we are allowing the front to move in a way yet to be specified. Mathematically, this will result in an additional rank-one term in the quadratic form arising in our energy estimate. We will see that this additional degree of freedom will allow us to get a positive definite quadratic form and a corresponding energy estimate. 

Substituting \eqref{def:x0} into \eqref{eqn:main} and recalling \eqref{eqn:phi0}, after some algebra, we arrive at
\begin{equation}\label{eqn:v}
v_t ={\mathcal L}v+v_{xx}-vv_x-\phi(x-\xi_0(t))v_x-\phi_x(x-\xi_0(t))v+\dot\xi_0(t) \phi_x(x-\xi_0(t)), 
\end{equation}
where $\dot\xi_0=\frac{d\xi_0}{dt}$. Multiplying \eqref{eqn:v} by $v$ and integrating, we obtain
\[
\int v v_t~dx =\int_{-\infty}^\infty \left(v {\mathcal{L}}v + v v_{xx} - \phi v v_x - \phi_x v^2 + \dot\xi_0(t) \phi_x v\right) dx, 
\]
whence
\begin{equation}\label{eqn:energy(v)}
\begin{aligned}
\frac12\frac{d}{dt}\int^{\infty}_{-\infty}v^2~dx= 
&-\int^{\infty}_{-\infty} \left(v_x^2(x,t)+\frac12\phi_x(x-\xi_0(t))v^2(x,t)\right)~dx \\
&+\dot\xi_0(t) \int^{\infty}_{-\infty} \phi_x(x-\xi_0(t))v(x,t)~dx.
\end{aligned}
\end{equation}
Here we use integration by parts and that $\int^{\infty}_{-\infty} v^2v_x~dx =0$ and $\int v{\mathcal{L}} v~dx =0$, by assumption.

Suppose 
\begin{equation}\label{eqn:x0}
\frac{d\xi_0}{dt} = -\gamma\int^{\infty}_{-\infty} \phi_x(x-\xi_0(t))v(x,t)~dx
\end{equation}
for some positive constant $\gamma$, which will be specified later in \eqref{eqn:gamma}. There are various ways to interpret \eqref{eqn:x0}, but perhaps the most useful is that we are imposing gradient flow dynamics on the translation parameter, where $\xi_0$ evolves so as to minimize $\int (u(x,t)-\phi(x-\xi_0(t))^2~dx$. To see this, note that if we define the $L^2$ distance between the solution $u(x,t)$ and the front $\phi(x-\xi)$ by
\[
d(\xi)=\int (u(x,t)-\phi(x-\xi_0(t))^2~dx
\]
and calculate that 
\[
\frac{\partial d}{\partial\xi} = 2 \int_{-\infty}^\infty (u(x,t)-\phi(x-\xi))\phi_x(x-\xi)~dx = 2 \int v(x,t) \phi_x(x-\xi)~dx. 
\]
Therefore, the evolution of $\xi_0(t)$ is given by the gradient flow $\frac{d\xi_0}{dt} = -\frac{\gamma}{2}\frac{\partial d}{\partial \xi}(\xi_0(t)).$

Returning to the proof, substituting \eqref{eqn:x0} into \eqref{eqn:energy(v)}, we obtain 
\begin{align}\label{eqn:H0+rank_1}
\frac12\frac{d}{dt}\|v(t)\|_{L^2(\mathbb{R})}^2
=&-\int\left(v_x^2(x,t) +\frac12\phi_x(x-\xi_0(t))v^2(x,t)\right)~dx \notag \\
&- \gamma\left(\int \phi_x(x-\xi_0(t))v(x,t)~dx\right)^2 \notag\\
=&-\int \left(v_x^2(x+\xi_0(t),t) +\frac12 \phi_x(x) v^2(x+\xi_0(t),t) \right)~dx \notag \\
&- \gamma\left(\int \phi_x(x)v(x+\xi_0(t),t)~dx\right)^2 \notag \\ 
=&: -\int \left(w_x^2 +\frac12  \phi_x(x)w^2 \right)~dx - \gamma\left(\int \phi_x(x)w~dx\right)^2,
\end{align}
where $w(x,t)=v(x+\xi_0(t),t)$. Comparing \eqref{eqn:H0+rank_1} to \eqref{eqn:EnergyEstimate3} we see that there is one additional term in \eqref{eqn:H0+rank_1}, the term proportional to $\gamma$, which is negative and rank one. 
As mentioned previously the fact that $\int^{\infty}_{-\infty}\phi_x(x)~dx<0$ by \eqref{eqn:decreasing} implies that there is at least one negative eigenvalue of $-\frac{d^2}{dx^2}+\frac12\phi_x(x).$ What we will show is that if $\frac{d^2}{dx^2}+\frac12\phi_x(x)$ has only one negative eigenvalue then for $\gamma$ sufficiently large the quadratic form 
\[
\int \left(w_x^2 +\frac12  \phi_x(x)w^2 \right)~dx + \gamma\left(\int \phi_x(x)w~dx\right)^2
\]
is positive definite. This approach has the advantages that it does not require monotonicity of $\phi$ and that one gets an energy estimate for perturbations of essentially arbitrary size as long as they lie in an appropriate Sobolev space.

Returning to the proof of Theorem~\ref{thm:2.1}, recall \eqref{eqn:energy(v)} and \eqref{eqn:x0}, where
\begin{equation}\label{eqn:gamma}
\gamma>\gamma_0:=-\frac{2}{\int^{\infty}_{-\infty}\phi_x(x)~dx}=1.
\end{equation}
Recall $w(x,t)=v(x+\xi_0(t),t)=u(x+\xi_0(t),t)-\phi(x)$, and we calculate
\begin{align*}
\frac12\frac{d}{dt}\|v\|_{L^2(\mathbb{R})}^2=&\frac12\frac{d}{dt}\|w\|_{L^2(\mathbb{R})}^2\\
=&-\int^{\infty}_{-\infty} \left(w_x^2+\frac12\phi_x(x)w^2\right)~dx-\gamma\left(\int^{\infty}_{-\infty} \phi_x(x)w~dx\right)^2 \\
=&-\int \left((1-\epsilon)w_x^2+\frac12\phi_x(x)w^2\right)-\gamma\left(\int \phi_x(x)w\right)^2-\epsilon\int w_x^2 \\
\leq&-\epsilon\int w_x^2
\end{align*}
for some $0<\epsilon<1$. The last inequality uses that \eqref{eqn:H0} has exactly one bound state, so that $-(1-\epsilon)\frac{d^2}{dx^2}+\frac12\phi_x(x)+\gamma \Pi_{\phi_x}$ 
is positive semi-definite by Proposition~\ref{prop:gamma}. This proves \eqref{eqn:stabiliy1}. Additionally,
\[
 \Vert w(0)\Vert_{L^2(\mathbb{R})}^2 - \Vert w(T)\Vert_{L^2(\mathbb{R})}^2 \geq 2\epsilon \int_0^T \int^{\infty}_{-\infty}w_x^2(x,t)~dx dt
\]
for any $T>0$. Therefore, $w_x\in L^2(\mathbb{R}\times [0,\infty))$ and, hence, $v_x\in L^2(\mathbb{R}\times [0,\infty))$.

We wish to show that $\|w_x(t)\|_{L^2(\mathbb{R})}\rightarrow 0$  as $t\rightarrow\infty$ and, consequently, $\|v_x(t)\|_{L^2(\mathbb{R})}\rightarrow 0$ as $t\rightarrow\infty$. Since $w$ is related to $v$ by a (time-dependent) translation we have that $w$ satisfies
\begin{equation}\label{eqn:w}
w_t ={\mathcal L}w+w_{xx}-ww_x-\phi(x-\xi_0(t))w_x-\phi_x(x-\xi_0(t))w+\dot\xi_0(t) \phi_x(x-\xi_0(t)), 
\end{equation}
Taking the inner product of \eqref{eqn:w} and $-w_{xx}$, after some algebra, we arrive at 
\begin{equation}\label{eqn:h1}
\begin{aligned}
\frac12 \frac{d}{dt} \| w_x \|_{L^2(\mathbb{R})}^2 = &  -\int w_{xx} \mathcal{L}~\!w~\!dx - \int w_{xx}^2 dx  - \int \phi w_x w_{xx}  &\\
& - \int \phi_x~\!w~\!w_{xx} + \dot\xi\int \phi_x~\!w_{xx} + \int w~\!w_x~\!w_{xx}~\!dx  & \\
\frac12 \frac{d}{dt} \| w_x \|_{L^2(\mathbb{R})}^2 = &
\gamma \int \phi(x) w_x~dx \int \phi_{xx}(x) w_x~dx - \frac32\int \phi_x(x)w_x^2~dx \\
&- \int \phi_{xx}(x)ww_x~dx + \int ww_xw_{xx}~dx - \int w_{xx}^2~dx.
\end{aligned}
\end{equation}
Note that by assumption the dispersive term is a Fourier multiplier that commutes with differentiation and thus $-\int w_{xx}\mathcal{L} w~\!dx=\int w_x\mathcal{L}w_x~\!dx =0$. 
Applying the Cauchy--Schwartz and Sobolev inequalities we find
\begin{align*}
&\begin{aligned}\left|\gamma \int \phi w_x \int \phi_{xx}w_x \right|&= \left|\gamma \int \phi_x w \int \phi_{xx}w_x \right| \\
&\leq \|\phi_x\|_{L^2(\mathbb{R})} \|w\|_{L^2(\mathbb{R})} \|\phi_{xx} \|_{L^2(\mathbb{R})} \|w_x\|_{L^2(\mathbb{R})}=: C_1 \|w_x \|_{L^2(\mathbb{R})},\end{aligned} \\ 
&\left|\frac32 \int \phi_x w_x^2 \right|\leq \frac32 \Vert \phi_x \Vert_{L^\infty(\mathbb{R})} \Vert w_x \Vert_{L^2(\mathbb{R})}^2 =: C_2 \Vert w_x\Vert_{L^2(\mathbb{R})}^2,  \\
&\left|\int \phi_{xx}ww_x \right| \leq \Vert\phi_{xx}\Vert_{L^2(\mathbb{R})}~\Vert w\Vert_{L^\infty(\mathbb{R})}~\Vert w_x \Vert_{L^2(\mathbb{R})} =: C_3 \Vert w_x \Vert_2^{3/2}, \\
&\left| \int ww_xw_{xx} \right| \leq \frac12 \Vert w_{xx}\Vert_{L^2(\mathbb{R})}^2 + \frac12\Vert w w_x\Vert_{L^2(\mathbb{R})}^2  \leq \frac12 \Vert w_{xx} \Vert_2^2 + \Vert w \Vert_{L^2(\mathbb{R})}\Vert w_x \Vert_{L^2(\mathbb{R})}^3,
\end{align*}
where $C_1$, $C_2$, $C_3$, and $C_4:= \Vert w(x,0)\Vert_{L^2(\mathbb{R})}$ are positive constants, recalling that $\Vert w\Vert_{L^2(\mathbb R)}$ is decreasing and thus bounded. In the last inequality we have used the fact that 
\[
\Vert w\Vert_\infty^2 \leq 2 \Vert w\Vert_{L_2} \Vert w_x \Vert_{L_2}
\]
This follows from applying Cauchy-Schwartz to $w^2(x)=\int_{-\infty}^x 2 w(y) w'(y) dy $

Substituting these inequalities into \eqref{eqn:h1}, 
\[
\frac12 \frac{d}{dt} \Vert w_x\Vert_{L^2(\mathbb{R})}^2 \leq C_1 \Vert w_x \Vert_{L^2(\mathbb{R})} + C_2 \Vert w_x\Vert_{L^2(\mathbb{R})}^2 + C_3 \Vert w_x\Vert_{L^2(\mathbb{R})}^{3/2} + C_4 \Vert w_x\Vert_{L^2(\mathbb{R})}^3.  
\]

To simplify notation, let 
\[
\zeta(t):=\|w_x(t)\|_{L^2(\mathbb{R})}^2, 
\]
so that
\begin{equation}\label{eqn:diff_ineq}
\frac12\left|\frac{d\zeta}{dt}\right| \leq C_1 \zeta^{1/2}+C_2 \zeta + C_3 \zeta^{3/4} +  C_4 \zeta^{3/2}.
\end{equation}
Since $\|w_x\|_{L^2(\mathbb{R})} \in L^2([0,\infty))$ or, equivalently, $\zeta\in L^1([0,\infty))$, for any $\delta>0$ there exists some $T>0$ sufficiently large such that 
\[
\zeta(T)<\delta\quad\text{and}\quad \int_T^\infty \zeta(t)~dt<\delta.
\]
Let $T'>T$, and 
\[
M:=\max_{t\in(T,T')}\zeta(t).
\]
Multiplying \eqref{eqn:diff_ineq} by $\zeta^{1/2}$ and integrating over $(T,T')$, 
\begin{align*}
\frac13 (M^{3/2} - \delta^{3/2}) 
\leq \frac12 \int_T^{T'} & \zeta^{1/2} \left|\frac{d\zeta}{dt}\right|~dt \\
&\leq \int_T^{T'} (C_1 \zeta + C_2 \zeta^{3/2} + C_3 \zeta^{5/4} dt + C_4 \zeta^2)~dt.
\end{align*}
We calculate
\begin{align*}
& \int_T^{T'} \zeta(t)~dt \leq \int_T^\infty \zeta(t)~dt = \delta,  \\
& \int_T^{T'} \zeta^{3/2}(t)~dt \leq M^{1/2}\int_T^{T'} \zeta(t)~dt \leq M^{1/2}\int_T^\infty \zeta(t)~dt = \delta M^{1/2}, \\
& \int_T^{T'} \zeta^{5/4}(t)~dt \leq M^{1/4}\int_T^{\infty} \zeta(t)~dt\leq  \delta M^{1/4}, \\
& \int_T^{T'} \zeta^{2}(t)~dt \leq M\int_T^\infty \zeta(t)~dt = \delta M.
\end{align*}
Consequently,
\[
M^{3/2} \leq  \delta^{3/2} +3C_1 \delta + 3C_2 \delta M^{1/2} + 3 C_3 \delta M^{1/4} + 3 C_4 \delta M,
\]
implying that when $\delta>0$ is sufficiently small, 
$M\lesssim C\delta^{2/3}$ for some constant $C$, which is  clearly independent of $T'$. Therefore, for any $\delta>0$ there exists some $T>0$ such that $\|w_x(t)\|_{L^2(\mathbb{R})} \lesssim C\delta^{2/3}$ for all $t>T$, thereby
\[
\|w_x(t)\|_{L^2(\mathbb{R})} \to 0 \quad\text{as $t\rightarrow\infty$}.
\]
This proves \eqref{eqn:stability2}. Since $\|w\|_{L^\infty(\mathbb{R})}^2 \leq 2 \|w\|_{L^2(\mathbb{R})}\|w_x\|_{L^2(\mathbb{R})}$, moreover,
\[
\|w(t)\|_{L^p(\mathbb{R})} \to 0 \quad\text{as $t\rightarrow\infty$}
\]
for all $2<p\leq\infty$ and, hence, $\|w(t)\|_{L^\infty(\mathbb{R})} \to 0$ as $t\rightarrow\infty$. 
This proves \eqref{eqn:stability3} and Theorem~\ref{thm:2.1}. 
\end{proof}

\begin{remark*}
If one requires that $-\frac{d^2}{dx^2}+\frac12\phi_x(x)$ has one negative eigenvalue, instead of $-(1-\epsilon)\frac{d^2}{dx^2}+\frac12\phi_x(x)$, then $\|v(t)\|_{L^2(\mathbb{R})}$ is a decreasing function of $t$, but it does not necessarily decay to zero in any function space. Here we require that $-(1-\epsilon)\frac{d^2}{dx^2}+\frac12\phi_x(x)$ has one negative eigenvalue, to eliminate the possibility that $-\frac{d^2}{dx^2}+\frac12\phi_x(x)$ has a zero eigenvalue. 

If $\gamma$ is greater than the critical value (see \eqref{eqn:gamma}), then $\int \phi_x w \in L^2([0,\infty))$, but one should generally not expect that $\int \phi_xw$ is in $L^1([0,\infty))$, which would guarantee that $\lim_{t\rightarrow\infty} \xi_0(t)$ exists. If one considers the inviscid Burgers equation, which should have similar behavior to the KdVB at large scales and long time, it is not hard to construct examples where there is slow decay to the limiting state that leads to an unbounded but subballistic drift of the shock. For instance, if one takes an initial condition that is equal to $1$ for $x<0$ and decays to $-1$ like $x^{-1}$ for $x$ large then the location of the shock grows like $\log|t|.$ We believe that if one imposes additional decay requirements on $u(x,0)-\phi(x)$ then $\lim_{t\rightarrow\infty} \xi_0(t)$ should exist but we have not verified this.  
\end{remark*}

\section{Sufficient condition for one bound state}\label{sec:suff condition}

Verifying that \eqref{eqn:H0} has precisely one bound state is a nontrivial task. Here we establish a sufficient condition for this purpose.

\begin{defn}\label{def:tau}
Let $q$ be a real-valued function such that ${\displaystyle \int^{\infty}_{-\infty} (1+|x|) q(x)~dx<\infty}$. We define
\[
\tau(q) =\max_{x_0 \in (-\infty,\infty)} \min\left(\int_{-\infty}^{x_0} (x_0-x) q^-(x)~dx, \int_{x_0}^{\infty}(x-x_0) q^-(x)~dx\right),
\]
where $q^-(x):= \max(-q(x),0)$ represents the negative part of $q$.
\end{defn}

Suppose $\phi$ is a traveling front solution of \eqref{eqn:main}-\eqref{eqn:BC}, satisfying \eqref{eqn:phi0}, and $(1+|\cdot|)\phi_x \in L^1(\mathbb{R})$. 
We say that $\phi$ is {\em sharp} if $\tau(\tfrac12\phi_x)<1$ or $\tau(\phi_x)<2$.
 
\begin{defn}\label{def:shock}
Suppose $\phi$ is a real-valued function, satisfying, without loss of generality, $\phi(x)\rightarrow \mp1$ as $x\rightarrow \pm\infty$.  

If $\phi(x)$ is monotonically decreasing for all $x\in\mathbb{R}$, we define  
\[
s_\phi(x) = \begin{cases}
+1\quad \text{for $x<0$}, \\
-1\quad \text{for $x>0$}.
\end{cases}
\]
If $\phi$ is not monotone, let $m_\phi$ denote a piece-wise smooth function such that 
\[
{m_\phi}_x=-(\phi_x)^-
\quad\text{and}\quad 
m_\phi(x)\rightarrow \mp m_{\infty} \quad\text{as $x\rightarrow \pm \infty$},
\]
where ${\displaystyle m_\infty=-\frac{1}{2}\int_{-\infty}^{\infty} (\phi_x)^-(x)~dx}$,
and we define
\[
s_\phi(x)=\begin{cases}
+m_\infty\quad\text{for $x<0$},\\
-m_\infty\quad\text{for $x>0$}.
\end{cases}
\]
\end{defn}

When $\phi$ is monotone, $m_\phi=\phi$ and $m_\infty=1$. When $\phi$ is not monotone, on the other hand, the graph of $m_\phi$ is constructed by taking the graph of $\phi$ and replacing it by a constant over an interval where $\phi$ is increasing, so that $m_\phi$ is a piece-wise smooth function, decreasing where $\phi$ is decreasing but is constant where $\phi$ is increasing. 

\begin{figure}[htbp]
\centerline{
\includegraphics[width=0.45\textwidth]{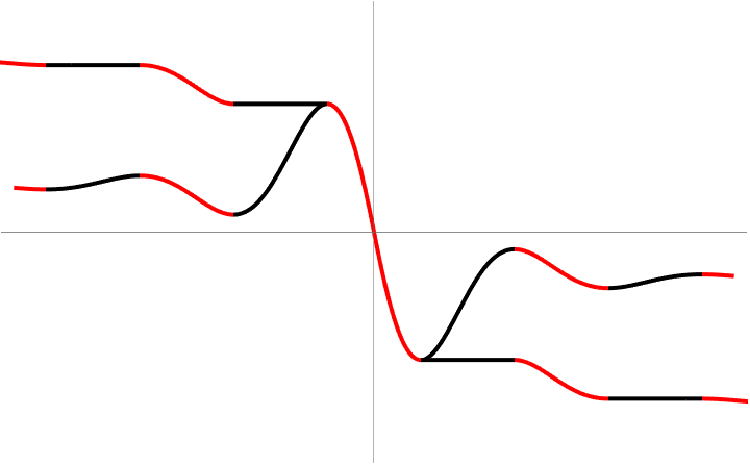}
\includegraphics[width=0.45\textwidth]{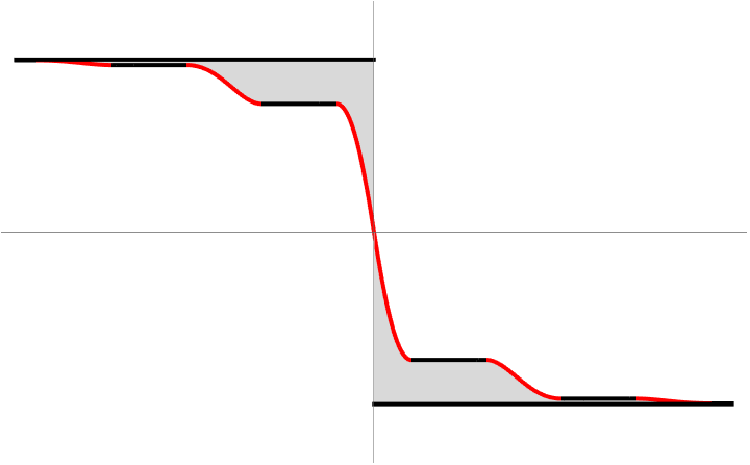}}
\caption{Left: The graphs of a non-monotone $\phi$ and its `monotonization' $m_\phi$. Red curves represent regions where $\phi$ is monotonically decreasing and $\phi_x={m_\phi}_x$, while black curves represent regions where $\phi$ is increasing and $m_\phi$ is constant. Right: The graphs of the monotonization $m_\phi$ and the corresponding ideal shock $s_\phi$. The shaded region, bounded by the profiles of $m_\phi$ and $s_\phi$, has an area~$=2\tau(\phi_x)$.} 
\label{fig:FrontPic}
\end{figure}

The left panel of Figure~\ref{fig:FrontPic} provides an example of a non-monotone profile and its `monotonization.' 

Suppose $\phi$ is a traveling front solution of \eqref{eqn:main}-\eqref{eqn:BC}. The content of the next lemma is that $\tau(\phi_x)$ can be interpreted as an $L^1$ distance to the corresponding piece-wise constant (ideal) shock or as the `width' of the profile. 

\begin{lemma}\label{lem:tau}
If $\phi$ is a real-valued function, satisfying $\phi(x)\rightarrow \mp1$ as $x\rightarrow\pm\infty$, and $(1+|\cdot|)\phi_x\in L^1(\mathbb{R})$, then
\[
\tau(\phi_x)=\frac12\int_{-\infty}^{\infty} |m_\phi(x) - s_\phi(x-x_0)|~dx,
\]
where $m_\phi$ and $s_\phi$ are in Definition~\ref{def:shock}, and $x_0$ is chosen so that
\begin{equation}\label{cond:area}
\int_{-\infty}^{x_0} (s_\phi(x-x_0)-m_\phi(x))~dx=\int^{\infty}_{x_0} (m_\phi(x)-s_\phi(x-x_0))~dx.
\end{equation}
\end{lemma}

Therefore, when $\phi$ has a monotone profile, $\tau(\phi_x)$ is the half of the area bounded by the graphs of $\phi$ and the corresponding piece-wise constant (ideal) shock, translated so that \eqref{cond:area} holds true. The right panel of Figure~\ref{fig:FrontPic} illustrates this. 

\begin{proof}
Let $x_0\in\mathbb{R}$, and
\[
f_-(x_0)=\int^{x_0}_{-\infty} (x_0-x) (\phi_x)^-(x)~dx\quad\text{and}\quad
f_+(x_0)=\int_{x_0}^{\infty} (x-x_0) (\phi_x)^-(x)~dx,
\]
so that 
\[
\tau(\phi_x) = \max_{x_0\in\mathbb{R}} \min(f_-(x_0),f_+(x_0)).
\]
It is easy to see that $f_-(x_0)$ is a positive increasing function of $x_0$, $f_-(x_0)\rightarrow 0$ as $x_0\rightarrow -\infty$ and $f_-(x_0)\rightarrow \infty$ as $x_0 \rightarrow\infty$. Also, $f_+$ is positive and decreasing, $f_+(x_0) \rightarrow 0$ as $x_0\rightarrow \infty$ and $f_+(x_0)\rightarrow \infty$ as $x_0 \rightarrow -\infty$. Consequently,  
\[
\tau(\phi_x)=f_-(x_0)=f_+(x_0)\quad\text{for a unique $x_0\in\mathbb{R}$}.
\] 

When $\phi$ has a monotone profile, $(\phi_x)^-=-\phi_x$ and integration by part leads to 
\begin{align*}
f_-(x_0) =& \int^{x_0}_{-\infty} (x-x_0)\phi_x(x)~dx \\
= &(x-x_0)(\phi(x)-1)\Big|^{x_0}_{-\infty} - \int^{x_0}_{-\infty} (\phi(x)-1)~dx = \int^{x_0}_{-\infty} (1-\phi(x))~dx. 
\end{align*} 
Similarly, ${\displaystyle f_+(x_0) = \int^{\infty}_{x_0} (1+\phi(x))~dx}$.
Since $s_\phi(x)=\begin{cases}+1\quad\text{for $x<0$}\\ 
-1\quad\text{for $x>0$}\end{cases}$, it follows that
\[
f_-(x_0)+f_+(x_0)=\int^{\infty}_{-\infty}|\phi(x)-s_\phi(x-x_0)|~dx.
\]
Moreover, $f_-(x_0)=f_+(x_0)$ provided that \eqref{cond:area} holds true. This proves the assertion. 

The proof is similar for non-monotone profiles. 
\end{proof}

\begin{remark*}
If $\phi$ or its monotonization is antisymmetric, the translate of the corresponding ideal shock $s_\phi(\cdot -x_0)$ minimizes the area between the front and $s_\phi(\cdot-x)$ over all $x\in\mathbb{R}$. However, for non-symmetric fronts, this is not necessarily true. 

It is noteworthy that \eqref{cond:area} resembles the Rankine--Hugoniot equal area condition. 
\end{remark*}

We state the main result of the section.

\begin{theorem}\label{thm:tau}
Suppose $\phi$ is a real-valued `front' like function, satisfying ${\displaystyle \int^{\infty}_{-\infty} (1+|x|)|\phi_x(x)|~dx<\infty}$ and ${\displaystyle \int^{\infty}_{-\infty} \phi_x(x)~dx<0}$. A sufficient condition for the operator 
\begin{equation}
-(1-\epsilon)\frac{d^2v}{dx^2}+\frac12\phi_x(x)v=\lambda v 
\label{eqn:Schro2}
\end{equation}
to have precisely one bound state for $\epsilon$ sufficiently small is that 
\[
\tau(\frac12\phi_x)<1.
\]
\end{theorem}

We will prove this in detail, which implies the asymptotic stability of sharp fronts.

\begin{corollary}\label{cor:stability}
Consider \eqref{eqn:main}-\eqref{eqn:BC}, satisfying $\mathcal{L}$ is defined in \eqref{def:L}, with \eqref{cond:L}, and $u_\mp=\pm1$. Let $\phi$ solve \eqref{eqn:phi0}. Suppose $\phi_x,\phi_{xx}\in L^2(\mathbb{R})$ and $(1+|x|)\phi_x\in L^1(\mathbb{R})$. If $\phi$ is sharp, that is, $\tau(\phi_x)<2$, then any solution of \eqref{eqn:main}-\eqref{eqn:BC}, subject to $u(\cdot,0)-\phi\in H^2(\mathbb{R})$, satisfies $\|u(\cdot,t)-\phi(\cdot-\xi_0(t))\|_{L^p(\mathbb{R})}\rightarrow0$ as $t\rightarrow\infty$ for some function $\xi_0(t)$ for all $2<p\leq\infty$. 
\end{corollary}

Therefore, for a broad class of nonlinear diffusive-dispersive equations of Burgers type, a traveling front solution is asymptotically, nonlinearly, and orbitally stable, provided that the profile lies within some $L^1$ distance from an appropriate translate of the corresponding ideal shock. Notably, our result does not require monotonicity of the front profile. We will show later that for the KdV--Burgers equation, this condition holds for strongly non-monotone shocks. Additionally, the result depends only on the front solution, and while the initial condition must differ from a front solution by something in $H^2(\mathbb{R})$, the difference need not be small in $H^2(\mathbb{R})$.

Corollary~\ref{cor:stability} follows from Theorems~\ref{thm:2.1} and \ref{thm:tau}. 
The proof of Theorem~\ref{thm:tau} relies on a well-known result concerning the non-existence of bound states of the Dirichlet boundary value problem for a Schr\"odinger equation. 

\begin{theorem}[Jost-Pais]\label{thm:JostPais}
Suppose $q$ is a real-valued function, satisfying ${\displaystyle \int_0^{\infty} xq^-(x)~dx < 1}$. There are no bound states to 
\[
\left\{
\begin{aligned}
&-\frac{d^2v}{dx^2} + q(x)v = \lambda v \quad \text{for $0<x<\infty$}, \\
& v(0)=0, \\
& v(x)\rightarrow 0\quad\text{as $x\rightarrow\infty$}.
\end{aligned}
\right.\]
\end{theorem} 

Jost and Pais \cite{JostPais} initially established Theorem~\ref{thm:JostPais} for Schr\"odinger equations in three dimensions with central fields of force when the angular momentum~$\ell=0$. Later, Bargmann \cite{Bargmann} extended the result to all values of the momentum.
The proof of Theorem~\ref{thm:JostPais}, as in \cite[pp. 94--97]{RS}, estimates the distance between the poles of the corresponding Riccati equation and is simpler than the proofs in \cite{JostPais} or \cite{Bargmann}.  

\begin{remark*}
Several other results address the non-existence of a second bound state for a Schr\"odinger equation in one dimension, and they could potentially be applicable here. For instance, Calogero \cite{Calogero} established various results of the kind. 
One well-known result by Calogero requires the monotonicity of the potential function, and it is difficult to incorporate such a result even when the front solution has a monotone profile because the location of the zero cannot be controlled. Nevertheless, such a Calogero bound will play a significant role later when complementing a computer-assisted proof. See, for instance, Proposition~\ref{prop:Calogero} and Lemma~\ref{lemma:Calogero} for details. There are other results that generalize the findings in \cite{JostPais}. See, for instance, \cite{Newton}. While any of these results could be useful here, the choice of \cite{JostPais} (or \cite{Bargmann}) is made because it is analytically simpler and has an interesting physical interpretation as the distance to a piece-wise constant shock.
\end{remark*}

We now present the proof of Theorem~\ref{thm:tau}.

\begin{proof}
Dividing \eqref{eqn:Schro2} by $1-\epsilon$ gives 
\begin{equation}
    -\frac{d^2v}{dx^2} + \frac{\phi_x(x)}{2(1-\epsilon)} v = \frac{\lambda}{1-\epsilon} v, \label{eqn:Schro3}
\end{equation}
so that \eqref{eqn:Schro2} and  \eqref{eqn:Schro3} have the same number of negative eigenvalues although of course the magnitudes will be scaled slightly relative to each other. Note that $\tau(\phi_x(x))<2$ is an open set condition, implying that $\tau(\frac{\phi_x}{1-\epsilon})=\frac{\tau(\phi_x)}{1-\epsilon}<2$ for $\epsilon$ sufficiently small. 

Suppose on the contrary that
\[
-\frac{d^2v}{dx^2} + q(x)v = \lambda v
\]
has two negative eigenvalues. Let $v_0$ and $v_1$ denote the ground state and the first excited state corresponding to the lowest and second lowest eigenvalues, respectively. It follows from the Strum oscillation theorem that $v_0(x)>0$ for all $x\in\mathbb{R}$, and there exists a unique $x_0\in\mathbb{R}$ such that $v_1(x_0)=0$ and $v_1(x)\neq0$ otherwise. Additionally, $v_1(x)\rightarrow 0$ as $|x|\rightarrow\infty$. Therefore, $v_1$ is a bound state to the Dirichlet problems to the left and right of $x_0$, satisfying 
\begin{align}
&-\frac{d^2v_1}{dx^2} + q(x)v_1 = \lambda v_1,&&v_1(x_0)=0,&& 
\text{$v_1(x)\rightarrow 0$ as $x\rightarrow-\infty$},\label{eqn:H0(-)}
\intertext{and}
&-\frac{d^2v_1}{dx^2} + q(x)v_1 = \lambda v_1,&&v_1(x_0)=0,&& 
\text{$v_1(x)\rightarrow 0$ as $x\rightarrow\infty$}.\label{eqn:H0(+)}
\end{align}
It then follows from Theorem~\ref{thm:JostPais} that 
\[
\int^{x_0}_{-\infty} (x-x_0)q^-(x)~dx\geq 1 \qquad\qquad
\int_{x_0}^{\infty} (x-x_0) q^-(x)~dx\geq 1.
\]
Although the value of $x_0$ which maximizes the smaller of the two is not known a priori, this implies that 
$\tau(q)\geq 1$. Therefore $\tau(q)<1$ implies that there can exist only a single bound state. Taking $q(x)=\frac{1}{2(1-\epsilon)} \phi_x$ shows 
\[
\tau\Big(\frac{\phi_x}{2(1-\epsilon)} \Big)<1 \quad \text{or, equivalently,}\quad
\tau(\phi_x) < 2(1-\epsilon)
\]
is sufficient to guarantee a single bound state. Since this is an open set condition, it suffices to require that $\tau(\phi_x)<2.$
\end{proof}

Concluding this section, we apply our result to the KdVB equation. 

\begin{proposition}\label{prop:KdVB}
Let
\begin{equation}\label{eqn:phi(KdVB)} 
\phi_x+\nu\phi_{xx}=\tfrac12(\phi^2-1)
\quad\text{and}\quad
\phi(x)\rightarrow \mp1 \quad\text{as $x\rightarrow\pm\infty$}.
\end{equation}
If $\phi_x(x)<0$ for all $x\in\mathbb{R}$ then 
\[
\tau(\tfrac12 \phi_x)
\leq \tau_{\max}\approx 0.8375. 
\]
In other words, a monotone traveling front solution of \eqref{eqn:KdVB} is sharp.
\end{proposition}

Therefore, Corollary~\ref{cor:stability} implies that a monotone traveling front solution of the KdVB equation is asymptotically, nonlinearly, and orbitally stable. Recall \cite{BRS} that the profile is monotone if and only if $|\nu|\leq\tfrac14$. Since $\tau$ depends continuously on $\nu$ (see Definition~\ref{def:tau}), a traveling front solution of \eqref{eqn:KdVB} is stable for $\nu$ in an open interval $\supset [-\tfrac14,\tfrac14]$. 

\begin{proof}
Recall from the proof of Lemma~\ref{lem:tau} that 
\begin{equation}\label{eqn:beta'}
\tau(\phi_x)=\max_{x_0\in(-\infty,\infty)}\min\left(\int_{-\infty}^{x_0} (1-\phi(x))~dx, \int^{\infty}_{x_0} (1+\phi(x))~dx\right).
\end{equation}
Assume for the moment that $\nu=0$, so that the former equation of \eqref{eqn:phi(KdVB)} becomes 
\[
\phi_x=\tfrac12(\phi^2-1),
\]
which can be rearranged as
\[
1-\phi=-2\frac{\phi_x}{1+\phi}.
\]
Integrating over the intervals $(-\infty,x_0)$ and $(x_0,\infty)$, where $x_0\in\mathbb{R}$, and using the latter equation of \eqref{eqn:phi(KdVB)}, we obtain 
\begin{align*}
&\int_{-\infty}^{x_0} (1-\phi(x))~dx = 2\log(2) - 2\log|1+\phi(x_0)|
\intertext{and}
&\int_{x_0}^{\infty} (1+\phi(x))~dx = 2\log(2) - 2\log|1-\phi(x_0)|. 
\end{align*}
Let $\phi_0:=\phi(x_0)$, and $\phi_0$ ranges over $(-1,1)$ as $x_0$ ranges over $(-\infty,\infty)$. Therefore, \eqref{eqn:beta'} simplifies to
\[
\tau(\phi_x)=\max_{\phi_0\in (-1,1)} \min(2\log(2) - 2\log|1+\phi_0|, 2\log(2) - 2\log|1-\phi_0|).
\]
Since $\log|1+\phi_0|$ is an increasing function of $\phi_0$, and $\log|1-\phi_0|$ decreasing we have that 
\[
\tau(\phi_x)=2\log(2)<2
\]

The proof for $\nu\neq 0$ is similar, but we will have to deal with inequalities rather than equalities. We assume without loss of generality that $\nu>0$. Indeed, \eqref{eqn:phi(KdVB)} remains invariant under
\[
x\mapsto -x, \quad \phi\mapsto -\phi \quad\text{and}\quad \nu\mapsto -\nu.
\]

We discuss some preliminaries. We claim that
\begin{equation}\label{eqn:phi'1}
\phi_x(x)\geq-\tfrac12\quad \text{for all $x\in\mathbb{R}$}.
\end{equation}
Indeed, $\phi_{xx}=0$ at a critical point of $\phi_x$, and the former equation of \eqref{eqn:phi(KdVB)} becomes $\phi_x=\frac12(\phi^2-1)$. Also, we claim that
\begin{equation}\label{eqn:phi'2}
\int^{\infty}_{-\infty}\phi_x(x)^2~dx=\frac23.
\end{equation}
Indeed, multiplying the former equation of \eqref{eqn:phi(KdVB)} by $\phi_x$ and integrating, 
\[
\int^{\infty}_{-\infty} (\phi_x)^2~dx=\frac12\int^{\infty}_{-\infty}\phi_x(\phi^2-1)~dx
=\left[\frac16\phi^3(x)-\frac12\phi(x)\right]^{\infty}_{-\infty}=\frac23.
\]

Following the argument for $\nu=0$, we rearrange the former equation of \eqref{eqn:phi(KdVB)} as 
\[
1-\phi=-2\frac{\phi_x}{1+\phi}-2\nu\frac{\phi_{xx}}{1+\phi},
\]
and we calculate 
\begin{align}
\int_{-\infty}^{x_0} (1-\phi(x))~dx = & -2\int_{-\infty}^{x_0} \frac{\phi_x(x)}{1+\phi(x)}~dx - 2\nu \int_{-\infty}^{x_0} \frac{\phi_{xx}(x)}{1+\phi(x)}~dx \notag\\
= & 2\log(2) - 2\log|1+\phi(x_0)| - 2\nu
\int_{-\infty}^{x_0} \frac{\phi_{xx}}{1+\phi}~dx \notag \\
= & 2\log(2) - 2\log|1+\phi(x_0)|\notag \\ &-
2\nu\left(\frac{\phi_x(x_0)}{1+\phi(x_0)} + 2\int_{-\infty}^{x_0} \left(\frac{\phi_x}{1+\phi}\right)^2~dx \right)\notag \\
\leq & -2 \log\left|\frac{1+\phi_0}{2}\right| + \frac{\nu}{1+\phi_0}=:f_-(\phi_0),\label{def:f-}
\end{align}
where $\phi_0:=\phi(x_0)$. The third equality follows from an integration by parts and the last inequality uses \eqref{eqn:phi'1}. 
Similarly, 
\begin{align}
\int_{x_0}^{\infty}(1+\phi(x))~dx &= -2\int_{x_0}^{\infty} \frac{\phi_x(x)}{1-\phi(x)}~dx - 2\nu \int_{x_0}^{\infty} \frac{\phi_{xx}(x)}{1-\phi(x)}~dx \notag\\ 
&=  -2 \log\left|\frac{1-\phi(x_0)}{2}\right| +2\nu \frac{\phi_x(x_0)}{1-\phi(x_0)} +
2\nu \int_{x_0}^\infty \left(\frac{\phi_x}{1-\phi}\right)^2~dx \notag\\
&\leq -2 \log\left|\frac{1-\phi_0}{2}\right|  + \frac{4 \nu}{3(1-\phi_0)^2}=:f_+(\phi_0),\label{def:f+} 
\end{align}
where the last inequality uses that $1-\phi(x_0)>0$, $\phi_x(x_0)<0$ and \eqref{eqn:phi'2}. 

We wish to show that for $\nu\in(0,\frac14]$, so that $\phi_x(x)<0$ for all $x\in\mathbb{R}$,
\[
\tau(\phi_x) = \max_{x_0\in\mathbb{R}} \min \left(\int_{-\infty}^{x_0}(1-\phi(x))~dx, \int^{\infty}_{x_0} (1+\phi(x))~dx \right) < 2.
\]
Since 
\[
\int_{-\infty}^{x_0}(1-\phi(x))~dx\leq f_-(\phi_0) 
\quad\text{and}\quad
\int^{\infty}_{x_0}(1+\phi(x))~dx \leq f_+(\phi_0),
\]
where $f_-(\phi_0)$ and $f_+(\phi_0)$ are in \eqref{def:f-} and \eqref{def:f+}, 
it suffices to show for $\nu\in(0,\frac14]$ that
\begin{equation}\label{eqn:beta(f)}
\max_{\phi_0\in (-1,1)} \min (f_-(\phi_0),f_+(\phi_0))<2.
\end{equation}
Since $f_-$ and $f_+$ are strictly increasing functions of $\nu$, moreover, it suffices to show that \eqref{eqn:beta(f)} holds at $\nu=\tfrac14$, the largest value of $\nu$ for which $\phi_x(x)<0$ for all $x\in\mathbb{R}$. 

Let $\nu=\tfrac14$. Numerical evaluation suggests that there is a unique $\phi_*\approx -0.017$ such that 
\[
f_-(\phi_*)=f_+(\phi_*)\approx 1.675,
\]
whence 
\[
\max_{\phi_0\in (-1,1)} \min (f_-(\phi_0),f_+(\phi_0))\leq 2\tau_{\max}\approx 1.675.
\]
We turn this observation into a proof. If there exist $\phi_-$ and $\phi_+$ such that 
\[
f_-(\phi_-)<f_+(\phi_-) < 2\quad \text{and} \quad  f_+(\phi_+)<f_-(\phi_+) < 2,
\]
then $f_-(\phi_*)=f_+(\phi_*)$ for some $\phi_*\in(\phi_-,\phi_+)$ and  
\[
f_-(\phi_*)=f_+(\phi_*)<2.
\]
Indeed, let $\phi_-=-\frac{1}{59}$ and $\phi_+=-\frac{4}{235}$. A straightforward calculation reveals that
\begin{align*}
f_-\left(-\frac{1}{59}\right)=&\frac{59}{232} + 2 \log(\frac{59}{29}) \approx 1.6748 \\
< &f_+\left(-\frac{1}{59}\right)=\frac{3481}{10800} + 2 \log(\frac{59}{30}) \approx 1.675 < 2
\intertext{and} 
f_+\left(-\frac4{235}\right)=&\frac{235}{924} + 2 \log(\frac{470}{231}) \approx 1.675 \\
< &f_-\left(-\frac4{235}\right)=\frac{55225}{171363} + 2 \log(\frac{470}{231}) \approx 1.6748 < 2.
\end{align*}
Therefore,
\[
\tau(\tfrac12\phi_x)<\tfrac12 f_-(\phi_*)=\tfrac12 f_+(\phi_*)\approx .8375. 
\]
This completes the proof.
\end{proof}

\section{Numerical experiments}\label{sec:numerics}

We perform numerical experiments for the eigenvalue problem associated to the energy argument for the KdVB equation. Recall that in the previous sections we established that
\begin{enumerate}
\item $-(1-\epsilon)\frac{d^2}{dx^2} + \frac12 \phi_x(x)$ having precisely one negative eigenvalue is a sufficient condition to guarantee global asymptotic stability,
\item $\tau(\tfrac12 \phi_x)<1$ is a sufficient condition to guarantee that $-(1-\epsilon)\frac{d^2}{dx^2} + \frac12 \phi_x(x)$ has precisely one negative eigenvalue,
\item $|\nu|\leq \frac14$ is a sufficient condition to guarantee that $\tau(\frac12\phi_x)<1$ holds. 
\end{enumerate}
where $\phi$ is a traveling front solution of \eqref{eqn:KdVB}, satisfying \eqref{eqn:phi(KdVB)}. In this section we will explore, through standard numerics, the range for which conditions 1. and 2. hold.

It is noteworthy that conditions 1. and 2. hold for intervals of $\nu$ that are considerably larger than the monotonicity interval $[-\frac14,\frac14]$,  and that the condition 1. holds in an interval considerably larger than condition 2. This will examine the second condition via some careful numerical experiments as it plays a crucial role later in setting expectations for a computer-assisted proof. 

To numerically experiment with the former condition, we begin by computing the traveling front solution of \eqref{eqn:KdVB}, by solving \eqref{eqn:phi(KdVB)} using the `{ode45}' native ODE solver in Matlab. We then numerically integrate
\[
f_-(x):=  \frac12 \int_{-\infty}^x (x-y) (\phi_x)^-(y)~dy \quad\text{and}\quad 
f_+(x):= \frac12 \int_x^{+\infty} (y-x) (\phi_x)^-(y)~dy. 
\]
We identify a unique point $x_*$ where $f_-(x_*)=f_+(x_*)$, implying $\tau(\tfrac12\phi_x)=f_-(x_*)=f_+(x_*)$.

\begin{figure}[htbp]
\centerline{\includegraphics[scale=0.45]{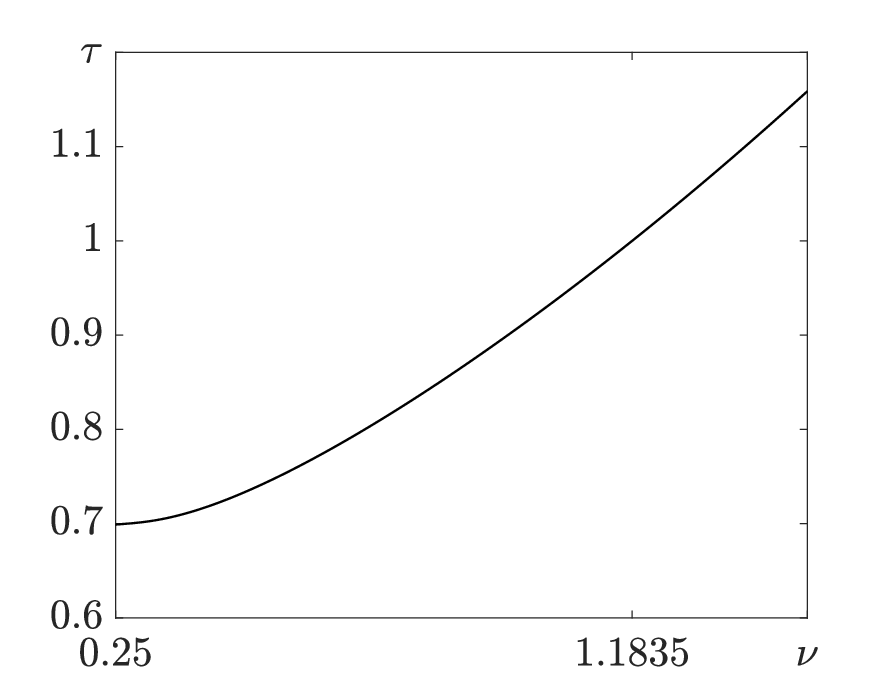}}
\caption{$\tau(\tfrac12\phi_x)$ versus $\nu$. We numerically find $\tau=1$ at $\nu\approx 1.1835$.}
\label{figure_evans}
\end{figure}

Our numerical findings suggest that $\tau(\tfrac12\phi_x)$ exhibits very weak dependence on $\nu$ within the range of $[-0.25,0.25]$, where $\phi$ is a monotonically deceasing function, and $\tau(\tfrac12\phi_x)$ appears to reach its minimum $\log(2)\approx 0.693$ when $\nu=0$. As $\nu$ increases from zero, $\tau(\tfrac12\phi_x)$ increases, reaching approximately $0.7$ when $|\nu|=0.25$. 

Figure~\ref{figure_evans} illustrates $\tau(\tfrac12\phi_x)$ versus $\nu$ for $\nu\in (0.25, 1.25)$. Our numerical results suggest that $\tau(\tfrac12\phi_x)$ is a convex function of $\nu$ with a local minimum at $\nu=0$. Furthermore, $\tau(\tfrac12\phi_x)$ surpasses the value of $1$ for the first time when $|\nu| \approx 1.1835$.   

We turn to numerical experiments with bound states of
\[
-\frac{d^2v}{dx^2}+\frac12\phi_x(x)=\lambda v,
\]
corresponding to negative eigenvalues.
We rewrite this as 
\begin{equation}\label{eigenproblemn}
\begin{pmatrix}v\\v_x\end{pmatrix}_x=\begin{pmatrix}0&1\\ \tfrac12\phi_x(x)-\lambda&0\end{pmatrix}\begin{pmatrix}v\\v_x\end{pmatrix}
=:\mathbf{\Lambda}(x,\lambda)\begin{pmatrix}v\\v_x\end{pmatrix}.
\end{equation} 
When $\lambda=k^2\geq 0$, $k\in\mathbb{R}$, the limiting matrix $\mathbf{\Lambda}(\pm\infty,\lambda)$ possesses a pair of purely imaginary eigenvalues $\pm ik$. When $\lambda<0$, on the other hand, $\mathbf{\Lambda}(\pm\infty,\lambda)$ has two real eigenvalues with opposite signs, $\pm\sqrt{-\lambda}$, and the corresponding eigenvectors are 
\[
\mathbf{v}_{1,+}:=\begin{pmatrix}-1\\\sqrt{-\lambda}\end{pmatrix}\quad\text{and}\quad 
\mathbf{v}_{2,-}:=\begin{pmatrix}1\\\sqrt{-\lambda}\end{pmatrix}.
\] 

Let $\lambda<0$, $\mathbf{v}_1(x,\lambda)$ and $\mathbf{v}_1(x,\lambda)$ denote solutions of \eqref{eigenproblemn} such that $\mathbf{v}_1(x,\lambda)\rightarrow0$ as $x\rightarrow\infty$, and $\mathbf{v}_2(x,\lambda)\rightarrow 0$ as $x\rightarrow -\infty$. We define the Evans function as
\begin{equation}\label{Evans}
\Delta(\lambda)=\det(\begin{pmatrix}\mathbf{v}_1(0,\lambda)& \mathbf{v}_2(0,\lambda)\end{pmatrix}).
\end{equation}
The significance of the Evans function is that $-\frac{d^2}{dx^2}+\frac12\phi_x(x)$ has precisely one negative eigenvalue if and only if $\Delta(\lambda)$ has precisely one negative root.
It follows from the gap lemma (see, for instance, \cite{GardnerZumbrun}) that 
\[\lim_{x\rightarrow \infty}\frac{\mathbf{v}_1(x,\lambda)}{|\mathbf{v}_1(x,\lambda)|}\in \Span\{\mathbf{v}_{1,+}\}\quad\text{and}\quad
\lim_{x\rightarrow -\infty}\frac{\mathbf{v}_2(x,\lambda)}{|\mathbf{v}_2(x,\lambda)|}\in \Span\{\mathbf{v}_{2,-}\}.
\] 

To compute the Evans function numerically, for $X\gg1$, consider
\[
\left\{\begin{aligned} &\mathbf{v}_{1,X}'=\mathbf{\Lambda}(x,\lambda)\mathbf{v}_{1,X}\\
&\mathbf{v}_{1,X}(X,\lambda)=\mathbf{v}_{1,+}\end{aligned} \right. \quad\text{and}\quad 
\left\{\begin{aligned} &\mathbf{v}_{2,X}'=\mathbf{\Lambda}(x,\lambda)\mathbf{v}_{2,X}\\
&\mathbf{v}_{2,X}(-X,\lambda)=\mathbf{v}_{2,-},\end{aligned}\right.
\]
and we define the numerical Evans function as 
\[
\Delta_X(\lambda)=\det(\begin{pmatrix}\mathbf{v}_{1,X}(0,\lambda) & \mathbf{v}_{2,X}(0,\lambda)\end{pmatrix}).
\]
To enhance numerical robustness, we scale out exponential rates. Introducing
\[
\mathbf{w}_1(x,\lambda):=e^{(x-X)\sqrt{-\lambda}}\mathbf{v}_{1,X}(x,\lambda)\quad\text{and}\quad 
\mathbf{w}_2(x,\lambda):=e^{-(x+X)\sqrt{-\lambda}}\mathbf{v}_{2,X}(x,\lambda),
\]
so that 
\[
\left\{\begin{aligned} &\mathbf{w}_1'=(\mathbf{\Lambda}(x,\lambda)+\sqrt{-\lambda}\mathbf{I})\mathbf{w}_1\\
&\mathbf{w}_1(X,\lambda)=\mathbf{v}_{1,+}\end{aligned} \right. \quad\text{and}\quad 
\left\{\begin{aligned} &\mathbf{w}_2'=(\mathbf{\Lambda}(x,\lambda)-\sqrt{-\lambda}\mathbf{I})\mathbf{w}_2\\
&\mathbf{w}_2(-X,\lambda)=\mathbf{v}_{2,-},\end{aligned}\right.
\]
we can redefine the numerical Evans function as
\[
\Delta_X(\lambda)=e^{2X\sqrt{-\lambda}}\det(\begin{pmatrix}\mathbf{w}_1(0,\lambda) & \mathbf{w}_2(0,\lambda)\end{pmatrix})=:e^{2X\sqrt{-\lambda}}\Delta_0(\lambda).
\]
Clearly, $\Delta_X(\lambda)=0$ if and only if $\Delta_0(\lambda)=0$. Hereafter we drop the subscript for simplicity.

It is crucial to exercise caution when interpreting the Evans function at $\lambda=0$. Note that $\mathbf{v}_1(x,0)$ and $\mathbf{v}_2(x,0)$ represent the Jost solutions to 
\[
-\frac{d^2v}{dx^2} + \frac12 \phi_x(x)v=0,
\]
which remain bounded as $|x|\rightarrow \pm \infty$. Therefore, when $\Delta(0)=0$, there exists a zero-energy resonance, signifying a bounded solution to $-\frac{d^2v}{dx^2} + \frac12 \phi_x(x)v=0$.

\begin{figure}[htbp]
\centerline{
\includegraphics[scale=0.6]{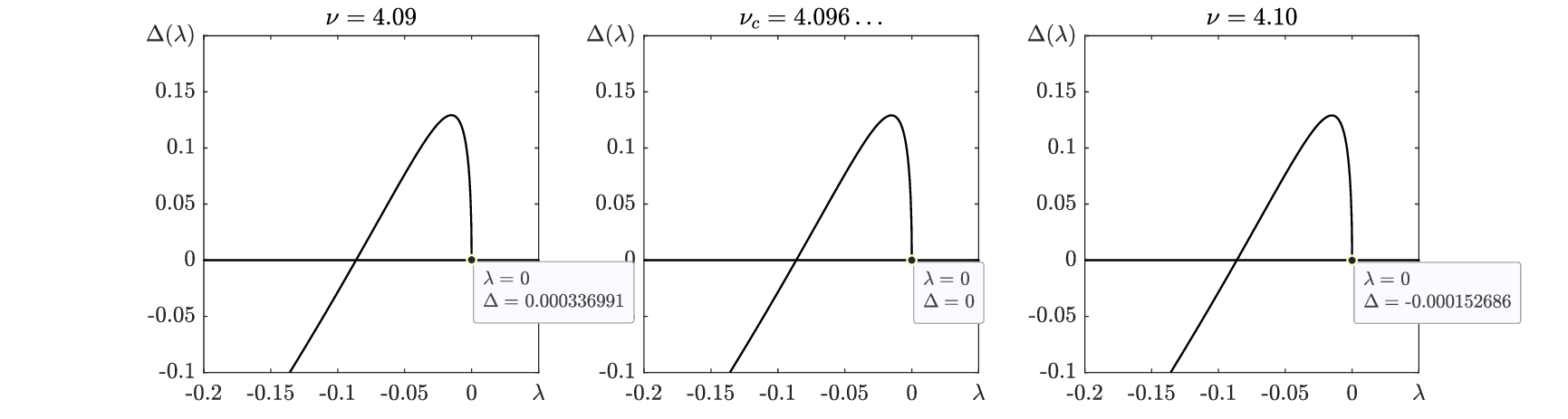}}
\caption{The graphs of $\Delta(\lambda)$ for $\nu=4.09$, $\nu_c\approx4.096$ and $4.10$. We numerically find that $\Delta(0)$ changes the sign passing through $\nu=\nu_c$.}
\label{figure_evans2}
\end{figure}

Utilizing Matlab, we numerically evaluate $\Delta(\lambda)$ for $\lambda<0$, to identify negative eigenvalues of the Schr\"odinger operator. Our numerical investigation suggests that $\Delta(\lambda)$ has at least two negative roots for $\nu>\nu_c$ for some $\nu_c$ and precisely one negative root for $0<\nu<\nu_c$. When $\nu=\nu_c$, we anticipate $\Delta(0)=0$. 
We compute the rescaled numerical Evans function at $\lambda=0$, and use the sign switch as a criterion for determining $\nu_c$. We find that $\Delta(0)$ changes its sign at $\nu_c\approx 4.096$. Figure~\ref{figure_evans2} illustrates the graphs of $\Delta(\lambda)$ for three different values of $\nu$. 

\section{Rigorous computation: an overview}\label{sec:BigPicture}

In summary, we have established two stability conditions for traveling front solutions of \eqref{eqn:main}-\eqref{eqn:BC}, applicable to a wide array of dispersive-diffusive operators. The first condition involves an auxiliary Schr\"odinger equation in one dimension. A traveling front solution, represented by $\phi$, is asymptotically, nonlinearly, and orbitally stable, provided that
\[
-(1-\epsilon)\frac{d^2v}{dx^2}+\frac{1}{2}\phi_x(x)v=\lambda v
\]
has precisely one bound state for some $\epsilon>0$ sufficiently small. A sufficient condition to ensure the non-existence of a second bound state is $\tau(\frac12 \phi_x)<1$, where $\tau$ is in Definition~\ref{def:tau}. Importantly, our findings are independent of the monotonicity of the profile and do not necessitate the initial condition to lie close to the front solution. Regarding the KdVB equation, recall from \cite{BS1985} that $\phi$ is monotone if and only if $|\nu|\leq\frac14$. We have analytically confirmed that $\tau(\frac12\phi_x)\leq \tau_{\max} \approx 0.8375$ for $|\nu|\leq\frac14$, thereby ensuring stability for $\nu$ within an open interval $\supset [-\frac14,\frac14]$. Numerical experiments suggest that the stability condition itself holds for $|\nu|\leq \nu_c\approx 4.09$.

In what follows, we focus on the KdVB equation and work to verify the stability condition through validated numerics or rigorous computation, with the ultimate goal of achieving a computer-assisted proof (CAP). Interested readers can find further details on CAP techniques for ODEs and PDEs in \cite{T2011,VL2015,HLMJ2016,CT2016,BJ2022,HLMJ2022,VdBMJR2016}. The code used to construct the CAP is available at \url{https://github.com/nonlinear-waves/stablab_matlab/tree/master/KdVB}. 

We state our main rigorous computation result.

\begin{theorem}\label{theorem1}
Suppose $\phi$ represents a traveling front solution of the KdVB equation, satisfying \eqref{eqn:phi(KdVB)}. If $|\nu|\in [0,0.25] \cup \ourint$, then \eqref{eqn:H0} has precisely one bound state for $\epsilon>0$ sufficiently small, thereby $\phi$ is asymptotically, nonlinearly, and orbitally stable. In other words, any solution of \eqref{eqn:KdVB} and \eqref{eqn:BC}, where $u_\mp=\pm1$, subject to $u(\cdot,0)-\phi\in H^2(\mathbb{R})$, satisfies $\|u(\cdot,t)-\phi(\cdot-\xi_0(t))\|_{L^p(\mathbb{R})}\rightarrow0$ as $t\rightarrow\infty$ for some function $\xi_0$ for all $2<p\leq\infty$.
\end{theorem}

Addressing the spectrum of an operator when the potential is not known explicitly poses a considerable challenge. To aid in verifying that \eqref{eqn:H0} has precisely one bound state, we introduced a proxy measure in Definition~\ref{def:tau}. This proxy provides a more analytically tractable method, yet it is important to note that its range of validity is more limited than that of the stability condition itself. The efficacy of CAP techniques lies in their capacity to directly determine the quantity of interest, thus extending our reach far beyond what was previously deemed attainable. Specifically, consider
\begin{equation}\label{eq:Riccati2}
-(1-\epsilon)\frac{d^2v}{dx^2}+\frac{1}{2}\phi_x(x)v=0 \quad\text{and}\quad v(x)\rightarrow 1~\text{as}~x\rightarrow \infty.
\end{equation}
It follows from the Sturm oscillation theorem that the number of zeros of a solution to \eqref{eq:Riccati2} is equal to the number of bound states of $-(1-\epsilon)\frac{d^2}{dx^2}+\frac12\phi_x(x)$. Our goal is to prove that the solution of \eqref{eq:Riccati2} has precisely one zero, which aligns well with the objectives of CAP techniques.

Given the nature of CAP, it suffices to show at $\epsilon=0$ that a solution of \eqref{eq:Riccati2} has precisely one zero. Throughout our CAP process, we verify open set conditions---strict inequalities are required. Since a zero of the solution to  \eqref{eq:Riccati2} varies continuously with $\epsilon$, if all the strict inequalities  hold at $\epsilon=0$, they will continue to hold for some $\epsilon>0$ sufficiently small. Hereafter, we proceed with $\epsilon=0$.

We present an overview of our CAP methodology. Let $\nu\neq0$, and we begin by recasting the former equations of \eqref{eqn:phi(KdVB)} and \eqref{eq:Riccati2} as 
\eq{
\begin{cases}
\phi' = \psi\\ 
\psi'= \frac{1}{\nu}\left(-\psi+ \frac{1}{2}(\phi^2-1)\right)
\end{cases}
\quad\text{and}\quad 
\begin{cases}
v'=w\\ 
w'=\tfrac{1}{2}\psi v.
\end{cases}
}{\label{eq:the_ode}}
Here $(\phi,\psi)$ corresponds to the traveling front solution of \eqref{eqn:KdVB} and its derivative, which connect to the fixed points $(-1,0)$ and $(1,0)$ at $x=\pm\infty$, respectively, that is, 
\eq{
\lim_{x\rightarrow -\infty} (\phi,\psi)=(1,0)
\quad \text{and} \quad \lim_{x\rightarrow \infty} (\phi,\psi)=(-1,0).
}{\label{eq:bd}}
Similarly, $(v,w)$ represents the solution of \eqref{eq:Riccati2} and its derivative, which connect to $(1,0)$ at $x=\infty$, that is,
\eq{
\lim_{x\rightarrow \infty} (v,w)=(1,0).
}{\label{eq:initial}}
If one can derive sufficiently accurate numerical estimates for the solutions of \eqref{eq:the_ode}, \eqref{eq:bd}, and \eqref{eq:initial}, complemented by additional analytical estimates, then one can rigorously determine the number of zeros in the solution of \eqref{eq:Riccati2}, thereby proving Theorem~\ref{theorem1}. All computations utilize interval arithmetic for rigorous control over rounding errors, specifically through Rump's INTLAB package for MATLAB \cite{Ru99a}.

\subsection{Rigorous verification for a fixed \(\nu\)}\label{sec:single}

We proceed by setting a fixed value for $\nu$, and outline the conceptual framework for the proof of Theorem~\ref{theorem1}, with further details elaborated in Section~\ref{sec:Supporting}. 

Our approach involves validated numerics or rigorous computation, which entails acquiring precise error bounds for a numerical approximation to the solution of \eqref{eq:the_ode}, \eqref{eq:bd}, and \eqref{eq:initial}. It is important to note that $(\phi,\psi)$ solves a boundary value problem, comprising the first two equations of \eqref{eq:the_ode} and \eqref{eq:bd}, over the interval $(-\infty,\infty)$, whereas $(v,w)$ solves an initial value problem, comprising the latter two equations of \eqref{eq:the_ode} and \eqref{eq:initial}. We partition $(-\infty,\infty)$ into intervals $(-\infty,x_0]$, $[x_N,\infty)$, and $[x_0,x_N]$ for some $-\infty<x_0<x_N<\infty$. We may choose without loss of generality $x_0=0$ and we select $x_N$ sufficiently large. We then proceed to construct a left solution manifold, comprised of solutions of the former two equations of \eqref{eq:the_ode} and the former equation of \eqref{eq:bd} over the interval $(-\infty,0]$ and a right solution manifold, composed of solutions of \eqref{eq:the_ode}, the latter equation of \eqref{eq:bd}, and \eqref{eq:initial} over the interval $[x_N,\infty)$. The initial-boundary value problem \eqref{eq:the_ode}-\eqref{eq:initial} over $(-\infty,\infty)$ then reduces to the bounded interval $[0,x_N]$, introducing boundary and initial conditions at $0$ and $x_N$. To address this, we partition the interval $[0,x_N]$ into subintervals and construct a middle solution manifold, composed of solutions of \eqref{eq:the_ode} with specified initial conditions at each partition point. This methodology transforms the task of obtaining a rigorous numerical solution over the interval $(-\infty,\infty)$ into a connection problem of the left, right, and middle solution manifolds.

We begin by constructing the left solution manifold, focusing on solutions over the interval $(-\infty,0]$.

\begin{proposition}\label{prop:left sol}
Consider $x\leq 0$. There exists a one-parameter family of solution, denoted as $(\phi,\psi)(x;\theta)$, $\theta\in\mathbb{R}$, of the first two equations of \eqref{eq:the_ode} and the former equation of \eqref{eq:bd}, which admits a convergent exponential series representation. Specifically,
\begin{equation}
\phi(x;\theta)=\sum_{n=0}^{\infty} \phi_{n}(\mu)\theta^n e^{n \mu x}\quad\text{and}\quad 
\psi(x;\theta) = \sum_{n=0}^{\infty} \psi_{n}(\mu)\theta^n e^{n \mu x},\label{eqn:LeftM}
\end{equation}
where 
\begin{equation}\label{def:mu}
\mu=-\frac{1}{2\nu}(1-\sqrt{1+4\nu}),
\end{equation}
$\phi_n$ and $\psi_n$ satisfy 
\begin{align}
& (\phi_0,\psi_0) =(1,0), \qquad 
(\phi_1,\psi_1) = ( -1, -\mu)\label{eqn:RecurLeft2}
\intertext{and}
& \phi_n= \frac{1}{2(n-1)(n+1-n\mu)}\sum_{n'=1}^{n-1}{\phi_{n'}\phi_{n-n'}},\qquad
\psi_n= n\mu \phi_n,\quad n\geq2. \label{eqn:RecurLeft1}
\end{align}
Additionally, when $\nu>0$ and $|\theta|< 2$,  
\[
\Big|\frac{\partial^k}{\partial \theta^k}\Big(\sum_{n=N+1}^{\infty} \phi_n\theta^n\Big)\Big| \leq \frac{1}{2^{1-k}}\frac{\partial^kf_N}{\partial q^k}\Big(\frac{\theta}{2}\Big),
\quad k=0,1,\ldots,N,
\]
where $f_N(q)= \frac{q^{N+1}}{1-q}$.
\end{proposition}

For $x\leq 0$, we define  
\begin{equation}\label{defU}
\mathbf{u}_L(\theta;x) = \begin{pmatrix} \phi(x;\theta) \\ \psi(x;\theta) \end{pmatrix},
\end{equation}
where $\phi(x;\theta)$ and $\psi(x;\theta)$ are in \eqref{eqn:LeftM}, \eqref{def:mu}, \eqref{eqn:RecurLeft2} and \eqref{eqn:RecurLeft1}. We remark that $\theta$ accounts for translation invariance. 

The result of Proposition~\ref{prop:left sol}, under the assumption of convergence, facilitates rigorous numerical estimates for the solutions of the first two equations of \eqref{eq:the_ode} and the former equation of \eqref{eq:bd} over the interval $(-\infty,0]$. For the proof, the infinite series in \eqref{eqn:LeftM} is divided into a finite sum plus a residual term. The finite sum is computed using interval arithmetic, specifically through Rump's INTLAB package for MATLAB \cite{Ru99a}, for rigorous control of rounding errors. Additionally, an estimate for the residual is derived and then computed by interval arithmetic. See Section~\ref{sec:Supporting} for further details.

We move on to constructing the right solution manifold, which consists of numerical solutions over the interval $[x_N,\infty)$ for some $x_N\gg1$. Our approach here is similar to the construction of the left solution manifold, but includes a few important distinctions. Invoking \eqref{eq:initial}, we are set to construct all of $\phi$, $\psi$, $v$, and $w$. Additionally, when $\nu>\frac14$ and the fixed point $(\phi,\psi)=(-1,0)$ is a spiral sink, the stable manifold is represented by a double exponential series. 

\begin{proposition}\label{prop:right sol}
Assume $\nu>\frac14$. For $x\geq x_N$ for some $x_N\gg1$, there exists a two-parameter family of solutions, represented as $(\phi,\psi,v,w)(x;\theta_i,\theta_r)$, where $\theta_i, \theta_r\in\mathbb{R}$, of \eqref{eq:the_ode}, the latter equation of \eqref{eq:bd}, and \eqref{eq:initial}, which admits a convergent double exponential series representation. Specifically,
\begin{equation}
\begin{aligned}
\label{eqn:RightM}
\phi(x;\theta_i,\theta_r)&=\sum_{m,n=0}^{\infty} \phi_{m,n}(\theta_r e^{\mu_- x - i\theta_i})^m(\theta_r e^{\mu_+ x +i\theta_i})^n,\\
\psi(x;\theta_i,\theta_r)&=\sum_{m,n=0}^{\infty} \psi_{m,n}(\theta_r e^{\mu_- x - i\theta_i})^m(\theta_r e^{\mu_+ x +i\theta_i})^n,\\
v(x;\theta_i,\theta_r)&=\sum_{m,n=0}^{\infty} v_{m,n}(\theta_r e^{\mu_- x - i\theta_i})^m(\theta_r e^{\mu_+ x +i\theta_i})^n, \\
w(x;\theta_i,\theta_r)&=\sum_{m,n=0}^{\infty} w_{m,n}(\theta_r e^{\mu_- x - i\theta_i})^m(\theta_r e^{\mu_+ x +i\theta_i})^n,
\end{aligned}
\end{equation}
where
\begin{equation}\label{def:alpha;mu_pm}
\mu_{\pm}= -\frac{2(1\pm i\nu_0)}{1+\nu_0^2},\qquad \nu_0=\sqrt{4 \nu-1},
\end{equation}
and $\phi_{m,n}$, $\psi_{m,n}$, $v_{m,n}$, $w_{m,n}$ satisfy
\begin{equation}
\mat{\phi_{0,0} \\ \psi_{0,0} \\ v_{0,0}\\ w_{0,0}}=\mat{-1\\0\\1\\0}, \quad 
\mat{\phi_{1,0} \\ \psi_{1,0} \\ v_{1,0}\\ w_{1,0}}=s\mat{1\\\mu_-\\\frac{1}{2\mu_-}\\\frac{1}{2}},\quad 
\mat{\phi_{0,1} \\ \psi_{0,1} \\ v_{0,1}\\ w_{0,1}}=s\mat{1\\\mu_+\\\frac{1}{2\mu_+}\\\frac{1}{2}}, \label{eqn:RecurRightIC}
\end{equation}
and
\begin{align}
&\phi_{m,n}= K(m,n)\sum_{m'=0}^{m}\sum_{n'=0}^{n} \delta_{m,m',n,n'} \phi_{m',n'}\phi_{m-m',n-n'},\label{eqn:RecurRightFirst}\\
&\psi_{m,n} = (m\mu_-+n\mu_+)\phi_{m,n},\\
&v_{m,n}= \frac{1}{2(m\mu_-+n\mu_+)^2} \sum_{m'=0}^{m}\sum_{n'=0}^{n} v_{m',n'}\psi_{m-m',n-n'},\label{eqn:RecurRightv}\\
&w_{m,n}= (m\mu_-+n\mu_+)v_{m,n}, \label{eqn:RecurRightLast}
\end{align} 
$m,n\geq1$. Here
\begin{equation}\label{def:Omega}
K(m,n)= \frac{1+\nu_0^2}{2((m(1-i\nu_0)+n(1+i\nu_0))^2-2m(1-i\nu_0)-2n(1+i\nu_0)+1+\nu_0^2)}
\end{equation}
\[
\delta_{m,m',n,n'} =\begin{cases}0&\quad \textrm{if}\ m' = n' = 0\\ 0&\quad \textrm{if}\  m' = m\ \textrm{and}\ n' = n\\
1&\quad \textrm{otherwise.}\end{cases}
\]
Additionally, 
\begin{align*}
& |\phi_{m,n}| \leq C_0 C^{m+n}, 
&& |\psi_{m,n}| \leq \frac{1}{\sqrt{\nu}}C_0 C^{m+n},  \\
&|v_{m,n}| \leq C_1 C_3^{m+n}, 
&& |w_{m,n}| \leq C_2 C_4^{m+n} 
\end{align*}
for some suitable constants $C$, $C_0$, $C_1$, $C_2$, $C_3$ and $C_4$, depending on $\nu$, which will be determined later. Furthermore,
\begin{align*}
\Big|\phi(x)-\sum_{N'=0}^N \sum_{m+n=N'} & \phi_{m,n} (\theta_r e^{\mu_- x - i\theta_i})^m(\theta_r e^{\mu_+ x +i\theta_i})^n \Big| \\
&\leq C_0\sum_{N'= N+1}^{\infty}(N'+1)(\theta_rC)^{N'}\\
& =C_0\frac{d}{d(\theta_rC)}\frac{(\theta_rC)^{N+2}}{1-\theta_rC}\\
&=C_0\left(\frac{(N+2)(\theta_rC)^{N+1}}{1-\theta_rC}+\frac{(\theta_rC)^{N+2}}{(1-\theta_rC)^2}\right)
\end{align*}
uniformly for $x\in[x_N,\infty)$. Similar bounds apply to $\psi$, $v$ and $w$.
\end{proposition}

For $x\geq x_N$ for some $x_N\gg1$, we define  
\begin{equation}\label{defS}
\mathbf{u}_R(\theta_i;x,\theta_r)=\begin{pmatrix}\phi(x;\theta_i,\theta_r) \\\psi(x;\theta_i,\theta_r) \\  v(x;\theta_i,\theta_r) \\ w(x;\theta_i,\theta_r)\end{pmatrix},   
\end{equation}
where the components $(\phi,\psi,v,w)(x,\theta_i,\theta_r)$ are in \eqref{eqn:RightM}, \eqref{def:alpha;mu_pm}, \eqref{eqn:RecurRightIC}, and \eqref{eqn:RecurRightFirst}-\eqref{eqn:RecurRightLast}.

It is noteworthy that since $v_{0,0}=0$, the right side of \eqref{eqn:RecurRightv} does not depend on $v_{m,n}$.

Here $\mu_\pm$ are the eigenvalues of the linearization of the first two equations of \eqref{eq:the_ode} around the fixed point $(\phi,\psi)=(-1,0)$, and $\theta_r\pm\theta_i$ parameterize the solution manifold. Since the fixed point $(-1,0)$ is a spiral sink, these eigenvalues are complex. 
For our purposes, $\theta_r$ can be taken to be a suitable fixed constant while $\theta_i$ is to be solved for. The convergence of these series representations is guaranteed by dominating them with a suitable geometric series, and the error in the representations can be refined by estimating the residual of the geometric series. This strategy falls under the umbrella of the parameterization method, a technique detailed in, for instance, \cite{MR1976079, MR1976080, MR2177465,Castelli2015,Berg2016,VdBMJR2016}. 
Further details of the proof are related discussions are in Section~\ref{sec:Supporting}.

\begin{remark*}
The estimates in Proposition~\ref{prop:right sol}, while serving our purposes, are not optimized to give the strongest result possible. We have made a deliberate choice to favor simplicity over achieving maximal possible decay rates. The choice allows us to avoid intricate calculations involved with sub-exponential terms.  

It is important to note that \eqref{eqn:RightM}-\eqref{eqn:RecurRightLast} do {\em not} offer a smooth extension to $\nu=\frac14.$ In cases where the linearization of a nonlinear ODE flow possesses stable eigenvalues $\lambda_1$, $\lambda_2$,\ldots, an analytic conjugacy between the flow on the stable manifold and the linearized flow exists, provided that $\sum k_j \lambda_j \neq 0$, where $k_j$ are integers, excluding the trivial scenario where all $k_j=0$. Such a non-resonance condition is met for all $\nu$ values except $\frac14$. When $\nu=\frac14$, the correct expansion necessitates a sum of products of polynomials and exponentials due to the non-trivial Jordan form of the linearized flow. While methodologies capable of handling such expansions exist, they are not implemented here. We direct interested reader to \cite{VdBMJR2016} and the references therein for a more thorough exploration of this phenomenon and its resolution.
\end{remark*}

Moving forward, we partition the interval $[0,x_N]$ into $\bigcup_{j=0}^{N-1}[x_j,x_{j+1}]$, where $0=x_0<x_1<\cdots<x_N$. In practice, $x_0$, \ldots, $x_N$ are evenly spaced. Let $\mathbf{\Phi}(\mathbf{u};x):\mathbb{R}^4\to \mathbb{R}^4$ denote the solution operator for \eqref{eq:the_ode} with the initial condition $\mathbf{u}$ and the step size $x$, that is,
\begin{equation}\label{definePhi}
{\bf u}_j:=\begin{pmatrix}\phi(x_j)\\ \psi(x_j)\\ v(x_j)\\ w(x_j)\end{pmatrix},\qquad 
\mathbf{\Phi}({\bf u}_j;x):=\begin{pmatrix}\phi(x_j+x)\\ \psi(x_j+x)\\ v(x_j+x)\\ w(x_j+x)\end{pmatrix},
\end{equation}
where $\phi(x_j+x)$, $\psi(x_j+x)$, $v(x_j+x)$, and $w(x_j+x)$ are represented by Taylor series. 

\begin{proposition}\label{prop:middle sol}
Consider \eqref{eq:the_ode} and the initial conditions $\phi(x_j)$, $\psi(x_j)$, $v(x_j)$, and $w(x_j)$ at $x_j$. The solution $\phi(x_j+x)$, $\psi(x_j+x)$, $v(x_j+x)$, and $w(x_j+x)$ at $x_j+x$ possesses a convergent Taylor series representation. Specifically,
\begin{equation}
\begin{aligned}\label{eqn:MiddleM}
\phi(x_j+x)&=\sum_{n=0}^{\infty} \phi_n x^n, \quad
&\psi(x_j+x) &= \sum_{n=0}^{\infty} \psi_n x^n, \\
v(x_j+x) &= \sum_{n=0}^{\infty} v_n x^n, \quad
&w(x_j+x)&= \sum_{n=0}^{\infty} w_n x^n,  
\end{aligned}
\end{equation}
where $(\phi_0,\psi_0,v_0,w_0)=(\phi,\psi,v,w)(x_j)$,
\begin{align}
\phi_1=\psi_0,\quad \psi_1=\frac{1}{\nu}\Big(-\psi_0+\frac12(\phi_0^2-1)\Big),\quad v_1=w_0,\quad w_1=\frac12\psi_0v_0,\label{rec:middlefirst}
\end{align}
and
\begin{equation}\label{rec:middlelast}
\begin{aligned}
\phi_{n+1} &= \frac{\psi_n}{n+1},\quad &\psi_{n+1} &= \frac{1}{\nu (n+1)}\Big( -\psi_n+\frac{1}{2}\sum_{n'=0}^n \phi_{n-n'}\phi_{n'}\Big),\\
v_{n+1} &= \frac{w_n}{n+1},&w_{n+1} &= \frac{1}{2(n+1)}\sum_{n'=0}^n \psi_{n-n'}v_{n'}
\end{aligned}
\end{equation}
for $n\geq1$.
\end{proposition}

It is noteworthy that when $n=0$ in the second equation of \eqref{rec:middlelast}, the recurrence relation does not simplify to the second equation of \eqref{rec:middlefirst}. 

The coefficients of the series in \eqref{eqn:MiddleM} are determined by recurrence relations of the corresponding ODEs. Since the first terms of the series are chosen to satisfy the initial conditions, the series offer solutions to the initial value problems, provided they converge. Throughout the process of establishing convergence, it is particularly beneficial that we also derive truncation error bounds for the residual. Further details are discussed in Section~\ref{sec:Supporting}. For further exploration of validated numerics for solving ODEs with Taylor series, see, for instance, \cite{corliss1982solving,Moore1979,moore1984survey,Nickel}. 

\begin{figure}[htbp]
\begin{center}
$
\begin{array}{lcr}
\includegraphics[scale=0.5]{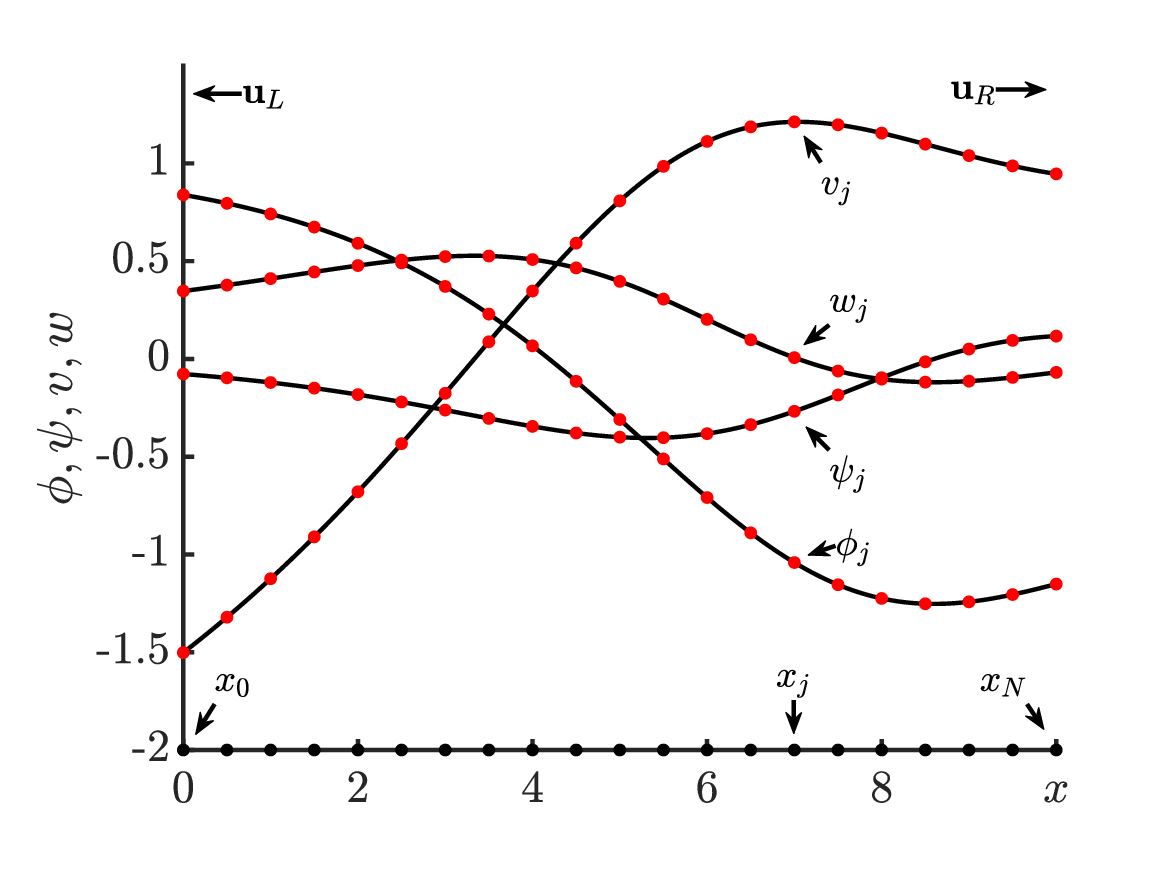}
\end{array}
$
\end{center}
\caption{Illustration demonstrating how to divide the computational domain and solve the ODEs with series solutions for various parts of the subdivided domains. The practical mesh we employed in our study is denser. 
}\label{fig6}
\end{figure}

Having established the solution manifolds for the intervals $(-\infty,0]$, $[x_N,\infty)$, and $[0,x_N]$, we aim to find a connecting solution whose segments on $(-\infty,0]$, $[x_N,\infty)$, and $[0,x_N]$ belong to the respective solution manifolds. This is accomplished by solving for the zero of ${\bf F}: {\mathbb R}^{4N+6} \rightarrow {\mathbb R}^{4N+6}$, defined as
\begin{equation}\label{def:boldF}
{\bf F}(\theta,{\bf u}_0,{\bf u}_1,\ldots, {\bf u}_{N},\theta_i) = \begin{pmatrix}\mathbf{u}_L(\theta;0) -\mathbf{\Pi} {\bf u}_0 \\
{\bf u}_1 - {\mathbf{\Phi}}({\bf u}_0;\tfrac{x_N}{N})\\
{\bf u}_2 - {\mathbf{\Phi}}({\bf u}_1;\tfrac{x_N}{N})\\
\vdots \\
{\bf u}_{N} - {\mathbf{\Phi}}({\bf u}_{N-1};\tfrac{x_N}{N})\\
\mathbf{u}_R(\theta_i;x_N,\theta_r) - {\bf u}_N \\
\end{pmatrix},
\end{equation}
where $\mathbf{u}_L$, $\mathbf{u}_R$, and $\mathbf{\Phi}$ are defined in \eqref{defU}, \eqref{defS}, and \eqref{definePhi}, respectively, and $\mathbf{\Pi}\begin{pmatrix}\phi\\ \psi\\ v\\ w\end{pmatrix}=\begin{pmatrix}\phi \\ \psi\end{pmatrix}$. The vanishing of the first two entries of the right side of \eqref{def:boldF} corresponds to the connection at $x_0=0$ of a solution from the left manifold with a solution from the middle manifold on the interval $[x_0,x_1]$. Similarly, the vanishing of the last two entries corresponds to the connection at $x_N$ of a solution from the right manifold with a solution from the middle manifold on the interval $[x_{N-1},x_N]$. The vanishing of the remaining entries establishes connections at $x_1$, \ldots, $x_{N-1}$ of solutions from middle manifolds on the intervals $[0,x_1]$, \ldots, $[x_{N-1},x_N]$. Therefore, a zero of \eqref{def:boldF} corresponds to a solution of the initial-boundary value problem \eqref{eq:the_ode}, \eqref{eq:bd}, and \eqref{eq:initial} over the interval $(-\infty,\infty)$. 

To establish the existence of an exact zero of \eqref{def:boldF} in the vicinity of an approximate zero, we employ the Newton--Kantorovich theorem, a technique extensively discussed in  \cite{NK1,NK2,Castelli2018,newton2,CHURCH2022133072}.

\begin{theorem}\label{thm:NewtK}
Consider $Y=\R^{d}$ for some integer $d>0$, and $y_{\rm ap}\in Y$. Suppose $\mathbf{F}:Y\to Y$ is continuously differentiable and $\mathbf{A}$, $\mathbf{A}^\dagger:Y \to Y$ are linear and bounded. Assume that:
\begin{itemize}[itemindent=11pt]
\item[{\rm (A1)}] $\mathbf{A}$ is one-to-one,
\item[{\rm (A2)}] $\|\mathbf{A}\mathbf{F}(y_{\rm ap})\|_Y \leq Y_0$ for some $Y_0>0$,
\item[{\rm (A3)}] $\|\mathbf{1}-\mathbf{AA}^{\dagger}\|_{B(Y)}\leq Z_0$ for some $Z_0>0$,
\item[{\rm (A4)}] $\|\mathbf{A}(\mathbf{A}^{\dagger}-D\mathbf{F}(y_{\rm ap}))\|_{B(Y)} \leq Z_1$ for some $Z_1>0$,
\item[{\rm (A5)}] There exists $Z_2:[0,\infty) \to [0,\infty)$ such that 
\[
\sup_{\|y-y_{\rm ap}\|\leq r} \|\mathbf{A}(D\mathbf{F}(y_{\rm ap})- D\mathbf{F}(y))\|_{B(Y)} \leq Z_2(r)r.
\]
\end{itemize}
If 
\begin{equation}\label{def:p}
p(r):=Z_2(r)r^2-(1-Z_0-Z_1)r+Y_0<0\quad\text{for some $r>0$},
\end{equation}
then there exists a unique $y_{\rm ex}\in Y$ satisfying $\|y_{\rm ex}-y_{\rm ap}\|_Y<r$ such that $\mathbf{F}(y_{\rm ex}) = 0$. 
\end{theorem}

It is important to note that while $Z_2(r)$ may be chosen to be arbitrarily large to ensure that {\rm (A5)} holds true, in practice, we select $r_0>0$ and show that {\rm (A5)} holds for some $r\in (0,r_0)$. Moreover, $\mathbf{A}$ and $\mathbf{A}^{\dagger}$ may be chosen arbitrarily as long as {\rm (A1)}-{\rm (A5)} are met. However, the most precise estimates are typically obtained when $\mathbf{A}^\dagger$ is close to the Jacobian $\nabla\mathbf{F}$ and $\mathbf{A}$ is an approximate inverse of $\nabla\mathbf{F}$. We employ numerical approximations of the Jacobian of $\mathbf{F}$ and its inverse for $\mathbf{A}^\dagger$ and $\mathbf{A}$. 

Our approach begins with computing an approximate zero $\theta_{\rm ap}$, ${\bf u}_{0,\rm ap}$, ${\bf u}_{1,\rm ap}$,\ldots, ${\bf u}_{N,\rm ap}$, and $\theta_{i,\rm ap}$ of \eqref{def:boldF} with sufficient accuracy. We then apply the Newton--Kantorovich theorem to demonstrate that this approximate zero is in close proximity to the exact zero $\theta_{\rm ex}$, ${\bf u}_{0, \rm ex}$, ${\bf u}_{1,\rm ap}$, \ldots, ${\bf u}_{N,\rm ex}$, and $\theta_{i, \rm ex}$ of \eqref{def:boldF}, along with explicit error bounds. To achieve this, all computations are carried out using interval arithmetic. We evaluate \eqref{def:boldF} and its Jacobian at the approximate zero of $\mathbf{F}$, using the series representations of the approximate solutions, and we bound its Hessian in an interval around the approximate zero of $\mathbf{F}$. The Newton--Kantorovich theorem guarantees the existence of a unique exact solution in a small-radius ball centered at the approximate zero. This provides control over the solution at $x_0$, \ldots, $x_{N}$ and control over the exact values of $\theta$ and $\theta_i$. We then utilize \eqref{definePhi} to propagate these bounds to the intervals $[x_0,x_1]$, \ldots, $[x_{N-1},x_N]$, and we leverage \eqref{defU} and \eqref{defS} to extend the bounds to the intervals $(-\infty,0]$ and $[x_N,\infty)$, respectively. This establishes rigorous upper and lower bounds of the solution of the former two equations of \eqref{eq:the_ode} and \eqref{eq:bd} over $(-\infty,\infty)$ and the solution of the latter two equations of \eqref{eq:the_ode} and \eqref{eq:initial} over $[0,\infty)$.


\begin{theorem}
\label{thm:verification1}
For each $|\nu| \in \ourint$, a numerical approximation, denoted as $(\phi_{\rm ap},\psi_{\rm ap},v_{\rm ap},w_{\rm ap})$, to the exact solution $(\phi_{\rm ex},\psi_{\rm ex},v_{\rm ex},w_{\rm ex})$ of \eqref{eq:the_ode}, \eqref{eq:bd}, and \eqref{eq:initial} can be rigorously validated. Specifically,  
\[
\|(\phi_{\rm ex}-\phi_{\rm ap},\psi_{\rm ex}-\psi_{\rm ap},v_{\rm ex}-v_{\rm ap},w_{\rm ex}-w_{\rm ap})\|_{L^\infty([x_j,x_j+1])}<\epsilon_j
\]
for some $\epsilon_j>0$ for some $0=x_0<x_1<\cdots<x_N$ and $x_N$ sufficiently large. Additionally,  
\begin{align*}
&\|(\phi_{\rm ex}-\phi_{\rm ap}, \psi_{\rm ex}-\psi_{\rm ap})\|_{L^\infty((-\infty,0])}<\epsilon_L
\intertext{and}
&\|(\phi_{\rm ex}-\phi_{\rm ap}, \psi_{\rm ex}-\psi_{\rm ap}, v_{\rm ex}-v_{\rm ap}, w_{\rm ex}-w_{\rm ap})\|_{L^\infty([x_N,\infty))}<\epsilon_R 
\end{align*}
for some $\epsilon_L$, $\epsilon_R>0$.
\end{theorem}

The rigorous enclosures of $\phi$ and $\psi$ over the interval $(-\infty,\infty)$, and $v$ and $w$ over $[0,\infty)$, enable us to determine the number of the zeros of $v$ over $[0,\infty)$. Since zeros of a solution of \eqref{eq:Riccati2} are necessarily simple, upon identifying an interval or intervals where $v$ has a zero, we only need to confirm that $w=v'$ is bounded away from zero to conclude that the zero is unique.

To ensure that $v$ has precisely one zero over the interval $(-\infty,\infty)$ and, correspondingly, \eqref{eqn:H0} has precisely one bound state, it remains to demonstrate that $v$ has no zeros over the interval $(-\infty,0]$. Since a rigorous enclosure of a solution of the latter two equations of \eqref{eq:the_ode} and \eqref{eq:initial} is lacking over the interval $(-\infty,0]$, we establish an analytic criterion that can be checked at $x=0$ that guarantees no zeros over $(-\infty,0]$. This is achieved in Lemmas \ref{lemma:BonaSchonbek}-\ref{lemma:Calogero}, establishing that if $\int_{-\infty}^0 \sqrt{-\phi_x(x)}~dx$ is sufficiently small compared to $v(0)$ and $w(0)$, then no more zeros of $v$ exist over the interval $(-\infty,0)$. This bears resemblance to Calogero's bound on the number of bound states of the Sch\"odinger equation for some potentials \cite{Calogero2,Calogero,RS}. We can estimate $\int_{-\infty}^0 \sqrt{-\phi_x(x)}~dx$ using  a rigorous enclosure of $\phi'$ over the interval $(-\infty,0]$, to ensure no zeros of $v$ over $(-\infty,0]$. 

\begin{proposition}\label{prop:Calogero} 
Consider
\[
-\frac{d^2v}{dx^2}+\frac12\phi_x(x)v=0.
\]
Suppose $\phi_x(x)<0$ and $\phi_{xx}(x)<0$ for all $x\in (-\infty,0)$. 
If $v(0)<0$, $v_x(0)>0$, and if
\begin{equation}
\arctan\left(-\sqrt{\frac{2}{-\phi_x(0)}}\frac{v_x(0)}{v(0)}\right) \geq \int_{-\infty}^0\sqrt{-\frac12\phi_x(x)}~dx \label{eqn:CalogeroBound}
\end{equation}
then $v$ has no zero over the interval $(-\infty,0)$.
\end{proposition}

Here $\|\phi_x\|_{L^{1/2}}$ can be bounded from above  using the result of Proposition~\ref{prop:left sol}, simplifying the matter to information at $x=0$. 

\subsection{Rigorous verification over an interval}\label{sec:interval}

To finalize our rigorous computation result for an interval in $\nu$ to complete the proof of Theorem~\ref{theorem1}, we will need to make necessary adjustments. The cornerstone result from the theory of analytic interpolation, essential for our methodology, provides a rigorous error bound for the interpolation of an analytic function using Chebyshev polynomials of the first kind. Further details can be found in \cite{Tadmor1986}, particularly Theorem 4.4 therein, or refer to \cite{reddy} for related literature.

\begin{theorem}\label{theorem:analytic_interpolation}
Let $\{T_m(x)\}_{m=0}^\infty$ denote the basis of  Chebyshev polynomials of the first kind, defined as
\[
T_0(x) = 1, \quad T_1(x) = x, \quad 
T_m(x) = 2xT_{m-1}(x)-T_{m-2}(x)\quad\text{for $m\geq 2$}.
\]
Consider $p_M(x) = \sum_{m=0}^{M-1} a_mT_m(x)$, an interpolating polynomial in the Chebyshev basis of degree $M-1$, with interpolation nodes at $\cos(\frac{\pi}{M}(m+\frac12))$ for $m=0,1,...,M-1$. 
If $f$ is a complex analytic function inside and on 
$E_{\rho}:=\{\frac{1}{2}(\rho e^{i\vartheta}+\rho^{-1} e^{-i\vartheta})\in\mathbb{C}: \rho>1~\text{and}~\vartheta \in [0,2\pi]\}$, then 
\[
\sup_{x\in[-1,1]}|f(x)-p_M(x)| \leq \frac{\sqrt{\rho^2+\rho^{-2}}}{ (\frac12(\rho+\rho^{-1})-1)}\frac{\max_{E_{\rho}}|f|}{\sinh(\log(\rho)(M+1))}.
\]
\end{theorem}

It is important to note that this result requires an a priori bound on the maximum modulus of $f$ over $E_\rho$. However, thanks to the exponential convergence in the number of collocation points, a relatively crude bound is sufficient if $M$ is chosen suitably large.


See Section~\ref{sec:interval2} for more details. 

\begin{remark}\rm
In Section~\ref{sec:suff condition}, we have demonstrated that the stability condition holds for $\nu$ in the range $[-\frac14,\frac14]$, where the profile is monotone. Additionally, Evans function computations have predicted that the condition continues to hold for $\nu<\nu_c\approx 4.09\ldots$.
For $\nu\in (\nuup,\nu_c)$, where $\nuup$ is the upper end of the range in Theorem~\ref{theorem1}, it becomes challenging to search for the optimal parameters so that \eqref{def:p} holds true. See Remark \ref{choiceofparameter} for more explanation.
\end{remark}

\section{Rigorous computation: proofs and supporting results}\label{sec:Supporting}

We provide further details of the proofs of the results in the previous section. We begin with an explanation of how to construct series solutions of \eqref{eq:the_ode} over various subintervals of $(-\infty,\infty)$. We employ induction arguments to demonstrate the convergence of these series wherever they are utilized, while leveraging CAP techniques to verify that the base cases of the induction arguments are valid. 

\subsection{Series solutions: the left, right, and middle solution manifolds}

\noindent{\it Proof of Proposition~\ref{prop:left sol}.}~(Construction of left solution manifold)~It is evident that \eqref{eqn:phi(KdVB)} has a saddle fixed point, $\phi=1$ and $\phi'=0$, corresponding to the eigenvalues $-\frac{1}{2\nu}(1\mp \sqrt{1+4\nu})$.
When constructing the one-dimensional unstable manifold for solutions of \eqref{eq:the_ode} near $x=-\infty$ for a range of $\nu$ of interest, it is advantageous to parameterize the solutions by $\mu= -\frac{1}{2\nu}(1- \sqrt{1+4\nu})$ rather than $\nu$. It is noteworthy that $\nu = \frac{1-\mu}{\mu^2}$, and $\mu\leq 1$ when $\nu\geq 0$. 

We assume that a solution over the interval $(-\infty,0]$ of the former two equations of \eqref{eq:the_ode} and the former equation of \eqref{eq:bd} takes the form \eqref{eqn:LeftM}, where $\mu$ is in \eqref{def:mu}.
Upon substitution, we verify \eqref{eqn:RecurLeft2} and \eqref{eqn:RecurLeft1}. We observe that the recurrence relation for $\phi_n$ admits a simple bound ensuring the convergence of \eqref{eqn:LeftM} for suitable values of $\mu$.  

\begin{lemma} \label{lemma:re_mu}
For $0<\Re(\mu)< 1$,
\[
|\phi_n|\leq \frac{1}{2^{n-1}},
\qquad n\geq 1.
\]
\end{lemma}

\begin{proof}
Recall that $\phi_1 = -1$. For $0<\Re(\mu)< 1$, 
\[
\left|\frac{1}{(n+1)-n\mu}\right|\leq 1
\quad\text{so that}\quad 
\left|\frac{1}{2(n-1)}\frac{1}{(n+1)-n\mu}\right| \leq \frac{1}{2(n-1)}.
\]
The assertion follows by induction.
\end{proof}

Clearly, $\nu \geq  \frac14$ implies $0<\mu \leq 2\sqrt{2}-2 (<1)$. Therefore, the result of Lemma~\ref{lemma:re_mu} holds for the range of $\nu$ of interest. 
We anticipate decay faster than exponential, although exploiting this can be challenging and has not proven necessary. 

To apply the Newton--Kantorovich theorem, recall that 
\[
\phi(\theta):=\phi(0;\theta)= \sum_{n=0}^{\infty} \phi_n \theta^n.
\] 
The truncation error arising from approximating $\phi(\theta)$ and its derivatives by finite sums can be controlled by the truncation error of geometric series, and can be rigorously justified by induction.

\begin{corollary}
For $|\theta|< 2$,
\[
\Big|\frac{\partial^k}{\partial \theta^k}\Big(\sum_{n=N+1}^{\infty} \phi_n\theta^n\Big)\Big|\leq \frac{1}{2^{1-k}}\frac{\partial^kf_N}{\partial q^k}\Big(\frac{\theta}{2}\Big), \quad k=0,1,..,N,
\]
where $f_N(q):= \frac{q^N}{1-q}$.
\end{corollary}

Additionally, to apply the Newton--Kantorovich theorem, we require a rigorous enclosure of 
\[
\psi(\theta):=\psi(0;\theta)= \sum_{n=0}^{\infty}\psi_n\theta^n,
\]
along with its first and second derivatives. This can be achieved using
\[
\psi(\theta) = \mu\theta\phi'(\theta),\quad \psi'(\theta) = \mu \phi'(\theta)+\mu\theta \phi''(\theta), \quad \psi''(\theta) = 2\mu \phi''(\theta)+\mu \theta \phi'''(\theta).
\]

\begin{tiny}
\begin{table}[ht]
\centering
\begin{tabular}{|c|c||c|c|c|c|c|c|c|c|c|c|c|c||} 
\hline 
$\nu_L$ & $\nu_R$ & $\theta$ & $N$ & $\rho$ & $M$ & Time\\
\hline
0.28 & 0.281 & 0.64 & 58 & 1.58e+03 & 7 & 4.9e+00 \\ 
0.5 & 0.501 & 0.64 & 58 & 3.46e+03 & 6 & 1.3e+01 \\ 
1 & 1.01 & 0.41 & 58 & 882 & 7 & 4.0e+00 \\ 
1.5 & 1.51 & 0.262 & 58 & 1.58e+03 & 7 & 1.4e+01 \\ 
2 & 2.01 & 0.168 & 58 & 2.34e+03 & 6 & 1.5e+01 \\ 
2.5 & 2.51 & 0.107 & 58 & 2.84e+03 & 6 & 1.4e+01 \\  
3 & 3.01 & 0.0687 & 58 & 3.29e+03 & 6 & 1.4e+01 \\ 
3.397 & 3.398 & 0.044 & 58 & 3.6e+04 & 5 & 4.0e+00 \\ \hline 
\end{tabular}
\caption{
Data concerning the rigorous enclosure of the one-dimensional unstable solution manifold, about $x=-\infty$, for the former two equations of \eqref{eq:the_ode} and the former equation of \eqref{eq:bd}, corresponding to the fixed point $(-1,0)$. We consider $\nu$ over the interval $[\nu_L,\nu_R]$ for various values of $\nu_L$ and $\nu_R$. Here $\theta$ is the parameter value at which the solution manifold is evaluated, $N$ is the index at which we truncate the infinite series in \eqref{eqn:LeftM}, $\rho$ is the radius of $E_\rho$ used in carrying out analytic interpolation of the coefficients, and $M$ is the number of Chebyshev coefficients used in the analytic interpolation. Additionally, `Time' indicates how long it took in seconds to compute the coefficients of the manifold.}\label{table:1}
\end{table}
\end{tiny}

\noindent {\it Proof of Proposition~\ref{prop:right sol}.}~(Construction of the right solution manifold)~We discuss the details of constructing the two-dimensional stable manifold for solutions of \eqref{eq:the_ode} near $x=\infty$, corresponding to the fixed point $(\phi,\psi,v,w)=(-1,0,1,0)$. Considering that later we will extend the result of Proposition~\ref{prop:right sol} to $\nu$ over an interval via analytic interpolation, we introduce $\nu_0= \sqrt{4\nu-1}$ for $\nu>\frac14$. We assume that a solution over the interval $[x_N,\infty)$ for some $x_N\gg1$ of \eqref{eq:the_ode}, the latter equation of \eqref{eq:bd}, and \eqref{eq:initial}, takes the form \eqref{eqn:RightM}, where $\mu_\pm$ are in \eqref{def:alpha;mu_pm}.
Upon substitution, we verify \eqref{eqn:RecurRightIC}-\eqref{eqn:RecurRightLast}.

We proceed to demonstrating that $\phi_{n,m}=\overline{\phi_{m,n}}$ for all $m,n=0,1,2,\ldots$, and determine the number of terms needed in the series to evaluate $\phi$, $\psi$, $v$, and $w$ within a desired error bound.

\begin{lemma}
For $\nu_0 \in \mathbb{R}$, $\phi_{n,m} = \overline{\phi_{m,n}}$, $m,n=0,1,2,\ldots$. 
\end{lemma}

\begin{proof}
Recall that $\phi_{0,0} = -1$ and $\phi_{1,0} = \phi_{0,1} = -s$. Let $m$, $n\in \{0\}\bigcup \N$ and $m+n\geq 1$, and assume that $\phi_{m',n'} = \overline{\phi_{n',m'}}$ for all $m',n'\in \{0\}\bigcup \N$ and $m'+n' \leq m+n$. Recall \eqref{def:Omega} and, clearly, $D(n,m)=\overline{K(m,n)}$. Recall \eqref{eqn:RecurRightFirst} and the assertion follows by induction.
\end{proof}

\begin{lemma}\label{lemma:Omega_bound}
For $m+n\geq 2$, $|K(m,n)|\leq {\displaystyle \frac{1+\nu_0^2}{2(m+n-1)^2}}$.
\end{lemma}

\begin{proof}
Recall \eqref{def:Omega}, and a straightforward calculation reveals
\begin{align*}
\Big|\frac{1}{K(m,n)}\Big|^2 &=\frac{4}{(1+\nu_0^2)^2}(\nu_0^{4}
((m - n)^{2} - 1)^{2} \\ 
&\qquad\qquad\qquad +2 \nu_0^{2}(m + n-1)^{2}((m - n)^2+1) + (m + n - 1)^{4})\\
&\geq \frac{4(m+n-1)^4}{(1+\nu_0^2)^2}.
\end{align*}
\end{proof}

\begin{lemma}\label{lemma:Omega_bound2}
For $m\geq 5$ or $n\geq 5$, 
${\displaystyle |K(m,n)| \leq \frac{1+\nu_0^2}{4 (mn + m+n -1)} }$. 
\end{lemma}

\begin{proof}
Note that the number of terms on the right side of \eqref{eqn:RecurRightFirst} is $(m+1)(n+1)-2= mn+m+n-1$, because $m'=0,1,\ldots, m$ and $n'=0,1,\ldots, n$, excluding $(0,0)$ and $(m,n)$. Recall from Lemma~\ref{lemma:Omega_bound} that $\frac{(1+\nu_0^2)}{2K(m,n)}\geq (m+n-1)^2$. Additionally, 
\[
(m+n-1)^2 - 2(mn+m+n-1)=(m-2)^2+(n-2)^2 -5>0
\]
for either $m\geq 5$ or $n\geq 5$. This completes the proof.
\end{proof}

\begin{lemma}\label{lem:phi(m,n)}
If ${\displaystyle |\phi_{m',n'}|\leq \Big(\frac{1+\nu_0^2}{4}\Big)^{m'+n'-1} C^{m'+n'}}$ for some $C>0$ for all $(m',n')$ such that $0\leq m'\leq m$ and $0\leq n' \leq n $ except $(m',n')=(0,0)$ and $(m,n)$, then 
\[
|\phi_{m,n}| \leq \Big(\frac{1+\nu_0^2}{4}\Big)^{m+n-1}C^{m+n}.
\]
\end{lemma}

\begin{proof}
Recall \eqref{eqn:RecurRightFirst}, \eqref{def:Omega}, and Lemma~\ref{lemma:Omega_bound2}, and we calculate
\begin{align*}
|\phi_{m,n}| \leq &\frac{1+\nu_0^2}{4 (mn + m+ n - 1)} \\ 
&\times \sum_{m'}\sum_{n'}\Big(\frac{1+\nu_0^2}{4}\Big)^{m'+n'-1}C^{m'+n'}
\Big(\frac{1+\nu_0^2}{4}\Big)^{m-m'+n-n'-1} C^{m-m'+n-n'}\\ 
= & \Big(\frac{1+\nu_0^2}{4}\Big)^{m+n-1}C^{m+n}.
\end{align*}
The assertion follows by induction.
\end{proof}

We utilize the result of Lemma~\ref{lem:phi(m,n)} to derive a rigorous bound of the form $|\phi_{m,n}|\leq C_0C^{m+n}$, and similar bounds for $\psi_{m,n}$, $v_{m,n}$, and $w_{m,n}$ via the following algorithm.

\begin{itemize}[itemindent=35pt]
\item[{\bf Step 1:}] Compute $\phi_{m,n}$ in $[0, M]\times [0, N]$, where $M$ and $N$ are at least $4$ (ensuring that all coefficients to be calculated have either $m$ or $n$ at least 5).

\item[{\bf Step 2:}] Determine the smallest value of $C_\phi$ such that $|\phi_{m,n}|\leq \big(\frac{1+\nu_0^2}{4}\big)^{m+n-1} C_\phi^{m+n}$ for all $(m,n)\in [0, M]\times [0, N]$. This involves finding $C_\phi = \max_{(m,n)} \left(\frac{4^{m+n-1}|\phi_{m,n}|}{(1+\nu_0^2)^{m+n-1}}\right)^{\frac{1}{m+n}}$.

\item[{\bf Step 3:}] By induction, $|\phi_{m,n}|\leq \big(\frac{1+\nu_0^2}{4}\big)^{m+n-1} C_\phi^{m+n}$ for all $m$ and $n$. 
    
\item[{\bf Step 4:}] Set $C = \big(\frac{1+\nu_0^2}{4}\big)C_\phi$ and $C_0 = \frac{4}{1+\nu_0^2}$ so that $|\phi_{m,n}|\leq C_0C^{m+n}$ for all $m$ and $n\geq 0$. 

\item[{\bf Step 5:}] Observe that $m\mu_-+n\mu_+ = -\frac{2(m+n)}{1+\nu_0^2} + \frac{2\nu_0(m-n)}{1+\nu_0^2}i$, whence $|\psi_{m,n}|\leq \frac{2}{\sqrt{1+\nu_0^2}}(m+n)C_0 C^{m+n}$. Set $C_\psi > C$ and $C_{\psi,0}=\frac{2}{\sqrt{1+\nu_0^2}}C_0 N(\frac{C}{C_\psi})^N$, so that $\frac{2}{\sqrt{1+\nu_0^2}}(m+n) C_0C^{m+n} \leq C_{\psi,0}C_\psi^{m+n}$ whenever $m+n\geq N$. Take $C_{\psi,0}$ larger if necessary to ensure $|\psi_{m,n}|\leq C_{\psi,0}C_\psi^{m+n}$ for all $0\leq m+n\leq N$. Thus, $|\psi_{m,n}|\leq C_{\psi,0}C_\psi^{m+n}$ for all $m$ and $n\geq 0$.

\item[{\bf Step 6:}] Consider 
\[
|v_{m,n}| \leq \frac{(1+\nu_0^2)^2}{8(m+n)^2}\sum_{m',n'=0}^{m,n} |v_{m',n'}||\psi_{m-m',n-n'}|.
\]
Choose $x_N>C_\psi$. For $N = m+n$,
\begin{align*}
|v_{m,n}| & \leq \frac{C_vx_N^N(1+\nu_0^2)^2}{8N^2}\sum_{m'=0,n'=0}^{m,n} \Big(\frac{C_\psi}{x_N}\Big)^{m'+n'}\\
&\leq \frac{C_vx_N^N(1+\nu_0^2)^2}{8N^2}\frac{1}{\big(1-\tfrac{C_\psi}{x_N}\big)^2}
\end{align*}
for some $C_v$. Take $C_v$ larger if necessary and choose $N$ and $x_N$ sufficiently large to ensure
$\frac{C_v(1+\nu_0^2)^2}{8N^2}\frac{1}{\big(1-\frac{C_\psi}{x_N}\big)^2}<1$.

\item[{\bf Step 7:}] Employ the same approach to establish a bound on $|w_{m,n}|$ as was done for $|v_{m,n}|$. Specifically, set $C_w>x_N$ and $C_{w,0}=KN(\frac{x_N}{C_w})^N$ for some $K>0$. Take $C_{w,0}$ larger if necessary to ensure $|w_{m,n}|\leq C_{w,0}C_w^{m+n}$ for all $0\leq m+n\leq N$. By induction, $|w_{m,n}|\leq C_{w,0}C_w^{m+n}$ for all $m,n\geq 0$.
\end{itemize}

\begin{tiny}
\begin{table}[!ht]
\centering
\begin{tabular}{|c|c||c|c|c|c|c|c|c|c|c|c|c|c||} 
\hline 
$\nu_L$ & $\nu_R$ & $\theta_i$ & $\theta_r$ & $R_\textrm{scale}$ & $N$ & $C_0$ & $C$ & $M$ & Time\\
\hline
0.28 & 0.281 & 3.62 & 0.0649 & 0.5 & 17 & 73.6 & 0.438 & 8 & 4.9e+00 \\ 
0.5 & 0.501 & 1.82 & 0.0825 & 0.5 & 17 & 55.1 & 0.585 & 6 & 1.3e+01 \\ 
1 & 1.01 & 4.86 & 0.0756 & 0.3 & 17 & 39 & 0.501 & 9 & 4.0e+00 \\ 
1.5 & 1.51 & 3.02 & 0.072 & 0.3 & 17 & 31.8 & 0.749 & 8 & 1.4e+01 \\ 
2 & 2.01 & 1.7 & 0.0745 & 0.3 & 17 & 27.6 & 0.997 & 8 & 1.5e+01 \\ 
2.5 & 2.51 & 0.671 & 0.0786 & 0.3 & 17 & 24.6 & 1.25 & 8 & 1.4e+01 \\ 
3 & 3.01 & 4.45 & 0.076 & 0.2 & 17 & 22.5 & 0.996 & 8 & 1.4e+01 \\ 
3.397 & 3.398 & 0.483 & 0.088 & 0.42 & 17 & 21.1 & 2.36 & 5 & 4.0e+00 \\ \hline 
\end{tabular}
\caption{Data concerning the rigorous enclosure of the two-dimensional stable solution manifold, about $x = \infty$, for \eqref{eq:the_ode}, the latter equation of \eqref{eq:bd}, and \eqref{eq:initial}, corresponding to the fixed point $(-1,0,1,0)$. We consider $\nu$ over the interval $[\nu_L,\nu_R]$ for various values of $\nu_L$ and $\nu_R$. Here $\theta_i$ and $\theta_r$ represent parameter values for the solution manifold as $\nu$ varies, and $R_\textrm{scale}$ denotes the scaling factor for the eigenvectors applied to initialize terms whose indices sum to one, $N$ indicates the index at which the infinite series are truncated, 
$C_0$ and $C$ are the constants for the dominating geometric series utilized to bound the truncation error, $M$ is the number of Chebyshev nodes in the parameter $\nu$ employed to approximate the coefficients of the stable manifold with polynomials, and `Time' is how long it took in seconds to compute the coefficients of the manifold and obtain the truncation error.}\label{table:2}
\end{table}
\end{tiny}

\noindent{\it Proof of Proposition~\ref{prop:middle sol}.}~(Construction of the middle solution manifold)~We assume
\[
\phi(x)=\sum_{n=0}^{\infty} \phi_n x^n,\quad 
\psi(x)=\sum_{n=0}^{\infty} \psi_n x^n,\quad 
v(x)=\sum_{n=0}^{\infty} v_n x^n,\quad 
w(x)=\sum_{n=0}^{\infty} w_n x^n,
\]
where $\phi_n$, $\psi_n$, $v_n$, $w_n$,  $n=0,1,2$,\ldots, are in \eqref{rec:middlefirst}-\eqref{rec:middlelast}. We establish bounds for $\phi_n$, $\psi_n$, $v_n$, $w_n$, and their first and second partial derivatives with respect to the first coefficients, $\phi_0$, $\psi_0$, $v_0$, $w_0$, along with the remainders and their first and second partial derivatives with respect to the first coefficients. 

Note that bounds for $\psi_n$ and $w_n$ can be derived from those for $\phi_n$ and $v_n$, whence it suffices to consider 
\begin{align*}
\phi_{n+2}=&\frac{1}{\nu(n+2)}\Big( -\phi_{n+1}+\frac{1}{2(n+1)}\sum_{n'=0}^n \phi_{n-n'}\phi_{n'}\Big)
\intertext{and}
v_{n+2}=&\frac{1}{2(n+1)(n+2)}\sum_{n'=0}^{n}(n-n'+1)\phi_{n-n'+1}v_{n'}.
\end{align*}
Although we may begin induction arguments with just two terms of recurrence relations, we consider situations where we know the first $m+2$ terms and derive bounds for all terms via induction.

\begin{lemma}\label{lemma_finite0}
Let $N_0\geq 0$ fixed, and  
\[
|\phi_n|\leq 2(N_0+1)\nu^{-n-1}\alpha^{n+2}\quad \text{for all $n\geq 0$},
\]
where 
\[
\alpha=\max\Big(\max_{0\leq n\leq N_0+1}\Big|\frac{\nu^{n+1}\phi_n}{2(N_0+1)}\Big|^{\frac{1}{n+2}},1\Big).
\]
Additionally, 
\[
|v_n|\leq C_0(\alpha\beta)^n\nu^{-n}\quad\text{for all $n\geq 0$},
\]
where $\beta>1$ is the root of $\alpha=\beta(\beta-1)$ and $C_0=\max_{0\leq n \leq N_0} |v_n|\nu^n(\alpha\beta)^{-n}$.
\end{lemma}

\begin{proof}
We define $\varphi_n=|\phi_n|$ for $0\leq n\leq N_0+1$ and 
\[
\varphi_{n}=\frac{1}{\nu n}\Big(\varphi_{n-1}+\frac{1}{2(N_0+1)}\sum_{n'=0}^{n-2}\varphi_{n-2-n'}\varphi_{n'}\Big)
\quad\text{for  $n\geq N_0+2$}.
\]
It follows from an induction argument that $|\phi_n|\leq \varphi_n$ for all $n\geq 0$, thus it suffices to prove the assertion for $\varphi_n$. Making the change of variables $\varphi_n = 2(N_0+1) \nu^{-n-1} \xi_n$, 
\[
\xi_{n}=\begin{cases}
{\displaystyle \frac{\nu^{n+1}}{2(N_0+1)}\varphi_n}&\text{for $0\leq n\leq N_0+1$},\\
{\displaystyle \frac{1}{n}\Big(\xi_{n-1}+\sum_{n'=0}^{n-2} \xi_{n-2-n'}\xi_{n'}\Big)} & \text{for $n \geq N_0+2$}.
\end{cases}
\]
For $0\leq n\leq N_0+1$, clearly, $|\xi_n| \leq \alpha^{n+2}$. Recalling the recurrence relation for $\xi_n$, we show by induction that
\begin{align*}
|\xi_{n}| \leq & \frac{1}{n}\Big(\alpha^{n+1} +\sum_{n'=0}^{n-2} \alpha^{n-n'} \alpha^{n'+2} \Big) \\
= & \frac{1}{n}(\alpha^{n+1} + (n-1) \alpha^{n+2}) \leq \alpha^{n+2}
\end{align*}
for all $n\geq N_0+2$.
Similarly, for $0\leq n\leq N_0$, clearly, $|v_n|\leq C_0(\alpha\beta)^n\nu^{-n}$. For $n\geq N_0+1$, we use the recurrence relation to show by induction that 
\begin{align*}
|v_n|&=\frac{1}{2n(n-1)}\Big|\sum^{n-2}_{n'=0}(n-n'-1)\phi_{n-n'-1}v_{n'}\Big|\\
&\leq \frac{1}{2n(n-1)}\sum^{n-2}_{n'=0}(n-n'-1)|\phi_{n-n'-1}|C_0(\alpha\beta)^{n'}\nu^{-n'}\\
&\leq \frac{1}{2n(n-1)}\sum^{n-2}_{n'=0}(n-n'-1)2(N_0+1)\nu^{n'-n}\alpha^{n-n'+1}C_0(\alpha\beta)^{n'}\nu^{-n'}\\
&=\frac{1}{2n(n-1)}\sum^{n-2}_{n'=0}(n-n'-1)2(N_0+1)\nu^{-n}\alpha^{n+1}C_0\beta^{n'}\\
&\leq C_0\frac{(N_0+1)\alpha}{n}\alpha^{n}\nu^{-n}\sum^{n-2}_{n'=0}\beta^{n'}
\leq C_0\frac{(N_0+1)\alpha}{n}\alpha^{n}\nu^{-n}\frac{\beta^{n-1}}{\beta-1}\\
&=C_0\frac{(N_0+1)\beta(\beta-1)}{n}\alpha^{n}\nu^{-n}\frac{\beta^{n-1}}{\beta-1}=C_0\frac{N_0+1}{n}\alpha^n\nu^{-n}\beta^n\\
&\leq C_0(\alpha\beta)^n\nu^{-n}.
\end{align*}
 This completes the proof.
\end{proof}

We now examine the first-order partial derivatives of $\phi_n$, $\psi_n$, $v_n$, and $w_n$ with respect to the first coefficients, $\phi_0$, $\psi_0$, $v_0$, and $w_0$, that is, 
\[
\mat{\partial_k\phi_{n}\\ \partial_k \psi_{n}\\ \partial_k v_{n}\\ \partial_k w_{n}}
=\mat{\frac{1}{n} \partial_k \psi_{n-1}\\\frac{1}{\nu n}(-\partial_k \psi_{n-1}+\sum_{n'=0}^{n-1}\phi_{n-1-n'}\partial_k \phi_{n'})\\
\frac{1}{n}\partial_k w_{n-1}\\ \frac{1}{2n}\sum_{n'=0}^{n-1}(v_{n'} \partial_k \psi_{n-1-n'}+\psi_{n-1-n'}\partial_k v_{n'})},
\]
where $\partial_k$, $k=1,2,3,4$, stand for $\frac{\partial}{\partial\phi_0}$, $\frac{\partial}{\partial \psi_0}$, $\frac{\partial}{\partial v_0}$, or $\frac{\partial}{\partial w_0}$.

\begin{lemma}\label{lemma_finite1}
For $N_0\geq0$ as in Lemma~\ref{lemma_finite0}, and
\[
|\partial_k\phi_n|\leq C_{1k}\alpha^n\nu^{-n} \quad \text{for all $n\geq 0$},
\]
$k=1,2,3,4$, where $C_{1k}=\max_{0\leq n\leq N_1}|\partial_k\phi_n|\nu^n\alpha^{-n}$, $N_1=\lceil\frac{1}{\alpha}+2N_0+1\rceil=2N_0+2$, and $\alpha$ is in Lemma~\ref{lemma_finite0}. Additionally, 
\[
|\partial_k v_n|\leq C_{2k}(\alpha\beta)^n\nu^{-n} \quad \text{for all $n\geq 0$},
\]
where $C_{2k}=\max\Big(\max_{0\leq n\leq N_2}|\partial_k v_n|\nu^n(\alpha\beta)^{-n},\frac{C_0C_{1k}\nu}{2\alpha^2}\Big)$, $N_2=N_0+1$, and $C_0$ and $\beta$ are in Lemma~\ref{lemma_finite0}.
\end{lemma}

\begin{proof}
For $0\leq n\leq N_1$, clearly, $|\partial_k\phi_n|\le C_{1k}\alpha^n\nu^{-n}$. For $n\geq N_1+1$, we use 
\[
\partial_k\phi_{n}=\frac{1}{\nu n(n-1)}\Big(-(n-1)\partial_k\phi_{n-1}+\sum_{n'=0}^{n-2}\phi_{n-2-n'}\partial_k\phi_{n'}\Big)
\]
and show by induction that 
\begin{align*}
|\partial_k\phi_{n}|\leq &\frac{1}{\nu n(n-1)}\Big((n-1)C_{1k}\frac{\alpha^{n-1}}{\nu^{n-1}}+\sum_{n'=0}^{n-2}2(N_0+1)\nu^{2-n+n'-1}\alpha^{n-2-n'+2}C_{1k}\alpha^{n'}\nu^{-n'}\Big)\\
=&C_{1k}\alpha^{n}\nu^{-n}\Big(\frac{1}{\alpha }+2N_0+2\Big)\frac{1}{n} \\
\leq & C_{1k}\alpha^n\nu^{-n}.
\end{align*}

Similarly, for $0\leq n\leq N_2$, clearly, $|\partial_k v_n|\leq C_{2k}(\alpha\beta)^n\nu^{-n}$. For $n\geq N_2+1$, we use 
\[
\partial_k v_{n}=\frac{1}{ 2n(n-1)}\sum_{n'=0}^{n-2}\Big(v_{n'}\partial_k \psi_{n-2-n'}+ \psi_{n-2-n'}\partial_k v_{n'}\Big)
\]
and show by induction that 
\begin{align*}
|\partial_k v_{n}|\leq & \frac{1}{2n(n-1)}\sum_{n'=0}^{n-2}(C_{0}\alpha^{n'}\beta^{n'}\nu^{-n'}(n-1-n')C_{1k}\nu^{-n+1+n'}\alpha^{n-1-n'}\\
&\qquad\qquad\qquad\quad+(n-1-n')2(N_0+1)\nu^{-n+1+n'-1}\alpha^{n-1-n'+2}C_{2k}\alpha^{n'}\beta^{n'}\nu^{-n'})\\
=&\frac{1}{2n(n-1)}\sum_{n'=0}^{n-2}(n-1-n')(C_{0}C_{1k}\alpha^{n-1}\beta^{n'}\nu^{-n+1}+2(N_0+1)\nu^{-n}\alpha^{n+1}C_{2k}\beta^{n'})\\
\leq&\frac{1}{2n}(C_{0}C_{1k}\alpha^{n-1}\nu^{-n+1}+2(N_0+1)\nu^{-n}\alpha^{n+1}C_{2k})\sum_{n'=0}^{n-2}\beta^{n'}\\
\leq&\frac{1}{2n}(C_{0}C_{1k}\alpha^{n-1}\nu^{-n+1}+2(N_0+1)\nu^{-n}\alpha^{n+1}C_{2k})\frac{\beta^{n-1}}{\beta-1}\\
=&\frac{1}{2n}(C_{0}C_{1k}\alpha^{-2}\nu+2(N_0+1)C_{2k})(\alpha\beta)^n\nu^{-n} \\
\leq& C_{2k}(\alpha\beta)^n\nu^{-n}.
\end{align*}
This completes the proof.
\end{proof}

Moving forward, consider the second-order partial derivatives of $\phi_n$, $\psi_n$, $v_n$, and $w_n$ with respect to the first coefficients $\phi_0$, $\psi_0$, $v_0$, and $w_0$, that is,
\[
\mat{\partial_{k_1}\partial_{k_2}\phi_{n}\\ \partial_{k_1}\partial_{k_2} \psi_{n}\\ \partial_{k_1}\partial_{k_2} v_{n}\\ \partial_{k_1}\partial_{k_2} w_{n}}=
\mat{\frac{1}{n} \partial_{k_1}\partial_{k_2} \psi_{n-1}\\ 
\begin{aligned}\tfrac{1}{\nu n}(&-\partial_{k_1}\partial_{k_2} \psi_{n-1} \\ &{\textstyle+\sum_{n'=0}^{n-1}(\partial_{k_2}\phi_{n-1-n'}\partial_{k_1} \phi_{n'}+\phi_{n-1-n'}\partial_{k_1} \partial_{k_2}\phi_{n'})})\end{aligned}\\
\frac{1}{n}\partial_{k_1}\partial_{k_2} w_{n-1}\\
\begin{aligned}{\textstyle \frac{1}{2n}\sum_{n'=0}^{n-1}}(&\partial_{k_2} v_k \partial_{k_1} \psi_{n-1-n'}+\partial_{k_1} v_{n'} \partial_{k_2} \psi_{n-1-n'} \\ &+ v_{n'} \partial_{k_1}\partial_{k_2} \psi_{n-1-n'}+\psi_{n-1-n'} \partial_{k_1}\partial_{k_2} v_{n'})\end{aligned}},
\]
where $\partial_{k_1}$ and $\partial_{k_2}$, $k_1, k_2=1,2,3,4$, stand for $\frac{\partial}{\partial\phi_0}$, $\frac{\partial}{\partial \psi_0}$, $\frac{\partial}{\partial v_0}$, or $\frac{\partial}{\partial w_0}$.

\begin{lemma}\label{lemma_finite2}
Let $N_0\geq0$ as in Lemma~\ref{lemma_finite0}, and
\[
|\partial_{k_1}\partial_{k_2}\phi_n|\leq C_{1k_1k_2}\alpha^n\nu^{-n}
\quad\text{for all $n\geq 0$},
\]
$k_1, k_2=1,2,3,4$, where $C_{1k_1k_2}=\max\Big(\max_{0\leq n\leq N_3}|\partial_{k_1}\partial_{k_2}\phi_n|\nu^n\alpha^{-n},\frac{C_{1k_1}C_{1k_2}\nu}{\alpha(2\alpha-1)}\Big)$, $N_3=\lceil\frac{1}{\alpha}+2N_0+2\rceil=2N_0+3$, $\alpha$ is in Lemma~\ref{lemma_finite0}, and $C_{1k_1}$ and $C_{1k_2}$ are in Lemma~\ref{lemma_finite1}. Similarly,
\[
|\partial_{k_1}\partial_{k_2}v_n|\leq C_{2k_1k_2}(\alpha\beta)^n\nu^{-n}\quad\text{for all $n\geq 0$},
\]
where 
\[
C_{2k_1k_2}=\max\Big(\max_{0\leq n\leq N_2}|\partial_{k_1} \partial_{k_2} v_n|\nu^n(\alpha\beta)^{-n},\frac{C_{1k_1}C_{2k_2}+C_{2k_1}C_{1k_2}+C_0C_{1k_1k_2}}{2\alpha^2}\Big),
\] $N_2=N_0+1$, and $C_0$ and $\beta$ are in Lemma~\ref{lemma_finite0}.
\end{lemma}

\begin{proof}
For $0\leq n\leq N_3$, clearly, $|\partial_{k_1}\partial_{k_2}\phi_n|\leq C_{1k_1k_2}\alpha^n\nu^{-n}$. For $n\geq N_3+1$, we use  
\begin{align*}
\partial_{k_1}\partial_{k_2}\phi_{n}=\frac{1}{\nu n(n-1)}\Big(&-(n-1)\partial_{k_1}\partial_{k_2}\phi_{n-1} \\ 
&+\sum_{n'=0}^{n-2}(\partial_{k_2}\phi_{n-2-n'}\partial_{k_1}\phi_{n'}+\phi_{n-2-n'}\partial_{k_1}\partial_{k_2}\phi_{n'})\Big)
\end{align*}
to show by induction that 
\begin{align*}
|\partial_{k_1}\partial_{k_2}\phi_{n}|\leq & \frac{1}{\nu n(n-1)}\Big((n-1)C_{1k_1k_2}\frac{\alpha^{n-1}}{\nu^{n-1}} \\ &\qquad\qquad\quad+\sum_{n'=0}^{n-2}(C_{1k_2}\alpha^{n-2-n'}\nu^{-n+2+n'}C_{1k_1}\alpha^{n'}\nu^{-n'} \\
&\qquad\qquad\qquad\qquad+2(N_0+1)\alpha^{n-2-n'+2}\nu^{-n+2+n'-1}C_{1k_1k_2}\alpha^{n'}\nu^{-n'})\Big)\\
=&\alpha^{n}\nu^{-n}\Big(\frac{C_{1k_1k_2}}{n\alpha}+\frac{2(N_0+1)C_{1k_1k_2}}{n}+\frac{C_{1k_1}C_{1k_2}\nu}{n\alpha^2}\Big)\\ 
\leq & C_{1k_1k_2}\alpha^n\nu^{-n}.
\end{align*}

Similarly, for $0\leq n\leq N_2$, clearly, $|\partial_{k_1}\partial_{k_2} v_n|\le C_{2k_1k_2}(\alpha\beta)^n\nu^{-n}$. For $n\geq N_2+1$, we use 
\begin{align*}
\partial_{k_1}\partial_{k_2} v_{n}=\frac{1}{ 2n(n-1)}\sum_{n'=0}^{n-2}(&\partial_{k_2}v_{n'}\partial_{k_1} \psi_{n-2-n'}
+\partial_{k_1}v_{n'}\partial_{k_2} \psi_{n-2-n'}\\
&+v_{n'}\partial_{k_1}\partial_{k_2} \psi_{n-2-n'}
+\psi_{n-2-n'}\partial_{k_1}\partial_{k_2} v_{n'})
\end{align*}
to show by induction that 
\begin{align*}
|\partial_{k_1}\partial_{k_2} v_{n}|\leq & 
\frac{1}{2n(n-1)}\big(C_{2k_2}C_{1k_1}\alpha^{n-1}\nu^{-n+1} \\
&\qquad\qquad\quad+C_{2k_1}C_{1k_2}\alpha^{n-1}\nu^{-n+1}+C_0C_{1k_1k_2}\alpha^{n-1}\nu^{-n+1}\\
&\qquad\qquad\qquad\qquad+2(N_0+1)C_{2k_1k_2}\nu^{-n}\alpha^{n+1}\big)\sum_{n'=0}^{n-2}(n-1-n')\beta^{n'}\\
\leq&\frac{\alpha^n\beta^n\nu^{-n}}{2n}\Big(\frac{\nu}{\alpha^2}(C_{1k_1}C_{2k_2}+C_{2k_1}C_{1k_2}+C_{0}C_{1k_1k_2})+2(N_0+1)C_{2k_1k_2}\Big)\\
\leq & C_{2k_1k_2}(\alpha\beta)^n\nu^{-n}.
\end{align*}
This completes the proof.
\end{proof}

Lastly, we address the remainders and their first and second partial derivatives with respect to the first coefficients. 

\begin{proposition}
For $x$ in a bounded interval of $(-\infty,\infty)$, let $x_1=|\alpha x\nu^{-1}|$ and $x_2=|\alpha \beta x\nu^{-1}|$, where $\alpha$ and $\beta$ are in Lemma~\ref{lemma_finite0}. For $x_2<1$, 
\begin{equation}\begin{aligned}\label{truncation_bounds}
&\begin{aligned}\frac{\nu}{2(N_0+1)\alpha^2}\Big|\sum_{n=N}^{\infty}\phi_nx^n\Big|,
&\frac{1}{C_{1k}}\Big|\sum_{n=N}^{\infty}\partial_k\phi_nx^n\Big|,\\
&\frac{1}{C_{1k_1k_2}}\Big|\sum_{n=N}^{\infty}\partial_{k_1}\partial_{k_2}\phi_nx^n\Big|\leq \frac{x_1^N}{1-x_1},\end{aligned}\\
&\begin{aligned}\frac{\nu^2}{2(N_0+1)\alpha^3}\Big|\sum_{n=N}^{\infty}\psi_nx^n\Big|,
&\frac{\nu}{C_{1k}\alpha}\Big|\sum_{n=N}^{\infty}\partial_k\psi_nx^n\Big|,\\ 
&\frac{\nu}{C_{1k_1k_2}\alpha}\Big|\sum_{n=N}^{\infty}\partial_{k_1}\partial_{k_2}\psi_nx^n\Big|\leq\frac{(N+1)x_1^N-Nx_1^{N+1}}{(1-x_1)^2},\end{aligned}\\ 
&\frac{1}{C_0}|\sum_{n=N}^{\infty}v_nx^n|,
\frac{1}{C_{2k}}|\sum_{n=N}^{\infty}\partial_kv_nx^n|,
\frac{1}{C_{2k_1k_2}}|\sum_{n=N}^{\infty}\partial_{k_1}\partial_{k_2}v_nx^n|\leq \frac{x_2^N}{1-x_2},\\
&\begin{aligned}\frac{\nu}{C_0\alpha\beta}|\sum_{n=N}^{\infty}w_nx^n|,
&\frac{\nu}{C_{2k}\alpha\beta}|\sum_{n=N}^{\infty}\partial_kw_nx^n|,\\
&\frac{\nu}{C_{2k_1k_2}\alpha\beta}|\sum_{n=N}^{\infty}\partial_{k_1}\partial_{k_2}w_nx^n|\leq \frac{(N+1)x_2^N-Nx_2^{N+1}}{(1-x_2)^2},\end{aligned}
\end{aligned}
\end{equation}
where $N_0$, $C_0$, $C_{1k}$, $C_{2k}$, $C_{1k_1k_2}$, and $C_{2k_1k_2}$ are defined in Lemmas \ref{lemma_finite0}-\ref{lemma_finite2}.
\end{proposition}

The proof follows from direct calculations. We omit the details.

\begin{remark}\label{choiceofparameter}\rm
In practice, setting $\alpha=1$ and $\beta=\frac{1+\sqrt{5}}{2}$ has proven effective. It is noteworthy that to ensure $x_2<1$, the maximum step size allowed is $\frac{2\nu}{1+\sqrt{5}}$. Interestingly, $\nu$ appears in the numerator, rather than denominator, so that as $\nu$ increases, the step size can be taken larger. This might seem counterintuitive, but it aligns well with the observation that the width of the transition region expands with increasing $\nu$, and is helpful for validated numerics for $\nu$ large. Consequently, as $\nu$ increases, the length of the interval $[0,x_N]$ will have to increase to ensure that the second part of Lemma~\ref{lemma:Calogero} holds. This also aligns with the widening of the transition region with increasing $\nu$.

This necessitates larger step size to control the dimension of \eqref{def:boldF} as well as $\mathbf{A}$ and $\mathbf{A}^\dagger$ involved in the  Newton--Kantorovich argument. 
However, this results in slower convergence of the tails in \eqref{truncation_bounds}, necessitating an increase in $N$ to manage truncation error. Computing more coefficients can lead to a worsening of the wrapping effect. Therefore, finding an optimal balance among the step size, $N$, and other parameters in both left and right solution manifolds is crucial for mitigating the wrapping effect while managing truncation and interpolation errors effectively.  

For $\nu>\nuup$, determining optimal parameters for which \eqref{def:p} holds true in the Newton--Kantorovich argument becomes challenging. Achieving this within the precision of \normalfont{INTLAB} proves cumbersome. It would be interesting to translate our package to Python, leveraging the arbitrary-precision interval arithmetic \normalfont{iv} available in \normalfont{mpmath}. This move would aim to push the upper end of range closer to $\nu_c\approx 4.09$.
\end{remark}

\begin{tiny}
\begin{table}[!ht]
\centering
\begin{tabular}{|c|c||c|c|c|c|c|c|c|c|} 
 \hline 
 $\nu_L$ & $\nu_R$ & $\Delta x$ & $\alpha$ & $\beta$ & $C_{1}$ & $C_{2}$ & $C_{2,1}$ & $C_{2,2}$ & Time \\
 \hline
0.28 & 0.281 & 5.0e-02 & 1.0e+00 & 1.6e+00 & 1.0e+00 & 1.0e+00 & 2.8e-01 & 5.0e-01 & 8.4e+01 \\ 
0.5 & 0.501 & 5.0e-02 & 1.0e+00 & 1.6e+00 & 1.0e+00 & 1.0e+00 & 5.0e-01 & 5.1e-01 & 5.6e+01 \\  
1 & 1.01 & 2.0e-01 & 1.0e+00 & 1.6e+00 & 1.0e+00 & 1.0e+00 & 1.0e+00 & 1.0e+00 & 6.0e+01 \\ 
1.5 & 1.51 & 2.0e-01 & 1.0e+00 & 1.6e+00 & 1.5e+00 & 1.1e+00 & 3.4e+00 & 3.4e+00 & 5.8e+01 \\ 
2 & 2.01 & 2.0e-01 & 1.0e+00 & 1.6e+00 & 2.0e+00 & 2.1e+00 & 8.1e+00 & 8.3e+00 & 5.8e+01 \\  
2.5 & 2.51 & 2.0e-01 & 1.0e+00 & 1.6e+00 & 2.5e+00 & 3.3e+00 & 1.6e+01 & 1.6e+01 & 5.7e+01 \\ 
3 & 3.01 & 2.0e-01 & 1.0e+00 & 1.6e+00 & 3.0e+00 & 4.4e+00 & 2.7e+01 & 2.7e+01 & 5.5e+01 \\ 
3.397 & 3.398 & 2.0e-01 & 1.0e+00 & 1.6e+00 & 3.4e+00 & 6.0e+00 & 3.9e+01 & 4.1e+01 & 1.6e+01 \\ \hline 
\end{tabular}
\caption{Data concerning the rigorous enclosure of solutions of \eqref{eq:the_ode} over bounded intervals, for $\nu$ within the interval $[\nu_L,\nu_R]$ for various values of $\nu_L$ and $\nu_R$. Here $\Delta x$ means the step size, $\alpha$ and $\beta$ are in Lemma~\ref{lemma_finite0},  $C_{1}$, $C_{2}$ are the maxima of those in Lemma~\ref{lemma_finite1}, and $C_{2,1}$ and $C_{2,2}$ are the maxima of those in Lemma \ref{lemma_finite2}, and we record values for the rightmost subinterval.
`Time' is how long it took in seconds to compute the coefficients, values for the rightmost subinterval.}\label{table:3}
\end{table}
\end{tiny}


\subsection{Applying the Newton--Kantorovich theorem} 

Applying Theorem~\ref{thm:NewtK} is a pivotal step in our CAP. It rigorously validates the existence and uniqueness of an exact solution in the vicinity of an approximate solution. Consider $\mathbf{F}:\mathbb{R}^{4N+6} \to \mathbb{R}^{4N+6}$, defined in \eqref{def:boldF}. 
Suppose $(\phi_{j,\rm ap}, \psi_{j,\rm ap}, v_{j,\rm ap}, w_{j,\rm ap})$ make numerical approximate solutions of \eqref{eq:the_ode} at $x_j$, $j=0,1,\ldots, N$, and are extended to $(-\infty,\infty)$ via $\mathbf{\Phi}$, $\mathbf{u}_L$, and $\mathbf{u}_R$, along with $\theta_{\rm ap}$ and $\theta_{i,\rm ap}$. Abusing notation slightly, let 
\[
y_{\rm ap}=(\theta_{\rm ap}, \{(\phi_{j,\rm ap}, \psi_{j,\rm ap}, v_{j,\rm ap}, w_{j, \rm ap})\}_{j=0}^N, \theta_{i, \rm ap}).
\]
We take $\mathbf{A}^\dagger$ to be the numerically computed Jacobian of $\mathbf{F}$, and $\mathbf{A}$ the numerically computed inverse. We verify (A1) through (A5) of Theorem~\ref{thm:NewtK}.

To ensure the operator norm $\|\cdot\|_{B(\mathbb{R}^{4N+6})}$ in (A3), (A4), and (A5) of Theorem~\ref{thm:NewtK} is easily computable, we take the Euclidean norm $\|\cdot\|_2$ on $\mathbb{R}^{4N+6}$ and the induced $2$-norm $\|\cdot\|_2$ on $B(\mathbb{R}^{4N+6})$. It is computationally simpler and preferable to work with the Hilbert-Schmidt norm, also known as the Frobenius norm, in Euclidean space. Since
\[
\|\mathbf{A}\|_{2}\leq \|\mathbf{A}\|_F\leq (4N+6)\|\mathbf{A}\|_\infty \quad\text{for any $\mathbf{A}\in B(\mathbb{R}^{4N+6})$},
\]
where $\|\cdot\|_F$ denotes the sub-multiplicative Frobenius norm and $\|\cdot\|_\infty$ the maximum norm, it follows that 
\[
\|\mathbf{1} -\mathbf{A}\mathbf{A}^\dagger\|_2 \leq\|\mathbf{1} -\mathbf{A}\mathbf{A}^\dagger\|_F\quad\text{and}\quad 
\|\mathbf{A}(\mathbf{A}^\dagger -D\mathbf{F}(y_{\rm ap}))\|_2\leq\|\mathbf{A}\|_F\|(\mathbf{A}^\dagger -D\mathbf{F}(y_{\rm ap}))\|_F.
\]
Moreover, if $\mathbf{F}$ is twice continuously differentiable, then for $y\in \mathbb{R}^{4N+6}$ with $\|y-y_{\rm ap}\|\leq r$ for some $r>0$, it follows that
\begin{align*}
\|\mathbf{A}&(D\mathbf{F}(y_{\rm ap})- D\mathbf{F}(y))\|_{2}\\
\leq &\|\mathbf{A}(D\mathbf{F}(y_{\rm ap})-D\mathbf{F}(y))\|_F
\leq \|\mathbf{A}\|_F\|D\mathbf{F}(y_{\rm ap})-D\mathbf{F}(y)\|_F \\
\leq & (4N+6)\|\mathbf{A}\|_F\|D\mathbf{F}(y_{\rm ap})-D\mathbf{F}(y)\|_\infty \\
\leq &  (4N+6)\|\mathbf{A}\|_F \max_{1\leq j,k\leq (4N+6)}\sup_{\|y_0-\frac12(y_{\rm ap}+y)\|\leq \frac12\|y_{\rm ap}-y\|}\| D(D\mathbf{F})_{jk}(y_0)\|_F\|y_{\rm ap}-y\| \\
\leq &  (4N+6)r \|\mathbf{A}\|_F \max_{1\leq j,k \leq 4N+6} \sup_{\|y-y_{\rm ap}\|\leq r}\|D(D\mathbf{F})_{jk}(y)\|_F,
\end{align*}
so that we can choose
\begin{align*}
&Z_0=\|\mathbf{A}(D\mathbf{F}(y_{\rm ap})-D\mathbf{F}(y))\|_F,\\ 
&Z_1=\|\mathbf{A}\|_F\|(\mathbf{A}^\dagger -D\mathbf{F}(y_{\rm ap})\|_F,
\intertext{and}
&Z_2(r)= (4N+6)\|\mathbf{A}\|_F \max_{1\leq j,k\leq 4N+6} \sup_{\|y-y_{\rm ap}\|\leq r}\|D(D\mathbf{F})_{jk}(y)\|_F.
\end{align*}
We observe that $D(D\mathbf{F})_{jk}(y)$ contains at most four non-zero entries, so that the Frobenius norm is the Euclidean $2$-norm. Through validated numerics or rigorous computation, we establish a rough yet uniform bound $R$ on the modulus of the entries of $D^2\mathbf{F}(y)$, where $y$ varies within a ball centered at $y_{\rm ap}$ with radius $r$. Consequently, 
\[
\sup_{\|y-y_{\rm ap}\|\leq r}\|D(D\mathbf{F})_{jk}(y)\|_2\leq \sqrt{4R^2}=2R.
\]
To establish a rough bound on the entries of $D^2\mathbf{F}$, we employ the series representations of solutions, detailed throughout the manuscript, initializing with the interval given by the ball centered at $y_{\rm ap}$ with radius $r$. This allows us to derive bounds for these quantities over the ball using interval arithmetic.

One can a priori apply the Newton--Kantorovich theorem when the interval in which \eqref{def:p} holds arbitrarily far out along the positive $r$ axis. However, for practical considerations, we fix some $r_0$ and evaluate 
\[
Z_2(r_0) r_0 = \sup_{\|y-y_{\rm ap}\|\leq r_0} \| \mathbf{A}(D\mathbf{F}(y_{\rm ap}) - D\mathbf{F}(y))\|,
\]
requiring that there exists a subinterval $[r_1,r_2] \subset (0,r_0)$ where \eqref{def:p} holds to proceed. 

Once (A1)-(A5) in Theorem~\ref{thm:NewtK} hold, along with other conditions, there exist an exact solution  
\[
y_{\rm ex}=(\theta_{\rm ex}, \{({\phi}_{j,\rm ex}, {\psi}_{j,\rm ex},{v}_{j,\rm ex},{w}_{j,\rm ex})\}_{j=0}^{N},\theta_{i,\rm ex}),\quad
\|y_{\rm ex}-y_{\rm ap}\|<r_1,
\]
to the connection problem. Clearly, 
\[
|{\phi}_{j,\rm ex}-\phi_{j,\rm ap}|, |{\psi}_{j,\rm ex}-\psi_{j,\rm ap}|, |v_{j,\rm ex}-v_{j,\rm ap}|, |w_{j,\rm ex}-w_{j,\rm ap}|<r_1.
\] 

Therefore, at each $x_j$, the exact solution $({\phi}_{j,\rm ex}, {\psi}_{j,\rm ex}, {v}_{j,\rm ex}, {w}_{j,\rm ex})$ lies within an interval, accompanied by similar control of $\theta_{\rm ex}$ and $\theta_{i,\rm ex}$. With additional control over $\mathbf{\Phi}$, $\mathbf{u}_R$, and $\mathbf{u}_L$, they extend to provide bounds on $({\phi}_{\rm ex}, {\psi}_{\rm ex}, {v}_{\rm ex}, {w}_{\rm ex})$ over the intervals $(x_j,x_{j+1})$, $j=0,1,\ldots, N$, as well as $[x_N,\infty)$ and $(-\infty,0],$ with some increase in $r_1$. This rigorous control allows us to validate the exact solution across the entire real line. 

Since $\|y_{\rm ex}-y_{\rm ap}\|_{\infty}\leq \|y_{\rm ex}-y_{\rm ap}\|_{2} $, once we validate 
\[
\|(\phi_{\rm ex}-\phi_{\rm ap}, \psi_{\rm ex}-\psi_{\rm ap}, v_{\rm ex}-v_{\rm ap}, w_{\rm ex}-w_{\rm ap})\|_{2}\leq r,
\]
we validate 
\[
|\phi_{\rm ex}-\phi_{\rm ap}|, 
|\psi_{\rm ex}-\psi_{\rm ap}|,
|v_{\rm ex}-v_{\rm ap}|, |w_{\rm ex}-w_{\rm ap}|\leq r,
\]
which proves useful in validating the signs of the entry-wise variables $\phi$, $\psi$, $v$, and $w$.

\begin{remark}\rm 
We pause here to explain how we obtain an approximate zero of \eqref{def:boldF} in order to apply the Newton--Kantorovich theorem. This is accomplished using standard non-rigorous numerical methods:

\begin{itemize}[itemindent=35pt]
\item[\bf Step 1:] Find an approximation to $\mathbf{u}_L(\theta;0)$ by summing a large number of terms of \eqref{eqn:LeftM} with $x=0$. We have provided typical values of the parameters in Table~\ref{table:1}. Typical values of $\theta$ are less than but of the order of $1$, and $N$ is approximately $60$.

\item[\bf Step 2:] Construct an approximation to $\mathbf{u}_R(\theta_i;x_N,\theta_r)$ by summing a large number of terms of \eqref{eqn:RightM} with $x=x_N$. 

\item[\bf Step 3:] Numerically evolve the first two equations of \eqref{eq:the_ode} from $x=0$ to $x_N$ with the initial condition $\mathbf{u}_L(\theta;0)$ obtained from Step 1 using the MATLAB's `ode15s' algorithm. We denote the final state by $(\phi(x_N),\psi(x_N))$.
    
\item[\bf Step 4:] Use $(\phi(x_N),\psi(x_N))$ from Step 3 and a nonlinear root solver to find an approximate zero of 
\[
(\theta_{r,\rm ap},\theta_{i,\rm ap})\mapsto \mathbf{\Pi}\mathbf{u}_R(\theta_i;x_N,\theta_r)-(\phi(x_N),\psi(x_N)).
\]
This is performed in MATLAB using `fsolve' with the Levenberg--Marquardt algorithm.

\item[\bf Step 5:] The values of $\theta_{r,\rm ap}$, and $\theta_{i,\rm ap}$ from Step 4 provide approximate boundary conditions for $(\phi(x_N),\psi(x_N))$. We numerically solve the boundary value problem using MATLAB's `bvp5c' boundary value problem solver.  
\end{itemize}
\end{remark}

\begin{tiny}
\begin{table}[h!]
\centering
\begin{tabular}{|c|c||c|c|c|c|c|c|} 
 \hline 
 $\nu_L$ & $\nu_R$ & $B_2$ & $Y_0$ & $Z_0$ & $Z_1$ & $Z_2$& $r$\\
 \hline

0.28 & 0.281 & 4.9e-01 & 1.4e-08 & 1.2e-12 & 1.6e-02 & 9.8e+04 & 1.0e-06 \\ 
0.5 & 0.501 & 4.4e-01 & 8.7e-09 & 1.5e-13 & 3.7e-03 & 6.8e+04 & 1.0e-06 \\ 
1 & 1.01 & 4.6e-01 & 4.0e-07 & 2.7e-13 & 3.0e-02 & 6.1e+04 & 1.0e-05 \\ 
1.5 & 1.51 & 4.5e-01 & 2.3e-09 & 6.6e-13 & 2.2e-02 & 1.8e+05 & 4.8e-06 \\ 
2 & 2.01 & 4.5e-01 & 2.3e-09 & 1.3e-12 & 3.5e-02 & 4.2e+05 & 2.1e-06 \\ 
2.5 & 2.51 & 4.5e-01 & 2.2e-09 & 2.4e-12 & 4.5e-02 & 8.0e+05 & 1.1e-06 \\ 
3 & 3.01 & 4.5e-01 & 3.2e-08 & 6.3e-12 & 5.8e-02 & 2.2e+06 & 3.5e-07 \\ 
3.397 & 3.398 & 4.5e-01 & 2.3e-07 & 3.2e-12 & 4.0e-03 & 1.0e+06 & 5.3e-07 \\ 
\hline
\end{tabular}
\caption{Data regarding the bounds obtained in the Newton--Kantorovich argument for $\nu$ in $[\nu_L,\nu_R]$ for various values of $\nu_L$ and $\nu_R$. Here $B_2$ is a bound on the second derivative.}
\label{table:4}
\end{table}
\end{tiny}

\begin{remark}\rm
The rigorous enclosures of $\phi$ and $\psi$ over the interval $(-\infty,\infty)$, and $v$ and $w$ over $[0,\infty)$, enable us to determine the number of the zeros of $v$ over $[0,\infty)$. Since zeros of a solution of \eqref{eq:Riccati2} are necessarily simple, upon identifying an interval or intervals where $v$ has a zero, we only need to confirm that $w=v'$ is bounded away from zero to conclude that the zero is unique. It is noteworthy that confining the zero to one subinterval $[x_j,x_{j+1}]$ is not necessary to establish its uniqueness; several adjacent intervals where the bounds on $v$ straddle zero are acceptable as long as $w$ remains bounded away from zero on all of these intervals. 

The basic procedure involves: 
\begin{itemize}
\item Computing the right manifold using interval arithmetic to manage rounding errors, particularly ensuring the truncation error is sufficiently small to ascertain that $v$ has no zeros over $(x_N,\infty)$.

\item Identifying interval(s) $(x_j,x_{j'})$ where a sign change of $v$ must occur. Typically, there will be an interval $(x_{j-1},x_j)$ where the bounds guarantee that $v$ is positive, and another interval $(x_{j'},x_{j'+1})$ where the bounds guarantee that $v$ is negative. It then follows from continuity that there exists at least one zero in $(x_{j},x_{j'})$.

\item Verifying that the bounds on $w=v'$ are such that $w$ cannot vanish on $(x_j,x_{j'}).$ If true, there exists a unique zero of $v$ in $(x_j,x_{j'})$.
\end{itemize}
\end{remark}

\subsection{Riccati equation and no zeros over $(-\infty,0]$}

In the final stages of our CAP, it is essential to establish analytical conditions ensuring that   
\[
-\frac{d^2v}{dx^2}+\frac12 \phi_x(x)v=0, \qquad \lim_{x\rightarrow \infty}v(x)\rightarrow 1 
\]
will have no zeros over the interval $(-\infty,0]$. This can be conveniently done by considering the associated Riccati equation, where the zeros of $v$ correspond to the poles of the Riccati equation. Let
\[
\zeta=\frac{v_x}{v},
\]
and $\zeta$ solves 
\begin{equation}\label{eqn:w}
\frac{d\zeta}{dx}=-\zeta^2+\frac12\phi_x(x).
\end{equation}
Since $-\frac{d^2}{dx^2}+\frac12\phi_x(x)$ has at least one negative eigenvalue, $v$ has at least one zero. 

To ensure some sign condition crucial for our analysis, we rely on results from \cite{BS1985}.

\begin{lemma}\label{lemma:BonaSchonbek}
Let $\nu>\frac14$, and $\phi$ solve \eqref{eqn:phi(KdVB)}. Suppose $\phi(x_0)$, where $x_0\in\mathbb{R}$, represents the minimum of $\phi$ over $x\in\mathbb{R}$. Then
\begin{itemize}
\item $\phi(x)$ in monotonically decreasing for $x\in(-\infty,x_0)$;
\item There exists a unique $x_1<x_0$ such that $\phi_{xx}(x_1)=0$ and $\phi_{xx}(x)<0$ for $x\in (-\infty,x_1)$.
\item If $\phi(x)>-1$ for $x\in(-\infty,x_2)$ for some $x_2\in\mathbb{R}$, then $\phi(x)$ is monotonically decreasing over $x\in(-\infty,x_2)$. Additionally, if $\phi_{xx}(x_2)<0$ then $\phi_{xx}(x)<0$ for $x\in(-\infty,x_2)$.
\end{itemize}
\end{lemma}

The first two assertions are from \cite{BS1985}. Note that our notation slightly differs; the oscillations are around $\phi=1$, not $-1$, corresponding to $\nu<0$ in our notation. We have slightly adjusted the results of \cite{BS1985} to align with our conventions. 

Since the global minimum of $\phi$ is $<-1$, we deduce that if $\phi(x)>-1$ for $x\in(-\infty,x_2)$ then the interval $(-\infty,x_2)$ does not contain the global minimum. Consequently, $(-\infty,x_2)\subset (-\infty,x_0)$ and, thus, $\phi$ is monotonically decreasing over $(-\infty,x_2)$. If $\phi_{xx}(x_2)<0$ then we can further conclude that $\phi_{xx}(x_2)<0$ over the interval $(-\infty,0)$ because it follows from \cite{BS1985} that there exists a unique zero of $\phi_{xx}$ to the left of $x_0$.  

It is important to emphasize that we aim to establish $\phi(x)>-1$ across a semi-infinite interval, rather than merely at a point. 

\begin{lemma}
If ${\displaystyle \Large\left\vert\frac{\theta}{1-\frac{\theta}{(6-4\mu)}}\Large\right\vert<2}$ and $\phi_{xx}(0)<0$,
then $\phi_x(x)<0$ and $\phi_{xx}(x)<0$ for $x\in (-\infty,0)$.
\end{lemma}

\begin{proof}
Recall that
\[
\phi(x) = 1+ \sum_{n=1}^\infty \phi_n \theta^n e^{n \mu x}, 
\quad \text{for $x\in(-\infty,0)$},
\]
where $|\phi_n|<(6 - 4 \mu)^{-(n-1)}$. Bounding the residual using the geometric series and summing, we obtain
\[
\phi(x) \geq 1 - \frac{\theta e^{\mu x}}{1 - \theta e^{\mu x}(6 - 4 \mu)^{-1}} \geq 1-\frac{\theta}{1-\frac{\theta}{(6-4\mu)}},
\]
where the second inequality follows from monotonicity of $ae^x(1-a e^x)^{-1}$. If 
\[
\Big|\frac{\theta}{1-\frac{\theta}{(6-4\mu)}}\Big|<2
\]
then $\phi(x)>-1$ for $x\in(-\infty,0)$, implying that $\phi$ must be monotone and concave down over the interval $(-\infty,0)$ by Lemma~\ref{lemma:BonaSchonbek}. 
\end{proof}

\begin{lemma}\label{lemma:Calogero}
Assume that zero is sufficiently far to the left such that $\phi_x(x)<0$ and $\phi_{xx}(x)<0$ for $x\in(-\infty,0]$. 
\begin{itemize}
\item If $\zeta(0)>0$ then there must exist $x^*<0$ such that $v(x^*)=0$. In other words, there are at least two zeros of $v$ over the interval $(-\infty,\infty)$.
\item If $\zeta(0)<0$ and 
\[
\arctan(-\sqrt{\frac{2 }{-\phi_x(0)}}\zeta(0)) \geq \frac{1}{\sqrt{\mu}}\log(\frac{1+\sqrt{\frac{\theta} {2}}}{1-\sqrt{\frac{\theta} {2}}}),
\]
then there can be no zero of $v$ over $(-\infty,0]$. Thus, verifying the absence of a second zero reduces to the finite interval between $0$ and the first zero of $v$.
\end{itemize}
\end{lemma}

\begin{proof}
For the first assertion, consider an arbitrary $x_0<0$ such that the solution of \eqref{eqn:w} is continuous over the interval $[x_0,0]$. Comparing the solution with the solution of 
\begin{equation}\label{eqn:w0}
\frac{d\zeta_0}{dx}=-\zeta_0^2,\qquad \zeta_0(0)=\zeta(0),
\end{equation}
we obtain that $\zeta(x)>\zeta_0(x)$. Solving \eqref{eqn:w0} explicitly, $\zeta_0(x)=(x+\frac{1}{\zeta(0)})^{-1}$. Since $\zeta_0$ blows up for finite $x<0$, so must $\zeta$, ensuring the existence of a zero of $v$.

For the second assertion, note that since $v(x)\rightarrow 1$ as $x\rightarrow \infty$ and since $v$ has one simple zero over the interval $(0,\infty)$, it follows that $v(0)<0$. Since $\zeta(0)<0$, moreover, $v_x(0)>0$. To have a zero of $v$ over $(-\infty,0]$, we must have $-x_0\in (-\infty,0]$ (not the same as in the proof of the first assertion), where $v_x(-x_0)=0$. We wish to show that some inequality must hold; failure of the inequality provides sufficient for the non-existence of a zero over the interval $(-\infty,0)$.

Making the change of variables $\zeta(x)=g(x)\sqrt{-\frac12\phi_x(x)}$, 
\[
\frac{dg}{dx} = -\frac12 \frac{\phi_{xx}(x)}{\phi_x(x)} g -\frac{\sqrt{2}}{2}\sqrt{-\phi_x(x)}(g^2+1).
\]
The first term on the right side is positive over the interval $(-x_0,0)$, by assumption, while the second is negative, so that
\[
\frac{1}{1+g^2}\frac{dg}{dx} \geq -\sqrt{-\frac12\phi_x(x)}.
\]
Integrating from $x=-x_0$ such that $g(-x_0)=0$, we obtain
\[
-\arctan(g(0)) \leq \int_{-x_0}^{0} \sqrt{-\frac12\phi_x(x)}~dx 
\leq \int_{-\infty}^{0} \sqrt{-\frac12\phi_x(x)}~dx,
\]
whence 
\[
\arctan(-g(0)) > \int_{-\infty}^{0} \sqrt{-\frac12\phi_x(x)}~dx
\]
is a sufficient condition for the non-existence of a second zero. 

We give a direct estimation of this using the series representation of the left solution manifold. Recall from Proposition~\ref{prop:left sol} that
\eq{
\phi_x(x;\theta)= \psi(x;\theta) = \sum_{n=0}^{\infty} \mu n \phi_n(\mu) \theta^ne^{n\mu x}.
}{\notag}
Recall from Lemma~\ref{lemma:re_mu} that $|\phi_n(\mu)|<\frac{1}{2^{n-1}}$, so that
\[
|\phi_x(x)| \leq \sum_{n=1}^\infty n \mu \frac{\theta^n e^{n \mu x}}{2^{n-1}}
\leq \mu \theta e^{\mu x}\sum_{n=1}^\infty n \Big(\frac{\theta}{2}\Big)^{n-1} e^{(n-1) \mu x}.
\]
Since $\sum_{n=1}^\infty n r^{n-1} = \frac{1}{(1-r)^2}$,
\[
\sqrt{|\phi_x(x)|}\leq \frac{\sqrt{\theta\mu}e^{ \frac{\mu x}{2}}}{1-\frac{\theta}{2}e^{\mu x}},
\]
so that
\[
\int_{-\infty}^{0}\sqrt{|\phi_x(x)|}~dx  \leq \int_{-\infty}^{0}\frac{\sqrt{\theta\mu}e^{ \frac{\mu x}{2}}}{1-\frac{\theta}{2}e^{\mu x}}~dx.
\]
Therefore, if 
\[
\arctan(-\sqrt{\frac{2 }{-\phi_x(0)}}\zeta(0)) \geq \frac{1}{\sqrt{\mu}}\log(\frac{1+\sqrt{\frac{\theta}{2}}}{1-\sqrt{\frac{\theta}{2}}})
\]
then there does not exist a second zero of $v$ over the interval $(-\infty,0]$. 
\end{proof}

We pause to remark that Lemma~\ref{lemma:Calogero} resembles Calogero's bound on the number of negative eigenvalues for a Schr\"odinger operator with a monotone central potential \cite{Calogero2,Calogero,RS}.

\subsection{Rigorous verification over an interval}\label{sec:interval2}

We have outlined our CAP approach to derive rigorous error bounds for the numerical approximation of solutions to \eqref{eq:the_ode} for a single value of the parameter $\nu$, in practice, within an interval of machine epsilon width. Extending the methodology to cover intervals in $\nu$ is achieved through rigorous Chebyshev interpolation. The basic strategy follows that for the single parameter scenario.

Step 1 involves selecting a small interval $[\nu_0,\nu_{M+1}]$ 
and dividing it into Chebyshev nodes $\nu_1$, $\nu_2$, \ldots, $\nu_M$. We then compute numerical approximations to $\phi_j$, $\psi_j$, $v_j$, $w_j$, as well as the left and right solution manifold parameters at each $\nu_m$. Utilizing these values, we establish polynomial approximations for each of these quantities as a function of $\nu$. Since these serve as inputs to the Newton--Kantorovich argument, approximate values suffice, and we need not error computation in the Chebyshev approximation. 

Step 2 involves applying the Newton--Kantorovich theorem to obtain error bounds for the approximate solution. While one could in principle compute the solution operator $\mathbf{\Phi}$ via the series solution
\[
\mathbf{\Phi}(\mathbf{u}_j(\nu);x) = \sum_{k=0}^\infty (\phi_j^k(\nu),\psi_j^k(\nu),w_j^k(\nu),z_j^k(\nu))(x)
\]
using a polynomial initial condition for the recurrence, this has proven to be inefficient or insufficient for rigorous verification due to computational costs and the wrapping effect. Instead we compute $\mathbf{\Phi}(\mathbf{u}_j(\nu);x)$ using rigorous computation for a discrete set of Chebyshev nodes $\{\nu_m\}$ in $[\nu_0,\nu_{M+1}]$. This discretization, determined by the desired analytic interpolation error bound, is generally different and finer than the one used in the previous step. By employing rigorous Chebyshev interpolation (see Theorem~\ref{theorem:analytic_interpolation}), we obtain a polynomial representation of $\mathbf{\Phi}(\mathbf{u}_j(\nu);x)$ with a guaranteed error bound. Similarly, We derive rigorous Chebyshev representations for the derivatives $\frac{d\mathbf{\Phi}}{d\phi}$ and $\frac{d^2\mathbf{\Phi}}{d^2\phi}$, as well as for the stable and unstable solution manifold maps $\mathbf{u}_R(\theta_i;x_N,\theta_r)$ and $\mathbf{u}_L(\theta;0)$ and their derivatives.   

In practice, we make a change of coordinates in the parameter $\nu$. The stable manifold map $\mathbf{u}_R(\theta_i;x_N,\theta_r)$ is initialized with values corresponding to the fixed point of \eqref{eq:the_ode}, along with eigenvectors of the Jacobian evaluated at the fixed point. These eigenvectors take the form $\frac{1}{4\nu}(1,-2(1\pm i\sqrt{4\nu-1})$. However, near $\nu = \frac14$, using analytic interpolation to approximate these eigenvectors as polynomials in $\nu$ is not feasible because $\nu = \frac14$ corresponds to a branch point. To address this issue, we instead utilize $\alpha$, where $\alpha^2 = 4\nu-1$, enabling analytic interpolation. Despite this, technical coding issues persist.

Additionally, we adopt the parameter $\nu = \frac{1-\mu}{\mu^2}$ for the unstable manifold map $\mathbf{u}_L(\theta,0)$. This simplifies the eigenvector used in initializing $\mathbf{u}_L(\theta,0)$ to $(1,\mu)$. We use $\mu$ as the parameter in the finite region as well. During the Newton--Kantorovich argument, we must convert polynomials in $\alpha$ to polynomials in $\mu$. This is accomplished via analytic interpolation as in Theorem~\ref{theorem:analytic_interpolation}. To ensure that $0<\mu<1$, so that Lemma \ref{lemma:re_mu} holds, we enforce
\eq{
-\frac{\mu_L+\mu_R}{\mu_R-\mu_L}< \mu < \frac{2-(\mu_L+\mu_R)}{\mu_R-\mu_L},
}{\label{eq:mu_ineq}}
where $[\mu_L,\mu_R]$ is the interval of parameters treated simultaneously.

Step 3 begins the Newton--Kantorovich argument, carried out uniformly across $\nu$. With quantities like $\mathbf{A}$ and $\mathbf{A}^\dagger$ approximately calculated as Chebyshev series as in the previous step, and $\mathbf{F}$, in \eqref{def:boldF}, rigorously calculated with an error in the previous step, it becomes necessary to compute $Y_0$, $Z_0$, $Z_1$, $Z_2$ as in Theorem~\ref{thm:NewtK}. For instance, 
\[
Y_0(\nu)=\|\mathbf{A}(\nu)\mathbf{F}(\nu,y_{\rm ap}(\nu))\|_Y
\]
is an explicit but intricate polynomial in $\nu$, accurate within some previously determined error bound. By dividing the $\nu$ interval into smaller segments and computing $Y_0(\nu)$ on each subinterval using interval arithmetic, we obtain 
\[
\overline{Y_0} = \sup_\nu \|\mathbf{A}(\nu)\mathbf{F}(\nu,y_{\rm ap}(\nu))\|_Y, 
\]
and similarly for $Z_0$, $Z_1$, and $Z_2$. These uniform parameters
$\overline{Y_0}$, $\overline{Z_0}$, $\overline{Z_1}$, $\overline{Z_2}$, upon application of the Newton--Kantorovich theorem, enable a uniform estimate of the solution across $\nu$ at spatial nodes. 

Step 4 is to take the values, with error bounds, of $\phi$, $\psi$, $v$, $w$ at spatial nodes, and extend them to intervals using the solution operator $\mathbf{\Phi}(\mathbf{u}_j(\nu);x)$, and the stable and unstable solution manifold operators. Since the bounds derived in this step are uniform across the range of $\nu$ under consideration, the remaining process mirrors the single parameter case. Zeroes can be located by identifying intervals in which $v$ changes its sign, and uniqueness of such a zero is guaranteed by non-vanishing of $w=v'$ within the same interval.  

\section{Concluding remarks}

In summary, we have successfully established the asymptotic, nonlinear, and orbital stability of traveling front solutions, a general class of nonlinear diffusive-dispersive equations of Burgers type, provided that $-\frac{d^2}{dx^2}+\frac12\phi_x(x)$ has precisely one negative eigenvalue, where $\phi$ denotes the front. Notably, our result transcends the monotonicity of the profile and does not rely on the proximity of the initial condition to the front. Furthermore, the dispersion operator plays a rather negligible role in our proof, thus enabling its extension to many other dispersion operators. 

While our proof heavily relies on the nonlinearity and the dissipation operator, it holds the promise for extensions to treat a broad class of equations. Of particular interest is 
\begin{equation}\label{eqn:fKdVB}
u_t+uu_x=\mathcal{D}^\alpha_-u_x+\nu u_{xxx},
\quad\text{where}\quad
(\mathcal{D}^\alpha_-v)(x)=\frac{1}{\Gamma(1-\alpha)} \int^x_{-\infty}\frac{v_x(y)~dy}{(x-y)^\alpha},
\end{equation}
$0<\alpha<1$, represents the left-sided Caputo fractional derivative \cite{KS2006, Uchaikin1, Uchaikin2}. Achleitner, Cuesta and Hittmeir \cite{ACH2014} established the existence of traveling front solutions of \eqref{eqn:fKdVB}, whose profile satisfies, say, $\lim_{x\to\mp\infty}\phi(x)=\pm1$, demonstrating the monotonicity of the profile for $|\nu|\ll1$. See also \cite{AHS2011} for $\nu=0$. Following a similar line of argument as in Section~\ref{sec:analysis}, we anticipate that we establish the asymptotic, nonlinear, and orbital stability, provided that $-\frac{d}{dx}\mathcal{D}^\alpha_-+\frac12\phi_x(x)$ has precisely one negative eigenvalue, where $\phi$ denotes the traveling front solution of \eqref{eqn:fKdVB}. It will be interesting to verify the stability condition for monotone profiles.

In the case of the KdVB equation (see \eqref{eqn:KdVB}), we have established that the stability condition holds for $|\nu|$ within the range $[0,\frac14]$ through rigorous analysis and, additionally, for $|\nu|$ in the range $\ourint$ through CAP techniques. Our current focus is directed toward extending this result to cover the entire targeted range $0.25\leq |\nu|\lesssim 4.09$. 

However, as we progress towards higher values of $|\nu|$, we anticipate that the effective width of the profile will increase, leading to computational time and memory issues. As $|\nu|\rightarrow \frac14$, we encounter non-analytic behavior in the eigenvalues of the linearization about the fixed point $(\phi,\phi')=(-1,0)$. 

While our computations suggest a second negative eigenvalue for $-\frac{d^2}{dx^2}+\frac12\phi_x(x)$ for $|\nu|\leq 4.1$, numerical evidence \cite{CanosaGazdag} supports spectral stability for all $\nu\in\mathbb{R}$. It will be interesting to employ CAP techniques to rigorously establish stability for $|\nu|\gtrsim 4.1$, particularly, in the KdV limit as $|\nu|\to\infty$.

The authors would like to thank the reviewer for a very thorough review.

 \section*{Data Availability:}
 The codes used in the validated numerical proof are available for download from \\
 {\em {https://github.com/nonlinear-waves/stablab\_matlab/tree/master/KdVB. }}

 \section*{Ethics Declarations}
 \subsection*{Conflict of Interest}
 On behalf of all of the authors the corresponding author states that there is no conflict of interest. 
 
\section*{Funding}
JCB acknowledges support under a Simons travel grant. 
JCB and VMH acknowledge support under DMS-2407358

\bibliographystyle{amsplain}
\bibliography{Front}

\end{document}